\documentclass[11pt,reqno]{amsart}
\usepackage{mathrsfs}
\usepackage{amsmath}
\usepackage{amssymb, amsmath}
\usepackage{graphicx}

\newcommand{\newcom}{\newcommand}

\newcom{\al}{\alpha}
\newcom{\be}{\beta}
\newcom{\eps}{\epsilon}
\newcom{\ga}{\gamma}
\newcom{\Ga}{\Gamma}
\newcom{\ka}{\kappa}
\newcom{\Lam}{\Lambda}
\newcom{\lam}{\lambda}
\newcom{\la}{\lambda}
\newcom{\Om}{\Omega}
\newcom{\om}{\omega}
\newcom{\Si}{\Sigma}
\newcom{\si}{\sigma}
\newcom{\tht}{\theta}
\newcom{\dtri}{\nabla}
\newcom{\tri}{\triangle}
\newcom{\oo}{\infty}
\newcom{\vphi}{\varphi}
\newcom{\cA}{{\mathcal A}}
\newcom{\cB}{{\mathcal B}}
\newcom{\cC}{{\mathcal C}}
\newcom{\cD}{{\mathcal D}}
\newcom{\cE}{{\mathcal E}}
\newcom{\cF}{{\mathcal F}}
\newcom{\cG}{{\mathcal G}}
\newcom{\cL}{{\mathcal L}}
\newcom{\cM}{{\mathcal M}}
\newcom{\cP}{{\mathcal P}}
\newcom{\cS}{{\mathcal S}}
\newcom{\cQ}{{\mathcal Q}}
\newcom{\caly}{{\mathcal Y}}
\newcom{\calZ}{{\mathcal Z}}
\newcom{\bfz}{{\bf Z}}
\newcom{\R}{\Bbb R}
\newcom{\N}{\Bbb N}
\newcom{\Z}{\Bbb Z}
\newcom{\C}{\Bbb C}
\newcom{\E}{\Bbb E}

\def\dv{\mbox{div}}
\def\bnu{\bar{\nu}}
\def\ba{\bar{a}}
\def\bd{\bar{d}}
\newcom{\hn}{{\bf H}^n}
\newcom{\hnn}{{\mathbf H}^{n'}}
\newcom{\ulzs}{u^\lam_{z,s}}
\newcom{\Hl}{{{\cal  H}_\lam}}
\newcom{\fal}{F_{\al, \lam}}
\newcom{\Dh}{\Delta_{{\mathbf H}^n}}
\newcom{\fgl}{F_{\g, \lam}}
\newcom{\f}{\frac}
\newcom{\di}{\displaystyle\int}
\newcom{\ds}{\displaystyle\sum}
\newcom{\dl}{\displaystyle\lim}
\newcom{\ov}{\overline}
\newcom{\sset}{\subset}
\newcom{\wt}{\widetilde}
\newcom{\pa}{\partial}
\newcom{\p}{\partial}
\newcom\na{\nabla}
\newcom{\co}{\cdot}
\newcom{\suml}{\sum\limits}
\newcom{\supl}{\sup\limits}
\newcom{\intl}{\int\limits}
\newcom{\infl}{\inf\limits}
\newcom{\disp}{\displaystyle}
\newcom{\non}{\nonumber}
\newcom{\no}{\noindent}
\newcom{\QED}{$\square$}

\def\ef{\hphantom{MM}\hfill\llap{$\square$}\goodbreak}

\def\eqdefa{\buildrel\hbox{\footnotesize def}\over {=\!\!=}}

\newtheorem{athm}{\bf \t}[section]
\newenvironment{thm} [1] {\def\t{#1}\begin{athm} \bf \rm} {\end{athm}}
\newcom{\bthm}{\begin{thm}}\newcom{\ethm}{\end{thm}}

\newcom{\beq}{\begin{equation}}
\newcom{\eeq}{\end{equation}}
\newcom{\ben}{\begin{eqnarray}}
\newcom{\een}{\end{eqnarray}}
\newcom{\beno}{\begin{eqnarray*}}
\newcom{\eeno}{\end{eqnarray*}}
\newcom{\bali}{\begin{aligned}}
\newcom{\eali}{\end{aligned}}

\numberwithin{equation}{section}

\begin{document}

\title[Global well-posedness of the compressible Navier-Stokes equations]
{Global well-posedness for the compressible Navier-Stokes equations
with the highly oscillating initial velocity}

\author{Qionglei Chen}
\address{Institute of Applied Physics and Computational Mathematics,P.O. Box 8009, Beijing 100088, P. R. China}
\email{chen\_qionglei@iapcm.ac.cn}

\author{Changxing Miao}
\address{Institute of Applied Physics and Computational Mathematics,P.O. Box 8009, Beijing 100088, P. R. China}
\email{miao\_changxing@iapcm.ac.cn}

\author{Zhifei Zhang}
\address{School of Mathematical Sciences, Peking University, 100871, P. R. China}
\email{zfzhang@math.pku.edu.cn}

%\thanks{}

\date{2,August, 2009}



\begin{abstract}
Cannone \cite{Cannone} proved the global well-posedness of
the incompressible Navier-Stokes equations for a class of highly oscillating data.
In this paper, we prove the global well-posedness  for the compressible
Navier-Stokes equations in the critical functional framework with the initial data close to a stable equilibrium.
Especially, this result allows us to construct global solutions for the highly oscillating initial velocity.
The proof relies on a new estimate for the hyperbolic/parabolic system with convection terms.

\end{abstract}

\maketitle

\section{Introduction}
We consider the compressible Navier-Stokes equations
in $\mathbb{R}^+\times \mathbb{R}^n(n=2,3)$
\begin{equation}\label{equ:cNS}
\left\{
\begin{array}{ll}
\p_t\rho+\textrm{div}(\rho u)=0,\\
\p_t(\rho u)+\textrm{div}(\rho u\otimes u)-\mu\Delta u-(\lambda+\mu)\na\textrm{div}u+\na P(\rho)=0, \\
(\rho,u)|_{t=0}=(\rho_0,u_0).
\end{array}
\right.
\end{equation}
Here $\rho(t,x)$ and $u(t,x)$ denote the density and velocity of the
fluid respectively. The pressure $P$ is a suitably smooth function of $\rho$,
and the Lam\'{e} coefficients $\mu$ and $\lambda$ satisfy
\begin{align}\label{assu:coeff}
\mu>0\quad \textrm{and} \quad \lambda+2\mu>0,
\end{align}
which ensures that the operator $\mu\Delta
+(\lambda+\mu)\na\textrm{div}$ is elliptic.

The local existence and uniqueness of smooth solutions for the
system (\ref{equ:cNS}) were proved by Nash \cite{Nash} for smooth
initial data without vacuum. Later on, Matsumura and Nishida\cite{Mat-Nish}
proved the global well-posedness for smooth data close to
equilibrium. Recently, Lions\cite{Lions} proved the global existence of weak solutions
for large initial data. However, the question of uniqueness of weak solutions remains open,
even in the two dimensional case.

A natural way of dealing with the problem of uniqueness is to
find a functional setting as large as possible in which existence and uniqueness hold.
For the incompressible Navier-Stokes equations
\begin{align} \left\{
\begin{array}{ll}
\p_tu-\nu\Delta u+u\cdot\na u+\na p=0, \\
\dv u=0,\\
u(x,0)=u_0,
\end{array}
\right.\label{equ:NS}
\end{align}
such an approach was initiated by Fujita and Kato\cite{Fuj-Kat}. They proved the global existence and uniqueness
of strong solution for small initial data in the homogeneous Sobolev space $\dot H^{\f n2-1}$.
Let's point out that the space $\dot H^{\f n2-1}$ is a critical space.
We give the precise meaning of critical spaces. It is well-known that if
$(u,p)$ is a solution of (\ref{equ:NS}), then
\begin{align}\label{scaling}
(u_\lambda(t,x),p_\lambda(t,x))\eqdefa (\lambda u(\lambda^2 t,\lambda x),\lambda^2 p(\lambda^2 t,\lambda x))
\end{align}
is also a solution of (\ref{equ:NS}). A functional space is called
critical if the corresponding norm is invariant under the scaling of
(\ref{scaling}). For the compressible Navier-Stokes equations, one
can check that if $(\rho,u)$ is a solution of (\ref{equ:cNS}), then
\begin{align*}
(\rho_\lambda(t,x), u_\lambda(t,x))\eqdefa
\bigl(\rho(\lambda^2 t,\lambda x),\lambda u(\lambda^2 t,\lambda x)\bigr),
\end{align*}
is also a solution of (\ref{equ:cNS}) provided the pressure law has been changed into $\lambda^2 P$.
This motivates the following definition.
\bthm{Definition}
A functional space is called critical
if the associated norm is invariant under the transformation
$
(\rho,u)\longrightarrow (\rho_\lambda,u_\lambda)
$(up to a constant independent of $\lambda$).
\ethm

Consequently, a natural candidate is the homogenous Sobolev space $\dot
H^{\f {n} 2}\times \bigl(\dot H^{\f {n} 2-1}\bigr)^n$. However, $\dot
H^{\f {n} 2}$ is not included in $L^\infty$ such that we cannot expect to obtain a
$L^\infty$ control of the density when $\rho_0-\bar \rho\in \dot H^{\f {n} 2}$.
Instead, one can choose the initial data $(\rho_0,u_0)$ such that for some
$\bar{\rho}$,
\begin{align*}
(\rho_0-\bar{\rho},u_0)\in\dot B^{\f n 2}_{2,1}\times \bigl(\dot
B^{\f n 2-1}_{2,1}\bigr)^n, \end{align*} since $\dot B^{\f n
2}_{2,1}$ is continuously embedded in $L^\infty$. In a seminal
paper, Danchin\cite{Danchin-inven} proved the global well-posedness
of (\ref{equ:cNS}) in the critical Besov space for small initial
data $u_0$ and $\rho_0$ close to $\bar{\rho}$. More precisely,

\bthm{Theorem}\label{thm:Danchin}
Let $\bar{\rho}>0$  be such that $P'(\bar{\rho})>0$. There exist two positive
constants $c$ and $M$ such that for all $(\rho_0,u_0)$ with $\rho_0-\bar{\rho}\in \dot{\cB}^{\f n2-1,\f n2}$,
$u_0\in \dot B^{\f n2-1}_{2,1}$ and
$$\|\rho_0-\bar{\rho}\|_{\dot{\cB}^{\f n2-1,\f n2}}+\|u_0\|_{\dot B^{\f n2-1}_{2,1}}
\le c$$
then the system \eqref{equ:cNS} has a global unique solution $(\rho-\bar{\rho},u)$ such that
\begin{align*}
\|(\rho_0-\bar{\rho},u)\|_{F^{\f n2}}\le M
\big(\|\rho_0-\bar{\rho}\|_{\dot{\cB}^{\f n2-1,\f n2}}+\|u_0\|_{\dot B^{\f n2-1}_{2,1}}\big),
\end{align*}
where
\begin{align*} F^{s}\eqdefa L^1(\R^+; \dot{\cB}^{s+1,s})
\cap C(\R^+; \dot{\cB}^{s-1,s})
\times\Big(L^1(\R^+; \dot  B^{s+1}_{2,1})\cap C(\R^+; \dot  B^{s-1
}_{2,1})\Big)^n.\end{align*}
Here and what follows, we refer to Section 2.2 for the definition of Besov spaces.
\ethm

On the other hand, Fujita-Kato result was generalized to the critical Lebesgue space $L^n$ by
Weissler \cite{Weissler}(see also \cite{Kato}), and to a class of function spaces  such as $\dot{B}^0_{n, \infty}$ by Cannone, Meyer and Planchon
\cite{Cannon-Meyer-Planchon} which allows to construct self-similar solutions.
 Recently, Cannone \cite{Cannone} generalized it to Besov spaces of negative index of regularity.  More precisely,
he showed that if the initial data satisfies
$$
\|u_0\|_{\dot B^{-1+\f np}_{p,\infty}}\le c,\quad p>n
$$
for some small constant $c$, then the Navier-Stokes equations (\ref{equ:NS}) is globally well-posed.
Let us emphasize that this result allows us to construct global solutions for highly oscillating initial data
which may have a large norm in $\dot H^{\f n2-1}$ or $L^n$. A typical example in three dimensions is
\begin{align*}
u_0(x)=\sin\bigl(\f {x_3} {\varepsilon}\bigr)(-\p_2\phi(x), \p_1\phi(x),0)
\end{align*}
where $\phi\in \cS(\R^3)$ and $\varepsilon>0$ is small enough.
We refer to \cite{Chemin-Gal,Chemin-Zhang, Chemin-Gal-Paicu, Paicu-Zhang} for some relevant results.
An important question is then to prove a theorem of this type for the compressible Navier-Stokes equations.

To solve this problem, the key point is to generalize the well-posed result in Theorem \ref{thm:Danchin} to
the critical functional spaces with negative index of regularity(the corresponding norm of the highly oscillating function is small).
A typical space is $\dot B^{\f n p}_{p,1}\times \bigl(\dot B^{\f n p-1}_{p,1}\bigr)^n
$ for some $p>n$. Let's recall the proof of Theorem \ref{thm:Danchin},
where the energy method is applied to obtain the smoothing effect and damping effect of the system.
This method works well for the well-posed problem in the $L^2$ framework.
However, it doesn't work for the well-posed problem in the $L^p$ framework for $p\neq 2$,
especially in the case when the coupling effect of the system is important for the problem.
In this paper, we will develop a new method based on Green's matrix of the linearized system
to study the well-posed problem in general Besov spaces.

In order to present our method, let us look at the linearized system of (\ref{equ:cNS})
with the convection terms:
\begin{equation}\label{equ:linerized system-intro}
\left\{
\begin{array}{ll}
\p_ta+\Lambda d=-v\cdot\na a+F,\\
\p_td-\bar{\nu}\Delta d-\Lambda a=-v\cdot\na d+G,\\
(a,d)|_{t=0}=(a_0,d_0),
\end{array}
\right.
\end{equation}where $\Lambda=(-\Delta)^{\f12}$.
Let $\cG(x,t)$ be the Green matrix of the linear system
\begin{equation}
\left\{
\begin{array}{ll}
\p_ta+\Lambda d=0,\\
\p_td-\bar{\nu}\Delta d-\Lambda a=0.\nonumber
\end{array}
\right.
\end{equation}
Then the solution $(a,d)$ of (\ref{equ:linerized system-intro}) can be expressed as
\begin{align}\label{equ:integrate}
\left(\begin{array}{ll} \!a \!\!\vspace{.15cm}\\ \!d\!\end{array}\right)=
\mathcal{G}(x,t)\ast\left(\begin{array}{ll} \!a_0 \!\vspace{.15cm}\\ \!d_0 \!\end{array}\right)
+\int_0^t\mathcal{G}(x,t-\tau)\ast\left(\begin{array}{ll}\!F-v\cdot\na a\!\vspace{.15cm}\\
\!G-v\cdot\na d\!\end{array}\right)d\tau.
\end{align}
To obtain the estimates of $(a,d)$ in the Besov spaces, it is natural to study the action of Green's matrix
$\mathcal{G}(x,t)$ on distributions whose Fourier transform
is supported in a ring. As we know, the semi-group of the heat
equation shows for any couple $(t,\lambda)$ of positive real numbers,
\begin{align*}{\rm supp}\, \hat{f}\subset\lambda \cC\Longrightarrow
 \|e^{t\Delta}f\|_{L^p}\le Ce^{-ct\lambda^2}\|f\|_{L^p},\quad 1\le p\le\infty,
\end{align*}
where $\cC$ is a ring, which plays an important role in a new proof given by Chemin\cite{Chemin-Lecture}
on the Cannone, Meyer and Planchon result.
To obtain similar properties for $\cG(x,t)$, we need to investigate the precise behaviour of $\widehat{\cG}(\xi,t)$.
Roughly speaking, for $|\xi|\lesssim 1$, $\widehat{\mathcal{G}}(\xi,t)$ behaves like the heat kernel:
$$|D^\alpha_{\xi}\widehat{\mathcal{G}}(\xi,t)|\le
Ce^{-\vartheta|\xi|^2t}(1+|\xi|)^{|\alpha|}(1+t)^{|\alpha|},$$
which allows us to obtain
\begin{align}\label{equ:Green-Lowest}
\|\mathcal{G}\ast f\|_{L^2}\le Ce^{-c
t\lambda^2}\|f\|_{L^2},
\end{align}
if ${\rm supp}\, \hat{f}\subset\lambda \cC$ for $\lambda\lesssim 1$.
For $|\xi|\gg 1$, $\widehat{\mathcal{G}}(\xi,t)$ has the following precise expansion:
\begin{align*}
\widehat{\mathcal{G}}(\xi,t)=&e^{-\bnu^{-1}t}\left[\begin{array}{ll}1\quad0\\0\quad0\end{array}\right]
+e^{-\bnu|\xi|^2t}\left[\begin{array}{ll}0\quad0\\0\quad1\end{array}\right]\nonumber\\
&+\widehat{\mathcal {G}}^1(\xi,t)\left[\begin{array}{ll}0\quad 1\\1\quad 0\end{array}\right]+
\widehat{\mathcal {G}}^2(\xi,t)\left[\begin{array}{ll}1\quad0\\0\quad1\end{array}\right],
\end{align*}
where $\widehat{\mathcal {G}}^1$ and $\widehat{\mathcal {G}}^2$ satisfy the estimates
\begin{align*}
&|\pa^\al_\xi\widehat{\mathcal {G}}^1|
\le C|\xi|^{-|\al|-1}\big(e^{-\frac{\bnu^{-1}t}{2}}+e^{-\vartheta|\xi|^2t}\big),\\
&|\pa^\al_\xi \widehat{\mathcal {G}}^2|
\le C|\xi|^{-|\al|-2}\big(e^{-\frac{\bnu^{-1}t}{2}}+e^{-\vartheta|\xi|^2t}\big),
\end{align*}
which allow us to obtain, for any $1\le p\le \infty$,
\begin{align}\label{equ:Green-Highest}
\|\mathcal{G}^1\ast f\|_{L^p}\le C\lambda^{-1}e^{-c
t}\|f\|_{L^p},\quad
\|\mathcal{G}^2\ast f\|_{L^p}\le C\lambda^{-2}e^{-c t}\|f\|_{L^p},
\end{align}
if ${\rm supp}\, \hat{f}\subset\lambda \cC$ for $\lambda\gg 1$.
Since $\cG$ has different behaviour for low frequency and high frequency, it is natural to
work with the hybrid-Besov spaces(see Definition \ref{Def:hybridBes}). Indeed, from
(\ref{equ:integrate}), (\ref{equ:Green-Lowest}) and (\ref{equ:Green-Highest}),
we can obtain the following estimate in the hybrid-Besov spaces: for any $1\le r\le\infty$,
\begin{align}\label{equ:hype-parabolic-est}
&\|a\|_{\widetilde{L}^r_T(\dot \cB_{2,p}^{s_p-1+\f 2r,s})}
+\|d\|_{\widetilde{L}^r_T(\dot \cB_{2,p}^{s_p-1+\f 2r,s-1+\f 2r})}\nonumber\\
&\quad\le C\Bigl[\|a_0\|_{\dot\cB^{s_p-1,s}_{2,p}}+\|d_0\|_{\dot\cB^{s_p-1,s-1}_{2,p}}+\|F\|_{L^1_T(\dot\cB_{2,p}^{s_p-1,s})}
+\|G\|_{L^1_T(\dot\cB_{2,p}^{s_p-1,s-1})}
\nonumber\\
&\qquad\qquad\qquad\qquad+\|v\cdot\na a\|_{L^1_T(\dot\cB_{2,p}^{s_p-1,s})}+\|v\cdot\na d\|_{L^1_T(\dot\cB_{2,p}^{s_p-1,s-1})}\Bigr],
\end{align}
where $s_p\eqdefa s-\f np+\f n2$. Let us point out that if we take $s=\f np$ for $p>n$,
the regularity index $s-1$ for the high frequency part of $d_0$ is negative,
which is enough to ensure that the norm of $d_0$ in the hybrid-Besov space $\dot\cB^{s_p-1,s-1}_{2,p}$ is small
if $d_0$ is highly oscillating, see Proposition \ref{Prop:oscillate}.

Unfortunately, the inequality (\ref{equ:hype-parabolic-est}) cannot be applied directly. In fact,
if  the convection term $v\cdot\na a$ is viewed as a perturbation term,
one derivative loss about the function $a$ will appear no matter how smooth is $v$.
To overcome this difficulty, we will work in the Lagrangian coordinate such that
the convection terms in (\ref{equ:linerized system-intro}) will disappear modulus some commutators.
Here the idea is partially motivated by the recent work \cite{Hmidi-Keraani-Advance,Hmidi-Keraani-ARMA}
for the incompressible flows.

We denote
\begin{align*}
\cE^{s}\eqdefa \bigg\{(a,d)\in \Big(&L^1(\R^+; \dot
\cB^{s_p+1,s}_{2,p})
\cap C(\R^+; \dot \cB^{s_p-1,s}_{2,p})\Big)\\
&\times\Big(L^1(\R^+; \dot \cB^{s_p+1,s+1}_{2,p})\cap C(\R^+; \dot \cB^{s_p-1,
s-1}_{2,p})\Big)^n\bigg\}
\end{align*}
with $s_p\eqdefa s-\f np+\f n2$. Our main result is stated as follows.

\bthm{Theorem}\label{thm:Global}
Let $\bar{\rho}$  be a positive constant such that $P'(\bar{\rho})>0$. There exist two positive
constants $c$ and $M$ such that for all $(\rho_0,u_0)$ with $\rho_0-\bar{\rho}\in \dot{\cB}^{\f n2-1,\f np}_{2,p}$,
$u_0\in \dot \cB^{\f n2-1,\f np-1}_{2,p}$ and
\begin{align}\label{equ:initial condition}
\|\rho_0-\bar{\rho}\|_{\dot \cB^{\f {n} 2-1,\f n p}_{2,p}}+
\|u_0\|_{\dot \cB^{\f {n} 2-1,\f n p-1}_{2,p}}\le\eta,
\end{align}
the following results hold

(1)\,{\bf Existence.}\, if $2\le p<2n,$ and $p\le\min(4,\f {2n}{n-2})$, the system (\ref{equ:cNS})
has a global solution $(\rho-\bar{\rho},u)\in \cE^{\f np}$ with
\begin{align*}
\|(\rho-\bar{\rho},u)\|_{\cE^{\f np}}\le
M\big(\|\rho_0-\bar{\rho}\|_{\dot \cB^{\f {n} 2-1,\f n p}_{2,p}}+
\|u_0\|_{\dot \cB^{\f {n} 2-1,\f n p-1}_{2,p}}\big).
\end{align*}

(2)\,{\bf Uniqueness.}\, if $2\le p\le n$, then the uniqueness holds in $\cE^{\f np}$.
\ethm

\bthm{Remark} Compared with Theorem \ref{thm:Danchin}, an important improvement of
Theorem \ref{thm:Global} is that it allows the regularity index for
the high frequency part of $u_0$ to be negative. Especially, this
allows us to obtain the global well-posedness of (\ref{equ:cNS}) for
the highly oscillating initial velocity $u_0$. For example,
$$
u_0(x)=\sin\bigl(\f {x_1} {\varepsilon}\bigr)\phi(x),\quad \phi(x)\in \cS(\R^n),
$$
which satisfies
$$
\|u_0\|_{\dot \cB^{\f n2-1,\f np-1}_{2,p}}\ll 1 \quad \textrm{for}\quad p>n
$$
if $\varepsilon$ is small enough, see Proposition \ref{Prop:oscillate}.
\ethm

\bthm{Remark}
In fact, Theorem \ref{thm:Global} also holds for general dimension.
But we need to restrict the dimension $n=2,3$ in order to obtain  global solutions for the
highly oscillating initial velocity.
\ethm

\bthm{Remark}
We believe that our method can be adapted to the other hyperbolic-parabolic models. This is the object of
our future work.
\ethm

The structure of this paper is as follows.\vspace{.1cm}

In Section 2, we recall some basic facts about Littlewood-Paley
theory and the functional spaces. In Section 3, we recall some
results concerning the classical parabolic regularizing effect.
Section 4 is devoted to the study of Green's matrix of the
linearized system. In Section 5, we establish the key estimates for
the the linearized system with convection terms. Section 6 is
devoted to the proof of Theorem \ref{thm:Global}. Finally, an
appendix is devoted to some estimates in the hybrid-Besov spaces.

\section{Littlewood-paley theory and the function spaces}
\subsection{Littlewood-Paley theory}
First of all, we introduce the Littlewood-Paley decomposition. Choose two
radial functions  $\varphi, \chi \in {\cS}(\mathbb{R}^n)$ supported in
${\cC}=\{\xi\in\mathbb{R}^n,\, \frac{3}{4}\le|\xi|\le\frac{8}{3}\}$,
${\cB}=\{\xi\in\mathbb{R}^n,\, |\xi|\le\frac{4}{3}\}$ respectively such
that
\begin{align*} \sum_{j\in\mathbb{Z}}\varphi(2^{-j}\xi)=1 \quad \textrm{for
all}\,\,\xi\neq 0.
\end{align*}
The frequency localization operators $\Delta_j$ and $S_j$ are defined by
\begin{align}
\Delta_jf=\varphi(2^{-j}D)f,\quad S_jf=\sum_{k\le
j-1}\Delta_kf=\chi(2^{-j}D)f,\quad\mbox{for}\,\, j\in \mathbb{Z}. \nonumber
\end{align}
With our choice of $\varphi$, it is easy to verify that
\begin{align}\label{orth}
\begin{aligned}
&\Delta_j\Delta_kf=0\quad \textrm{if}\quad|j-k|\ge 2\quad
\textrm{and}
\quad \\
&\Delta_j(S_{k-1}f\Delta_k f)=0\quad \textrm{if}\quad|j-k|\ge 5.
\end{aligned}
\end{align}

In the sequel, we will constantly use the Bony's decomposition from \cite{Bony}:
\begin{align}\label{Bonydecom}
fg=T_fg+T_gf+R(f,g), \end{align} with
$$T_fg=\sum_{j\in\mathbb{Z}}S_{j-1}f\Delta_jg,
\quad R(f,g)=\sum_{j\in\mathbb{Z}}\Delta_jf \widetilde{\Delta}_{j}g,
\quad \widetilde{\Delta}_{j}g=\sum_{|j'-j|\le1}\Delta_{j'}g.$$

Next, let us introduce some useful lemmas which will be repeatedly used throughout this paper.

\bthm{Lemma}\cite{Chemin-book}\label{Lem:Bernstein}
Let $1\le p\le q\le+\infty$. Then for any $\gamma\in(\mathbb{N}\cup\{0\})^n$, there exists a constant $C$
independent of $f$, $j$ such that
\begin{align*} &{\rm supp}\hat f\subseteq
\{|\xi|\le A_02^{j}\}\Rightarrow \|\partial^\gamma f\|_{L^q}\le
C2^{j{|\gamma|}+j n(\frac{1}{p}-\frac{1}{q})}\|f\|_{L^p},
\\
&{\rm supp}\hat f\subseteq \{A_12^{j}\le|\xi|\le
A_22^{j}\}\Rightarrow \|f\|_{L^p}\le
C2^{-j|\gamma|}\sup_{|\beta|=|\gamma|}\|\partial^\beta f\|_{L^p}.
\end{align*}
\ethm

Let a vector field $v(t,x)\in L^1_{\textrm{lo}c}(\R;\textrm{Lip}(\R^n))$. We define $\psi$ and $\psi_j$ as follows
\begin{align*}
&&\f {d} {dt}\psi(t,x)=v(t,\psi(t,x)), \quad \psi(0,x)=x,\\
&&\f {d} {dt}\psi_j(t,x)=S_{j-1}v(t,\psi_j(t,x)), \quad \psi_j(0,x)=x,\quad j\in \Z.
\end{align*}
In what follows, we denote $V(t)\eqdefa \int_0^t\|\na v(\tau)\|_{L^\infty}d\tau$.
Now we are position to state the following lemmas.

\bthm{Lemma}\label{Lem:com-est} Let $1\le p<\infty$. Then there hold
\begin{align*}
&\|f\circ\psi(t,x)\|_{L^p}
\le e^{V(t)}\|f\|_{L^p},\quad \|f\circ\psi(t,x)\|_{L^\infty}\le \|f\|_{L^\infty}.
\end{align*}
\ethm

\noindent{\bf Proof.} Lemma \ref{Lem:com-est} can be deduced from the fact
$$
\p_t \textrm{det}(\na \psi)=\dv\, v \textrm{det}(\na \psi),
$$
and a change of variables. \ef

\bthm{Lemma}\label{Lem:com-localest}
Let $1\le p\le \infty$. Then for any $j, k\in\mathbb{Z}$, there hold
\begin{align*}
&\|\Delta_j(\Delta_kf\circ\psi_k)\|_{L^p}\le C2^{-(j-k)}e^{CV(t)}\|\Delta_kf\|_{L^p},\\
&\|S_j(\Delta_kf\circ\psi_k)\|_{L^p}\le Ce^{CV(t)}\big(V(t)+2^{j-k}
\big)\|\Delta_kf\|_{L^p}.
\end{align*}
Here $C$ is a constant independent of $j,k$.
\ethm

\noindent {\bf Proof.}\, Let $\{\alpha_\ell\}_{\ell=1}^n$ be a smooth unity decomposition of ${\Bbb S}^{n-1}$ such that
$\xi_\ell' \neq 0$ on  $\textrm{supp}\, \alpha_\ell(\xi')$, here $\xi'=(\xi_1',\cdots,\xi'_n)\in {\Bbb S}^{n-1}$.
Then we write
$${\varphi}(\xi)=\sum_{\ell=1}^n i \xi_\ell\hat{\theta}_\ell(\xi),\,\,
\hat{\theta}_\ell(\xi)\eqdefa (i\xi_\ell)^{-1}{\alpha}_\ell(\xi/|\xi|){\varphi}(\xi),
$$
consequently,
\begin{align}\label{equ:unit-decomp}
\Delta_jf=2^{-j}\sum_{\ell=1}^n  \pa_\ell(2^{jn}\theta_\ell(2^j\,\cdot)*f)
\eqdefa 2^{-j}\sum_{\ell=1}^n  \pa_\ell(\bar{\Delta}_{j\ell}f),
\quad j\in\mathbb{Z},
\end{align}
from which, we get by integration by parts that
\begin{align*}
\Delta_j(\Delta_kf\circ\psi_k)
=&\sum_{\ell=1}^n2^{nj-j}\int_{\mathbb{R}^n}\theta_\ell(2^{j}(x-y))\pa_{\ell}((\Delta_kf)(\psi_k(y)))dy.
\end{align*}
Thus by Young's inequality, Lemma \ref{Lem:com-est}  and Lemma \ref{Lem:Bernstein}, we have
\begin{align}\label{eq:deltajesti}
\|\Delta_j(\Delta_kf\circ\psi_k)\|_{L^p}&\le C2^{-j}\|(\na \Delta_kf)(\psi_k)\|_{L^p}\|\na \psi_k\|_{L^\infty}\nonumber\\
&\le C2^{-j}e^{V(t)}\|\na\Delta_kf\|_{L^p}\|\na \psi_k\|_{L^\infty}\nonumber\\
&\le C2^{k-j}e^{V(t)}\|\Delta_kf\|_{L^p}\|\na \psi_k\|_{L^\infty}.
\end{align}

On the other hand, thanks to (\ref{equ:unit-decomp}), we have
\begin{align*}
S_j(\Delta_kf\circ\psi_k)=2^{-k}\sum_{\ell=1}^nS_j(\pa_\ell(\bar{\Delta}_{k\ell}f)\circ\psi_k).
\end{align*}
Set $g(x)\eqdefa(\cF^{-1}\chi)(x)$.  Making a change of variables: $y\mapsto \psi^{-1}_k(y)$,
and then integrating by parts yields
\begin{align*}
S_j(\Delta_kf\circ\psi_k)=&2^{nj-k}\sum_{\ell=1}^n\int_{\mathbb{R}^n}g(2^j(x-y))
\big(\pa_\ell(\bar{\Delta}_{k\ell}f) (\psi_k(y)\big)dy\\
=&2^{nj-k}\sum_{\ell=1}^n\int_{\mathbb{R}^n}g(2^j(x-\psi^{-1}_k(y)))
\pa_\ell(\bar{\Delta}_{k\ell}f)(y) \det(\na \psi^{-1}_k) dy\\
=&2^{nj-k}\sum_{\ell=1}^n\int_{\mathbb{R}^n}\Big\{2^j\pa_l\psi^{-1}_k(\pa_{y_\ell}g)(2^j(x-\psi^{-1}_k(y)))
\det(\na
\psi^{-1}_k)\nonumber\\&\qquad\qquad+g(2^j(x-\psi^{-1}_k(y)))
\pa_\ell\det(\na \psi^{-1}_k)\Big\}\bar{\Delta}_{k\ell}f(y)dy.
\end{align*}
Thus, we get by Young's inequality that
\begin{align}\label{eq:Sjesti}
\|S_j(\Delta_kf\circ\psi_k)\|_{L^p}\le &C2^{j-k}\|\na
\psi_k^{-1}\|^{n+1}_{L^\infty}\|\Delta_kf\|_{L^p}\nonumber\\&
\qquad+C2^{-k}\|\na^2\psi_k^{-1}\|_{L^\infty}\|\na
\psi_k^{-1}\|^{n-1}_{L^\infty}\|\Delta_kf\|_{L^p}.
\end{align}
Now, Lemma \ref{Lem:com-localest} follows from Proposition 3.1 in \cite{Danchin-JHDE}, \eqref{eq:deltajesti} and  \eqref{eq:Sjesti}.
\ef

\bthm{Lemma}\cite{Cannon-Miao-Wu}\label{Lem:com-commutator}
 Let $1\le p\le \infty$. Then for any $j\in \Z$, there hold
\begin{align*}
&\|\Delta(\Delta_j f\circ\psi_j)-(\Delta\Delta_jf)\circ\psi_j\|_{L^p}\le C2^{2j} e^{CV(t)}V(t)\|\Delta_jf\|_{L^p},\\
&\|\Lambda(\Delta_j f\circ\psi_j)-(\Lambda\Delta_jf)\circ\psi_j\|_{L^p}\le C2^je^{CV(t)}V(t)^{\f12}\|\Delta_jf\|_{L^p}.
\end{align*}
Here $C$ is a constant independent of $j,k$.
\ethm

\subsection{The hybrid-Besov space}

We denote the space ${\calZ'}(\mathbb{R}^n)$ by the dual space of
${\calZ}(\mathbb{R}^n)=\{f\in {\cS}(\mathbb{R}^n);\,D^\alpha \hat{f}(0)=0;
\forall\alpha\in\big(\mathbb{N}\cup 0\big)^n \,\mbox{multi-index}\}$.
Let us first recall the definition of general Besov space.

\bthm{Definition}\label{Def:Bes} Let $s\in\mathbb{R}$, $1\le p,
q\le+\infty$. The homogeneous Besov space $\dot{B}^{s}_{p,q}$ is
defined by
$$\dot{B}^{s}_{p,q}\eqdefa\big\{f\in {\calZ'}(\mathbb{R}^n):\,\|f\|_{\dot{B}^{s}_{p,q}}<+\infty\big\},$$
where
\begin{align*}
\|f\|_{\dot{B}^{s}_{p,q}}\eqdefa \Bigl\|2^{ks}
\|\Delta_kf(t)\|_{L^p}\Bigr\|_{\ell^q}.\end{align*}
\ethm

Now we introduce the hybrid-Besov space  we will work with in this
paper. Let $R_0>0$  be as in Proposition \ref{Prop:Green-Lpest}.
\bthm{Definition}\label{Def:hybridBes} Let $s$,
$\sigma\in\mathbb{R}$, $1\le p\le +\infty$. The hybrid-Besov space
$\dot{\cB}^{s,\sigma}_{2,p}$ is defined by
\begin{align*}
&\dot{\cB}^{s,\sigma}_{2,p}\eqdefa\big\{f\in{\calZ'}(\mathbb{R}^n):
\|f\|_{\dot{B}^{s,\sigma}_{2,p}}<+\infty\big\},
\end{align*}
where
$$\|f\|_{\dot{\cB}^{s,\sigma}_{2,p}}\eqdefa\sum_{2^{k}\le R_0}2^{ks}\|\Delta_k f\|_{L^2}
+\sum_{2^k>R_0}2^{k\sigma}\|\Delta_k f\|_{L^p}.
$$
We also denote $\dot \cB^{s,\sigma }$ by $\dot{\cB}^{s,\sigma}_{2,2}$.
\ethm

The norm of the space $\widetilde{L}^r_T(\dot{\cB}^{s,\sigma}_{2,p})$ is defined by
$$
\|f\|_{\widetilde{L}^r_T(\dot{\cB}^{s,\sigma}_{2,p})}
\eqdefa\sum_{2^k\le R_0}2^{ks}\|\Delta_k f\|_{L^r_TL^2}
+\sum_{2^k>R_0}2^{k\sigma}\|\Delta_k f\|_{L^r_TL^p}.
$$
It is easy to check that
$
 \widetilde{L}^1_T(\dot{\cB}^{s,\sigma}_{2,p})=L^1_T(\dot{\cB}^{s,\sigma}_{2,p}) \textrm{ and }
\widetilde{L}^r_T(\dot{\cB}^{s,\sigma}_{2,p})\subseteq {L}^r_T(\dot{\cB}^{s,\sigma}_{2,p})$ for $r>1$.

The following lemma is a direct consequence of the definition of Besov space
and Lemma \ref{Lem:Bernstein}.

\bthm{Lemma}\label{Lem:Besov-properties}
The following properties hold:\vspace{0.1cm}

(a)\,\,$\dot{\cB}^{s_2,\sigma}_{2,p}\subseteq \dot{\cB}^{s_1,\sigma}_{2,p}\quad\textrm{ if }s_1\ge s_2
\quad \textrm{and} \quad\dot{\cB}^{s,\sigma_2}_{2,p} \subseteq \dot{\cB}^{s,\sigma_1}_{2,p}\quad\textrm{ if }\sigma_1\le \sigma_2;
$

(b)\,\,interpolation: for
$s_1, s_2, \sigma_1, \sigma_2\in\mathbb{R}$ and $\theta\in[0,1]$,
we have $$\|f\|_{\dot \cB^{\theta s_1+(1-\theta)s_2,\, \theta \sigma_1+(1-\theta)\sigma_2}_{2,p}}\le
\|f\|^\theta_{\dot \cB^{s_1,\sigma_1}_{2,p}}\|f\|^{(1-\theta)}_{\dot \cB^{s_2,\sigma_2}_{2,p}};$$

(c)\,\,embedding: $\|f\|_{L^\infty}\le C\|f\|_{\dot{\cB}^{\f n2,\f n p}_{2,p}};$

(d)\,\,inclusion relation:
$\dot B^s_{2,1}\subseteq \dot \cB^{s,s-\f n2+\f n p}_{2,p}\subseteq \dot B^{s-\f n2+\f n p}_{p,1}$ for $p\ge 2$.
\ethm

\bthm{Lemma}\label{Lem:Product}\cite{Danchin-NDEA} Let $1\le p\le \infty$. Then there hold

(a)\; if $s_1, s_2\le \frac{n}{p}$ and $s_1+s_2>n\max
(0,\frac2p-1)$, then
\begin{align*}\|fg\|_{\dot B^{s_1+s_2-
\frac{n}{p}}_{p,1}} \le C\|f\|_{\dot B^{s_1}_{p,1}}\|g\|_{\dot B^{s_2}_{p,1}}.
\end{align*}

(b)\; if $s_1\le \f n p, s_2<\f n p$, and $s_1+s_2\ge n\max
(0,\frac2p-1)$, then
\begin{align*}\|fg\|_{\dot B^{s_1+s_2-
\frac{n}{p}}_{p,\infty}} \le C\|f\|_{\dot B^{s_1}_{p,1}}\|g\|_{\dot B^{s_2}_{p,\infty}}.
\end{align*}
\ethm

\bthm{Proposition}\label{Prop:oscillate} Let $\phi\in \cS(\R^n), p>n$.
If $\phi_\varepsilon(x)\eqdefa e^{i\f {x_1} \varepsilon}\phi(x)$, then for any $\varepsilon>0$,
\begin{align*}
\|\phi_\varepsilon\|_{\dot\cB^{\f n2-1,\f n p-1}_{2,p}}\le C\varepsilon^{1-\f np},
\end{align*}
here $C$ is a constant independent of $\varepsilon$.
\ethm

\noindent{\bf Proof.}\, Fix $j_0\in\N$ to be chosen later. By Lemma \ref{Lem:Bernstein}, we have
\begin{align*}
\sum_{j\ge j_0}2^{(\f np-1)j}\|\Delta_j \phi_\varepsilon\|_{L^p}\le C2^{(\f np-1)j_0}.
\end{align*}
Noting that $e^{i\f {x_1} \varepsilon}=(-i\varepsilon\p_1)^Ne^{i\f {x_1} \varepsilon}$ for any $N\in\N$,
we get by integration by parts that
\begin{align*}
\Delta_j\phi_\varepsilon(x)=(i\varepsilon)^N2^{jn}\int_{\R^n}
e^{i\f {y_1} \varepsilon}\p_{y_1}^N(h(2^j(x-y))\phi(y))dy,\quad h(x)\eqdefa (\cF^{-1}\varphi)(x),
\end{align*}
from which and Young's inequality, we infer that
\begin{align*}
\|\Delta_j\phi_\varepsilon\|_{L^q}\le C\varepsilon^N\max(2^{Nj},2^{(1-\f 1 q)nj}), \quad 1\le q\le \infty,
\end{align*}
which implies that by choosing $N>\frac{n}{p'}$
\begin{align*}
&\sum_{0 \le j<j_0}2^{(\f np-1)j}\|\Delta_j \phi_\varepsilon\|_{L^p}\le C\varepsilon^N 2^{(N-1+\f np)j_0},\\
&\sum_{j\le \ln R_0}2^{(\f n2-1)j}\|\Delta_j \phi_\varepsilon\|_{L^2}\le C\varepsilon^N.
\end{align*}
Taking $j_0$ such that $2^{j_0}\sim \varepsilon^{-1}$ gives
\begin{align*}
\|\phi_\varepsilon\|_{\dot\cB^{\f n2-1,\f n p-1}_{2,p}}\le C\varepsilon^{1-\f np}.
\end{align*}

This completes the proof of Proposition \ref{Prop:oscillate}.\ef

\section{Regularizing effect of the heat equation}

Let us recall the classical parabolic regularizing effect.
\bthm{Lemma}\cite{Chemin-Lecture}\label{Lem:parabolic}
Let $\phi$ be a smooth function supported in the annulus $\{\xi\in
\R^n: A_1\leq |\xi|\leq A_2\}$. Then there exist two positive
constants $c$ and $C$  such that for any $1\le p\le \infty$ and $\lambda>0$, we have
\beno
\|\phi(\lambda^{-1}D)e^{\nu t\Delta}f\|_{L^p}\le
Ce^{-ct\lambda^2}\|\phi(\lambda^{-1}D)f\|_{L^p}.\label{2.5} \eeno
\ethm

\bthm{Proposition}\label{Prop:Parabolic-est}
Let $s, \sigma\in\mathbb{R}$, and $p,r\in[1,\infty]$.
Assume that $u_0\in \dot \cB^{s,\sigma}_{2,p}, f\in L^1_T\dot \cB^{s,\sigma}_{2,p}$. Then the heat equation
\begin{align*}
\left\{
\begin{array}{ll}
\p_tu-\nu\Delta u=f,\\
u|_{t=0}=u_0,
\end{array}
\right.
\end{align*}
has a unique solution $u$ satisfying
$$\|u\|_{\widetilde{L}^r_T\dot \cB^{s+\f2 r,\sigma+\f2r}_{2,p}}
\le C\big(\|u_0\|_{\dot \cB^{s,\sigma}_{2,p}}+
\|f\|_{ L^1_T\dot \cB^{s,\sigma}_{2,p}}\big).$$
\ethm

\noindent{\bf Proof.}\,\,The solution $u$ can be written as
\begin{align*}
u(t)=e^{\nu t\Delta}u_0+\int_0^te^{\nu (t-s)\Delta}f(s)ds,
\end{align*}
from which and Lemma \ref{Lem:parabolic}, it follows that for any $1\le q\le \infty$,
\begin{align*}
\|\Delta_j u(t)\|_{L^q}\le Ce^{-ct2^{2j}}\|\Delta_ju_0\|_{L^q}
+C\int_0^te^{-c(t-s)2^{2j}}\|\Delta_jf(s)\|_{L^q}ds.
\end{align*}
Taking $L^r$ norm with respect to $t$ gives
\begin{align*}
\|\Delta_j u(t)\|_{L^r_TL^q}\le C2^{-\f 2 {r} j}\|\Delta_ju_0\|_{L^q}
+C2^{-\f 2 {r} j}\|\Delta_jf(s)\|_{L^1_TL^q}.
\end{align*}
Using $q=2$ for $2^j\le R_0$ and $q=p$ for $2^j> R_0$,
then Lemma follows from the definition of the hybrid-Besov space.\ef

To prove the uniqueness of the solution, we need the following two propositions.

\bthm{Proposition}\cite{Danchin-cpde01}\label{Prop:momentum} Let $p, q, r\in [1,\infty],s\in \R$.
Assume that $u_0\in \dot B^{s-1}_{p,q}, f\in L^1_T(\dot B^{s-1}_{p,q})$.
Let $u$ be a solution of the following equation
\begin{align*}
\p_t u-\bar{\mu}\Delta u+(\bar{\lambda}+\bar{\mu})\na \dv u=f,\quad u(0,x)=u_0(x),
\end{align*}
where $\bar\mu>0, \bar{\lambda}+\bar{\mu}>0$. Then for $t\in [0,T]$, there holds
\begin{align*}
\|u\|_{\widetilde{L}^r_t(\dot B^{s-1+2/q}_{p,q})} \le C\bigl(\|u_0\|_{\dot
B^{s-1}_{p,q}}+\|f(\tau)\|_{\widetilde{L}^1_t(\dot B^{s-1}_{p,q})}\bigr).
\end{align*}
\ethm

\bthm{Proposition}\cite{Danchin-NDEA}\label{Prop:transport}
Let $s\in (-n\min(\frac1p,\frac1{p'}), 1+\frac np)$, and $1\le
p,q\le+\infty$. Let $v$ be a vector field such that $\nabla v\in
L^1_T\dot{B}^{\frac{n}{p}}_{p,1}$. Assume that
$f_0\in \dot{B}^{s}_{p,q},$ $g\in L^1_T(\dot{B}^{s}_{p,q})$ and $f$
is a solution of the transport equation
\begin{align*}
\partial_t f+v\cdot \nabla f =g,\quad f(0,x)=f_0.
\end{align*}
Then  for $t\in[0,T]$, there holds
\begin{align*}
\|f\|_{\widetilde{L}^\infty_t(\dot{B}^{s}_{p,q})}\le e^{C\int_0^t\|\nabla
v(\tau)\|_{\dot{B}^{\frac{n}{p}}_{p,1}}d\tau}\Big(
\|f_0\|_{\dot{B}^{s}_{p,q}}+\int_0^t\|g(\tau)\|_{\dot{B}^{s}_{p,q}}d\tau\Big).
\end{align*}
\ethm

\section{Green's matrix of the linearized system}

We study the linearized system of the compressible Navier-Stokes system
\begin{align}\label{equ:linearized system}
\left\{
\begin{array}{ll}
\p_ta+\Lambda d=0,\\
\p_td-\bar{\nu}\Delta d-\Lambda a=0,\\
(a,d)|_{t=0}=(a_0,d_0).
\end{array}
\right.
\end{align}

We firstly give the explicit expression of Green's matrix for the system.
\bthm{Lemma}\label{Lem:Greenmatrix}
Let $\mathcal {G}$ be the Green matrix of the system (\ref{equ:linearized system}).
Then we have the following explicit expression for $\widehat{\mathcal {G}}$:
\begin{equation}\label{equ:Green matrix}
\widehat{\mathcal {G}}(\xi,t)=
\left[
\begin{array}{ll}
\frac{\lambda_+e^{\lambda_-t}-\lambda_-e^{\lambda_+t}}{\lambda_+-\lambda_-}
\qquad-\Big(\frac{e^{\lambda_+t}-e^{\lambda_-t}}{\lambda_+-\lambda_-}\Big)|\xi|\vspace{.3cm}\\
-\Big(\frac{e^{\lambda_+t}-e^{\lambda_-t}}{\lambda_+-\lambda_-}\Big)|\xi|\quad
\qquad \frac{\lambda_+e^{\lambda_+t}-\lambda_-e^{\lambda_-t}}{\lambda_+-\lambda_-}
\end{array}
\right],
\end{equation}where $$\lambda_\pm=-\frac12\bnu|\xi|^2\pm\frac12\sqrt{\bnu^2|\xi|^4-4|\xi|^2}.$$
\ethm

\noindent{\bf Proof.}\, We follow the proof of Lemma 3.1 in
\cite{Hoff-Zum}.  Firstly,
applying the operator $\Lambda$ to the second equation of (\ref{equ:linearized system}) gives
\begin{align*}
(\Lambda d)_t+\Delta a=\bnu\Delta (\Lambda d).
\end{align*}
Combining it with the first equation of (\ref{equ:linearized system}),
we get
\begin{align*}
a_{tt}=-(\Lambda d)_t=\Delta a+\bnu\Delta a_t.
\end{align*}
Taking Fourier transform to the above equation yields that
\begin{align}\label{equ:Fourierequ of a}
\left\{\begin{array}{ll}
\hat{a}_{tt}+\bnu|\xi|^2\hat{a}_{t}+|\xi|^2\hat{a}=0,\\
\hat{a}(\xi,0)=\hat{a}_0(\xi),\quad
\hat{a}_t(\xi,0)=-|\xi|\hat{d}_0(\xi).
\end{array}
\right.
\end{align}
It is easy to check that $\lambda_{\pm}$ are two roots of the corresponding
indicial equation of (\ref{equ:Fourierequ of a}).
Thus, we may assume that the solution of (\ref{equ:Fourierequ of a}) has the form
\begin{align*}
\hat{a}(\xi,t)=A(\xi)e^{\lambda_{-}(\xi)t}+B(\xi)e^{\lambda_{+}(\xi)t}.
\end{align*}
Using the initial conditions, we obtain
\begin{align*}
A=\f{\lambda_+\hat{a}_0+|\xi|\hat{d}_0}{\lambda_+-\lambda_{-}}\quad\hbox{and}\quad
B=\f{-|\xi|\hat{d}_0-\lambda_{-}\hat{a}_0}{\lambda_+-\lambda_{-}},
\end{align*}
which imply
\begin{align}\label{eq:dertg1112}
\hat{a}(\xi,t)=\bigg(\f{\lambda_+e^{\lambda_-t}-\lambda_-e^{\lambda_+t}}{\lambda_+-\lambda_{-}}\bigg)
\hat{a}_0(\xi)-\bigg(\f{e^{\lambda_+t}-e^{\lambda_-t}}{\lambda_+-\lambda_{-}}\bigg)|\xi|\hat{d}_0(\xi).
\end{align}
This determines $\widehat{\cG^{11}}$ and $\widehat{\cG^{12}}$.

On the other hand, taking Fourier transform to the second equation of (\ref{equ:linearized system}) gives
\begin{align*}
\hat{d}_t=-\bnu|\xi|^2\hat{d}+|\xi|\hat{a}.
\end{align*}
Thus, we have
\begin{align*}\label{}
\hat{d}(\xi,t)=e^{-\bnu|\xi|^2t}\Big[\hat{d}(\xi,0)+
|\xi|\int_0^te^{\bnu|\xi|^2\tau}\hat{a}(\xi,\tau)d\tau\Big].
\end{align*}
Plugging \eqref{eq:dertg1112} into the above equality and using the relations
$$\lambda_{\pm}+\bnu|\xi|^2=-\lambda_{\mp},\quad \lambda_-\lambda_+=|\xi|^2,$$
we finally get
\begin{align*}\label{}
\hat{d}(\xi,t)=-\bigg(\f{e^{\lambda_+t}-e^{\lambda_-t}}{\lambda_+-\lambda_{-}}\bigg)|\xi|\hat{a}_0(\xi)
+\bigg(\f{\lambda_+e^{\lambda_+t}-\lambda_-e^{\lambda_-t}}{\lambda_+-\lambda_{-}}\bigg)
\hat{d}_0(\xi),
\end{align*}
which determines $\widehat{\cG^{21}}$ and $\widehat{\cG^{22}}$.\ef

We have the following pointwise estimates and expansion for $\widehat{\mathcal {G}}$.

\bthm{Lemma}\label{Lem:pointwise estimate of Green function}
(a)\, Given $R>0$, there is a positive number $\vartheta$ depending on $R$ such that,
for any multi-indices $\alpha$ and $|\xi|\le R$,
\begin{align}\label{equ:lowfrequency of G}
|D_\xi^\alpha\widehat{\mathcal {G}}(\xi,t)|\le C e^{-\vartheta|\xi|^2t}(1+|\xi|)^{|\alpha|}(1+t)^{|\alpha|},
\end{align}where $C=C(R, |\al|)$;

(b)\, There exist positive constants $R, \vartheta$ depending on $\bar \nu$ such that the following  expansion is valid for $|\xi|\ge R$,
\begin{align}\label{equ:expansion of G}
\widehat{\mathcal{G}}(\xi,t)=&e^{-\bnu^{-1}t}\left[\begin{array}{ll}1\quad0\\0\quad0\end{array}\right]
+e^{-\bnu|\xi|^2t}\left[\begin{array}{ll}0\quad0\\0\quad1\end{array}\right]\nonumber\\
&+\widehat{\mathcal {G}}^1(\xi,t)\left[\begin{array}{ll}0\quad 1\\1\quad 0\end{array}\right]+
\widehat{\mathcal {G}}^2(\xi,t)\left[\begin{array}{ll}1\quad0\\0\quad1\end{array}\right],
\end{align}
where $\widehat{\mathcal {G}}^1$ and $\widehat{\mathcal {G}}^2$ satisfy the estimates
\begin{align}
&|\pa^\al_\xi\widehat{\mathcal {G}}^1|
\le C|\xi|^{-|\al|-1}\big(e^{-\frac{\bnu^{-1}t}{2}}+e^{-\vartheta|\xi|^2t}\big),\label{equ:highfrequency of G1}\\
&|\pa^\al_\xi \widehat{\mathcal {G}}^2|
\le C|\xi|^{-|\al|-2}\big(e^{-\frac{\bnu^{-1}t}{2}}+e^{-\vartheta|\xi|^2t}\big),\label{equ:highfrequency of G2}
\end{align}for a positive constant $C$ depending on $|\alpha|$, $\bnu$.
\ethm

\bthm{Remark}
In fact, $\widehat{\mathcal {G}}^2$ is a diagonal matrix, see \eqref{4-add1}. Since both nonzero elements in $\widehat{\mathcal {G}}^2$
can be estimated by the right side of \eqref{equ:highfrequency of G2}, we don't care about its explicit expression,
and  we view it as a scalar function.
\ethm
\noindent{\bf Proof.}\, We follow the proof of Lemma 3.2 in \cite{Hoff-Zum}.
Since the explicit expansion of $\widehat{\mathcal {G}}$ is important for us,
here we will present a proof. We refer to \cite{Hoff-Zum} for more details.

The inequality (\ref{equ:lowfrequency of G}) can be deduced from Lemma \ref{Lem:Greenmatrix} by routine calculations.
We present the proof of (b) below. Let $p$ be the symbol $$p(z,r)=z^2+\bnu r^2z+r^2.$$
Define contours $\Gamma_+, \Gamma_-$
and $\Gamma_0$ in the complex plane to be the circles of radius $\f{\bnu^{-1}}{2}$ centered respectively
at $-\bnu^{-1}$, $-\bnu r^2+\bnu^{-1}$ and 0, and
\begin{align*}
\left\{\begin{aligned}
&E(r,t)=\frac1{2\pi i}\int_{\Gamma_-\cup \Gamma_+}\f{e^{tz}}{p(z,r)}dz,\\
&F(r,t)=E_t(r,t)+\bnu r^2E(r,t),
\end{aligned}\right.
\end{align*}for large enough such that $r>10\bnu^{-1}$.
Then it is easy to verify that
\begin{align*}
\left\{\begin{array}{ll}
p(\pa_t,r)E=p(\pa_t,r)F=0,\\
E(r,0)=0,\quad F(r,0)=1,\\
E_t(r,0)=1,\quad  F_t(r,0)=0.
\end{array}\right.
\end{align*}
Thus, we may assume that $E$ and $F$ have the form
\begin{align*}
E(\xi,t)=A_E(\xi)e^{\lambda_{-}t}+B_E(\xi)e^{\lambda_{+}t},\\
F(\xi,t)=A_F(\xi)e^{\lambda_{-}t}+B_F(\xi)e^{\lambda_{+}t}.
\end{align*}
Thanks to the initial conditions on $E$ and $F$, we obtain
\begin{align*}
A_E=\f{-1}{\lambda_+-\lambda_{-}},\quad
B_E=\f{1}{\lambda_+-\lambda_{-}},\quad
A_F=\f{\lambda_+}{\lambda_+-\lambda_{-}},\quad
B_F=\f{-\lambda_-}{\lambda_+-\lambda_{-}}.
\end{align*}

Let $\hat{a}, \hat{d}$ be as in the proof of Lemma \ref{Lem:Greenmatrix}. By \eqref{equ:Fourierequ of a}, we find
\begin{align}
\hat{a}(\xi,t)=F(|\xi|,t)\hat{a}_0(\xi)-|\xi|E(|\xi|,t)\hat{d}_0(\xi),\quad
\widehat{\cG^{21}}=\widehat{\cG^{12}},\nonumber
\end{align}
which implies that
\begin{align}
\widehat{\cG^{11}}(\xi,t)=F(|\xi|,t),\quad\widehat{\cG^{12}}(\xi,t)=-|\xi|E(|\xi|,t),
\quad
\widehat{\cG^{21}}(\xi,t)=-|\xi|E(|\xi|,t).\nonumber
\end{align}

Introducing the function
\begin{align*}
H(r,t)&=e^{-\bnu r^2t}\int_0^te^{\bnu r^2\tau}E(r,\tau)d\tau
=\f{\lambda_-e^{\lambda_-t}-\lambda_+e^{\lambda_+t}}{|\xi|^2(\lambda_+-\lambda_-)},
\end{align*}
then we find
\begin{align}
\widehat{\cG^{22}}(\xi,t)=e^{-\bnu|\xi|^2t}-|\xi|^2H(\xi,t).\nonumber
\end{align}
Thus we have
\begin{align*}
\widehat{\cG}=&\left[\begin{array}{cc}e^{-\bnu^{-1}t}&0\\0&e^{-\bnu|\xi|^2t} \end{array}\right]
-|\xi|E(|\xi|,t)\left[\begin{array}{cc}0&1\\1&0 \end{array}\right]
\\&\quad+\left[\begin{array}{cc}F(|\xi|,t)-e^{-\bnu^{-1}t}&0\\0&-|\xi|^2H(\xi,t) \end{array}\right].
\end{align*}
Consequently, we obtain the required expansion \eqref{equ:expansion of G} of $\widehat{G}$ with
\begin{equation}\label{4-add1}
\widehat{\mathcal {G}}^1(\xi,t)=-|\xi|E(|\xi|,t),
\quad\widehat{\mathcal {G}}^2(\xi,t)
=\left[\begin{array}{cc}F(|\xi|,t)-e^{-\bnu^{-1}t}&0\\0&-|\xi|^2H(\xi,t) \end{array}\right].
\end{equation}
It remains to prove (\ref{equ:highfrequency of G1}) and (\ref{equ:highfrequency of G2}).
By making a change $w=z+\bnu r^2-\bnu^{-1}$, we  get
\begin{align*}
r\int_{\Gamma_-(r)}\f{e^{tz}}{p(z,r)}dz&=-e^{(\bnu^{-1}-\bnu r^2)t}\int_{\Gamma_0}\f{e^{tw}}{\bnu rw}
\Big[1-\f{(w+\bnu^{-1})^2}{\bnu r^2w}\Big]^{-1}dw\\
&=\sum_{j=0}^\infty r^{-2j-1}\Big[-e^{(\bnu^{-1}-\bnu r^2)t}
\int_{\Gamma_0}\f{(w+\bnu^{-1})^{2j}}{(\bnu w)^{j+1}}e^{tw}dw\Big].
\end{align*}for large enough $r>10\bnu^{-1}$, here $C$ is an absolute  constant.
Let $b_j(t)$ be the contour integral on the right hand side, we have
\begin{align*}|b_j(t)|\le C(\bnu)e^{\f{\bnu^{-1}t}2}(\widetilde{C}\bnu^{-1})^{2j},\end{align*}
where $\widetilde{C}$ is an absolutely constant. Then
if $r>2\bnu^{-1}$ is large enough, the term in bracket is bounded by
$$C(\bnu)e^{(\bnu^{-1}-\bnu r^2)t}e^{\f{\bnu^{-1}t}2}(\widetilde{C}\bnu^{-1})^{2j}
\le C(\bnu) e^{-\f{\bnu r^2t}{2}}(\widetilde{C}\bnu^{-1})^{2j}.$$
Taking $k$-derivatives with respect to $r$, we have
\begin{align*}\Big(\f{\pa}{\pa r}\Big)^k\int_{\Gamma_-}\f{re^{tz}}{p(z,r)}dz
=-e^{\bnu^{-1}t}\sum_{j=0}^\infty b_j(t)
\sum_{\ell=0}^k\Big(\begin{array}{ll}k\\\ell\end{array}\Big)\Big[\Big(\f{\pa}{\pa r}\Big)^{\ell}
e^{-\bnu r^2t}\Big]\Big[\Big(\f{\pa}{\pa r}\Big)^{k-\ell}r^{-2j-1}\Big].
\end{align*}
Noting that $\pa_r^{\ell}e^{-\bnu r^2t}\le C(\ell)r^{-\ell} e^{-\f23\bnu r^2t}$,
so if $r>2\widetilde{C}\bnu^{-1}$ is  large enough, we have
\begin{align}
&\bigg|\Big(\f{\pa}{\pa r}\Big)^k\int_{\Gamma_-}\f{re^{tz}}{p(z,r)}dz\bigg|\nonumber\\&
\le Ce^{(\bnu^{-1}-\f23\bnu r^2)t}\sum_{j=0}^\infty |b_j(t)|r^{-2j-1-k}\sum_{\ell=0}^k
\Big(\begin{array}{ll}k\\\ell\end{array}\Big)\f{(2j+k-\ell)!}{(2j)!}\nonumber\\&\le
C(\bnu,k)e^{-\f12\bnu r^2t}r^{-k-1}\sum_{j=0}^\infty\Big(\f{\widetilde{C}\bnu^{-1}}{r}\Big)^{2j}
\f{(2j+k-\ell)!}{(2j)!}\nonumber\\&\le C(\bnu,k)e^{-\f12\bnu r^2t}r^{-k-1}.\label{equ:Gamma-minus}
\end{align}

On the other hand, by making a change $w=z+\bnu^{-1}$, we have
\begin{align*}
r\int_{\Gamma_+(r)}\f{e^{tz}}{p(z,r)}dz&=e^{-\bnu^{-1}t}\int_{\Gamma_0}\f{e^{tw}}{\bnu r^2w}
\Big[1-\f{(w-\bnu^{-1})^2}{\bnu rw}\Big]^{-1}dw\\
&=\sum_{j=0}^\infty r^{-2j-1}(-1)^j\Big[e^{-\bnu^{-1}t}
\int_{\Gamma_0}\f{(w-\bnu^{-1})^{2j}}{(\bnu w)^{j+1}}e^{tw}dw\Big].
\end{align*}
The term in the bracket  is dominated by
$$C(\bnu)e^{-\f{\bnu^{-1}}{2}t}(\widetilde{C}\bnu^{-1})^{2j},$$
then by the same argument as leading to \eqref{equ:Gamma-minus}, we get
\begin{align}\label{equ:Gamma-plus}
&\bigg|\Big(\f{\pa}{\pa r}\Big)^k\int_{\Gamma_+}\f{re^{tz}}{p(z,r)}dz\bigg|
\le C(\bnu)e^{-\f{\bnu^{-1}}{2}t}r^{-k-1}.
\end{align}

With \eqref{equ:Gamma-minus} and \eqref{equ:Gamma-plus},
we can obtain (\ref{equ:highfrequency of G1}) and (\ref{equ:highfrequency of G2}) by
using the definitions of the functions $E, F$ and $H$.
\ef

Using Lemma \ref{Lem:pointwise estimate of Green function},
we can obtain the following smoothing effect of Green's matrix $\cG$, which will play an important role in this paper.

\bthm{Proposition}\label{Prop:Green-Lpest} %Given $R_0>0$, we decompose
% $$\cG=\cG_{\ell}+\cG_{h},$$
%where $\widehat{\cG}_{\ell}=\widehat{\cG}(\xi,t)\psi(\xi)$, and $\psi(\xi)$ is a smooth function
%supported on $|\xi|\le 2R_0$ and equals to 1 on $|\xi|\le R_0$.
Let ${\cC}$ be a ring  centered  at 0 in $\R^n$ . Then there exist  positive  constants
$R_0, C, c$ depending on $\bar \nu$ such that, if ${\rm
supp}\,\hat{u}\subset\lambda{\cC}$, then we have

(a) if $\lambda\le R_0$, then
\begin{align}\label{equ:lowGreen-L2}
&\|\mathcal{G}\ast u\|_{L^2}\le Ce^{-c\lambda^2
t}\|u\|_{L^2};
\end{align}

(b)\, if $b\le \lambda\le R_0$, then for any $1\le p\le \infty$,
\begin{align}\label{equ:lowGreen-Lp}
&\|\mathcal{G}\ast u\|_{L^p}\le C(1+b^{-n-1})e^{-c\lambda^2
t}\|u\|_{L^p};
\end{align}

(c)\, if $\lambda>R_0$, then for any $1\le p\le \infty$,
\begin{align}
&\|\mathcal{G}^1\ast u\|_{L^p}\le C\lambda^{-1}e^{-c
t}\|u\|_{L^p},\label{equ:highGreen1-Lp}\\
&\|\mathcal{G}^2\ast u\|_{L^p}\le C\lambda^{-2}e^{-c t}\|u\|_{L^p}.\label{equ:highGreen2-Lp}
\end{align}
\ethm

\noindent{\bf Proof.}\,(a)\, Thanks to Plancherel theorem, we get
\begin{align*}
\|\cG\ast u\|_{L^2}=\|\widehat{\cG}(\xi)\hat{u}(\xi)\|_{L^2}\le C
\|e^{-\vartheta|\xi|^2 t}\hat{u}(\xi)\|_2\le Ce^{-c\lambda^2 t}\|u\|_2,
\end{align*}
where  we have used \eqref{equ:lowfrequency of G} and the support property of $\hat{u}(\xi)$.

(b) Let $\phi\in {\cD}(\mathbb{R}^n\setminus\{0\})$, which equals to 1 near the  ring ${\cC}$.
Set
$$g(t,x)\eqdefa (2\pi)^{-n}\int_{\mathbb{R}^n}e^{ix\cdot\xi}
\phi(\lambda^{-1}\xi)\widehat{\mathcal{G}}(\xi,t)d\xi.
$$
To prove \eqref{equ:lowGreen-Lp}, it suffices to show
\begin{align}\label{equ:gL1}
\|g(x,t)\|_{L^1}\le C(1+b^{-n-1})e^{-c\lambda^2 t}.
\end{align}
Thanks to \eqref{equ:lowfrequency of G}, we infer that
\begin{align}\label{equ:gL2}
\int_{|x|\le\lam^{-1}}|g(x,t)|dx&\le C
\int_{|x|\le\lam^{-1}}\int_{\mathbb{R}^n}|\phi(\lambda^{-1}\xi)||\widehat{\mathcal{G}}(\xi,t)|d\xi dx
\le Ce^{-c\lambda^2 t}.
\end{align}
Set $L\eqdefa \frac {x\cdot \na_\xi} {i|x|^2}$. Noting that $L(e^{ix\cdot\xi})=e^{ix\cdot\xi}$,
we get by integration by part  that
\begin{align*}
g(x,t)=& \int_{\mathbb{R}^n}L^{n+1}(e^{ix\cdot\xi})\phi(\lambda^{-1}\xi)
\widehat{\mathcal{G}}(\xi,t)d\xi\nonumber\\=&(-1)^{n+1}
\int_{\mathbb{R}^n}e^{ix\cdot\xi}(L^*)^{n+1}\big(\phi(\lambda^{-1}\xi)
\widehat{\mathcal{G}}(\xi,t)\big)d\xi.
\end{align*}
From the Leibnitz's formula and \eqref{equ:lowfrequency of G},
\begin{align}
&\big|(L^*)^{n+1}\big(\phi(\lambda^{-1}\xi)\widehat{\mathcal{G}}(\xi,t)\big)\big|\nonumber\\&
\le C|\lambda x|^{-(n+1)}\sum_{|\gamma|= n+1,
|\beta|\le|\gamma|}\lam^{|\beta|}|(\na^{\gamma-\beta}\phi)(\lambda^{-1}\xi)|
e^{-\vartheta|\xi|^2t}(1+|\xi|)^{|\beta|}(1+t)^{|\beta|}.\nonumber
\end{align}
Thanks to the estimate
\begin{align*}
e^{-\vartheta|\xi|^2t}(1+|\xi|)^{|\beta|}(1+t)^{|\beta|}&\le
Ce^{-\vartheta|\xi|^2t}(1+t^{|\beta|}+|\xi|^{|\beta|}+t^{|\beta|}|\xi|^{|\beta|})\nonumber\\
&\le Ce^{-\f \vartheta 2|\xi|^2t}(1+|\xi|^{-2|\beta|}+|\xi|^{|\beta|}+|\xi|^{-|\beta|}),
\end{align*}
we obtain, for any $\xi$ with $b\lesssim |\xi|\sim \lambda\le R_0$,
\begin{align*}
\big|(L^*)^{n+1}\big(\phi(\lambda^{-1}\xi)\widehat{\mathcal{G}}(\xi,t)\big)\big|\le C(1+b^{-n-1})
|\lambda x|^{-(n+1)}e^{-\f \vartheta 2|\xi|^2t},
\end{align*}
which implies that
\begin{align*}
\int_{|x|\ge\f 1\lam}|g(x,t)|dx&\le C(1+b^{-n-1})
e^{-c\lambda^2t}\lambda^n\int_{|x|\ge\f 1\lam}|\lambda
x|^{-n-1}dx\nonumber\\&\le C(1+b^{-n-1})
e^{-c\lambda^2t},
\end{align*}
which together with \eqref{equ:gL2} gives \eqref{equ:gL1}.

(c)\, We set
\begin{align*}
g^1(t,x)\eqdefa (2\pi)^{-n}\int_{\mathbb{R}^n}e^{ix\cdot\xi}\phi(\lambda^{-1}\xi)
\widehat{\mathcal{G}}^1(\xi,t)d\xi.\end{align*}
Since the integrand is supported in
$\{\xi: |\xi|\sim \lam>R_0\}$,
then  we get by \eqref{equ:highfrequency of G1} that
\begin{align}\label{equ:gH1}
\int_{|x|\le \lam^{-1}}|g^1|dx&\le C
\int_{|x|\le \lam^{-1}}\int|\phi(\lambda^{-1}\xi)||\xi|^{-1}\big(e^{-\frac{\bnu^{-1}t}{2}}+e^{-\vartheta|\xi|^2t}\big)d\xi
dx\nonumber\\ &\le C\lambda^{-1}e^{-ct}.
\end{align}
We get by integration by parts that
\begin{align*}
g^1(x,t)=(-1)^{n+1}\int_{\mathbb{R}^n}e^{ix\cdot\xi}(L^*)^{n+1}\bigl(\phi(\lambda^{-1}\xi)
\widehat{\mathcal{G}}^1(\xi,t)\bigr)d\xi.
\end{align*}
While thanks to \eqref{equ:highfrequency of G1}, the  integrand is dominated by
\begin{align*}
|\lambda x|^{-n-1}\lambda^{-1}e^{-ct},
\end{align*}
which implies that
\begin{align}\label{Preq:linGm29}
\int_{|x|\ge\lam^{-1}}|g^1(x,t)|dx
&\le C\lambda^{-1}e^{-ct}\lambda^{n}\int_{|x|\ge\lam^{-1}}|\lambda
x|^{-n-1}dx
\nonumber\\
&\le C\lambda^{-1}e^{-ct},
\end{align}
from which and \eqref{equ:gH1}, it follows that
$$
\|g^1(x,t)\|_{L^1}\le C\lambda^{-1}e^{-ct},
$$
which implies \eqref{equ:highGreen1-Lp}. The inequality \eqref{equ:highGreen2-Lp} can be similarly proved.\ef

\section{The linearized  system with convection term}

In this section we consider the linearized system with convection terms
\begin{equation}\label{equ:linerized coupling system}
\left\{
\begin{array}{ll}
\p_ta+\Lambda d=-v\cdot\na a+F,\\
\p_td-\bar{\nu}\Delta d-\Lambda a=-v\cdot\na d+G, \\
(a,d)|_{t=0}=(a_0,d_0).
\end{array}
\right.
\end{equation}

Let $T>0$, $s\in\mathbb{R}, p\ge 2,$ and $s_p\eqdefa s-\f np+\f n2$. We introduce the functional space $\cE^s_T$
which is defined as follows
\begin{align*}
\cE^{s}_T\eqdefa \bigg\{(a,d)\in \Big(&L^1(0,T; \dot
\cB^{s_p+1,s}_{2,p})
\cap \widetilde{L}^\infty(0,T; \dot \cB^{s_p-1,s}_{2,p})\Big)\\
&\times\Big(L^1(0,T; \dot \cB^{s_p+1,s+1}_{2,p})\cap \widetilde{L}^\infty(0,T; \dot \cB^{s_p-1,
s-1}_{2,p})\Big)^n\bigg\}\end{align*}
with the norm
$$\|(a,d)\|_{\cE^{s}_T}\eqdefa \|a\|_{\widetilde{L}^\infty_T\dot \cB^{s_p-1,s}_{2,p}}+
\|d\|_{\widetilde{L}^\infty_T\dot \cB^{s_p-1,s-1}_{2,p}}+\|a\|_{L^1_T\dot \cB^{s_p+1,s}_{2,p}}+
\|d\|_{L^1_T\dot \cB^{s_p+1,s+1}_{2,p}}.$$
We denote
$$
\overline{V}(t)\eqdefa \|v\|_{{L}^1_t\dot \cB^{\f n2+1, \f np+1}_{2,p}}+
\|v\|_{\widetilde{L}^\infty_t\dot \cB^{ \f n2-1,\f np-1}_{2,p}}.
$$

\bthm{Theorem}\label{thm:hyper-prabolic estimate}
Let  $2\le p<2n$, $p\le\min\big(
4,\f {2n}{n-2}\big)$, and $1-\f np<s\le1+\f{2n}p-\f n2$.
Assume that
$v\in \widetilde{L}^\infty_T\dot \cB^{ \f n2-1,\f np-1}_{2,p}\cap{L}^1_T\dot \cB^{ \f n2+1,\f np+1}_{2,p},
F\in {L}^1_T\dot{\cB}^{s_p-1, s}_{2,p},
G\in {L}^1_T\dot{\cB}^{s_p-1, s-1}_{2,p}$.
Let $(a,d)$ be a solution of \eqref{equ:linerized coupling system} on $[0,T]$.
Then there exists a constant $C$ independent of $T$ such that
\begin{align*}
\|(a,d)\|_{\cE^{s}_T}\le Ce^{C\overline{V}(T)}
\Big\{&\|(a_0,d_0)\|_{\cE^{s}_0}
+\big(\overline{V}(T)+\overline{V}(T)^\f12
\big)\|(a,d)\|_{\cE^{s}_T}\\&\qquad\qquad\qquad+
\|F\|_{{L}^1_T\dot{\cB}^{s_p-1, s}_{2,p}}+\|G\|_{{L}^1_T\dot{\cB}^{s_p-1, s-1}_{2,p}}
\Big\}.
\end{align*}
Here $\|(a_0,d_0)\|_{\cE^{s}_0}\eqdefa \|a_0\|_{\dot \cB^{s_p-1, s}_{2,p}}+\|d_0\|_{\dot \cB^{s_p-1, s-1}_{2,p}}$.
\ethm

We introduce another functional space $E^s_T$ which is defined by
\begin{align*}
E^{s}_T\eqdefa \bigg\{(a,d)\in \Big(&L^1(0,T; \dot\cB^{s+1,s})
\cap \widetilde{L}^\infty(0,T; \dot\cB^{s-1,s})\Big)\\
&\times\Big(L^1(0,T; \dot B^{s+1}_{2,1})\cap \widetilde{L}^\infty(0,T; \dot  B^{s-1
}_{2,1})\Big)^n\bigg\}
\end{align*}
with the norm
$$\|(a,d)\|_{E^{s}_T}\eqdefa \|a\|_{\widetilde{L}^\infty_T{\dot\cB}^{s-1,s}}+
\|d\|_{\widetilde{L}^\infty_T\dot  B^{s-1}_{2,1}}+\|a\|_{L^1_T{\dot\cB}^{s+1,s}}+
\|d\|_{L^1_T\dot  B^{s+1}_{2,1}}.$$

\bthm{Theorem}\label{thm:hyper-prabolic estimate-2}
Let $2\le p<2n$, $p\le\min\big(
4,\f {2n}{n-2}\big)$, and $1-\f np<s\le1+\f np$. Assume that
$v\in \widetilde{L}^\infty_T\dot \cB^{ \f n2-1,\f np-1}_{2,p}\cap L^1_T\dot \cB^{ \f n2+1,\f np+1}_{2,p},
F\in  L^1_T\dot{\cB}^{s-1, s},
G\in  L^1_T\dot{B}^{s-1}_{2,1}$. Let  $(a,d)$ be a solution of \eqref{equ:linerized coupling system}
on $[0,T]$. Then there exists a constant $C$ independent of $T$ such that
\begin{align*}
\|(a,d)\|_{E^{s}_T}\le& Ce^{C\overline{V}(T)}
\Big\{\|(a_0,d_0)\|_{E^{s}_0}
+\big(\overline{V}(T)+\overline{V}(T)^\f12
\big)\|(a,d)\|_{E^{s}_T}\\
&\quad+\|v\|_{\widetilde{L}^2_T\dot B^{s}_{2,1}}\|(a,d)\|_{\widetilde{L}^2_T\dot\cB^{\f n2,\f {n}p}_{2,p}}
+\|F\|_{L^1_T{\dot\cB}^{s-1, s}_{2,1}}+\|G\|_{L^1_T\dot{B}^{s-1}_{2,1}}
\Big\}.
\end{align*}
Here $ \|(a_0,d_0)\|_{E^{s}_0}\eqdefa \|a_0\|_{\dot \cB^{s-1, s}}+\|d_0\|_{\dot B^{s-1}_{2,1}}$.
\ethm

We denote
$$a_j\eqdefa\Delta_j a,\quad d_j\eqdefa\Delta_j
d,\quad F_j\eqdefa\Delta_j F, \quad G_j\eqdefa\Delta_j G.
$$
The proof of Theorem \ref{thm:hyper-prabolic estimate}-\ref{thm:hyper-prabolic estimate-2}
is based on the following frequency localized system:
\begin{equation}\label{equ:linsys-local}
\left\{
\begin{array}{ll}
\p_ta_j+\Lambda d_j=-\Delta_j(v\cdot\na a)+F_j,\\
\p_td_j-\bar{\nu}\Delta d_j-\Lambda a_j=-\Delta_j(v\cdot\na d)+G_j, \\
(a_j,d_j)|_{t=0}=(\Delta_ja_0,\Delta_jd_0)\eqdefa (a_j^0,d_j^0).
\end{array}
\right.
\end{equation}

In what follows, we define $\widetilde{p}$ as $\f 1{\tilde{p}}={\f n{2p}-\f n4+\f12}$, and
denote $p'$, $\tilde{p}'$ by the conjugate index of $p'$, $\tilde{p}$ respectively.

\subsection{Low frequency estimate: $2^j\le R_0$}
In this case, the Green matrix of the linearized system behaves as the heat kernel.
Using the smoothing effect of the Green matrix, the terms $v\cdot\na a$ and $v\cdot\na d$ can be handled as the perturbation terms.

In terms of the Green matrix, the solution of (\ref{equ:linsys-local}) can be expressed as
\begin{align*}
\left(\begin{array}{ll} \!a_j \!\!\vspace{.15cm}\\ \!d_j \!\end{array}\right)=
\mathcal{G}(x,t)\ast\left(\begin{array}{ll} \!a_j^0 \!\vspace{.15cm}\\ \!d_j^0 \!\end{array}\right)
+\int_0^t\mathcal{G}(x,t-\tau)\ast\left(\begin{array}{ll}\!F_j-\Delta_j(v\cdot\na a)\!\vspace{.15cm}\\
\!G_j-\Delta_j(v\cdot\na d)\!\end{array}\right)d\tau.
\end{align*}
From Proposition \ref{Prop:Green-Lpest} (a) and Young's inequality, we infer that
\begin{align*}
\|a_j(t)\|_{L^2}&+\|d_j(t)\|_{L^2}\le Ce^{-c2^{2j} t}\bigl(\|a_j^0\|_{L^2}+\|d_j^0\|_{L^2}\bigr)\\
&+C\int_0^te^{-c2^{2j} (t-\tau)}
\bigl(\|F_j(\tau)\|_{L^2}+\|G_j(\tau)\|_{L^2}\bigr)d\tau\\&
+C\int_0^te^{-c2^{2j} (t-\tau)}\bigl(\|\Delta_j(v\cdot\na a)(\tau)\|_{L^2}+\|\Delta_j(v\cdot\na d)(\tau)\|_{L^2}\bigr)d\tau.
\end{align*}
Taking $L^r$ norm with respect to $t$ gives
\begin{align*}
\|a_j\|_{L^r_TL^2}+\|d_j\|_{L^r_TL^2}\le C2^{-\f {2j} r}&
\bigl(\|a_j^0\|_{L^2}+\|d_j^0\|_{L^2}
+\|{F}_j\|_{L^1_TL^2}+\|\Delta_j(v\cdot\na a)\|_{L^1_TL^2}\\
&+\|{G}_j\|_{L^1_TL^2}+\|\Delta_j(v\cdot\na d)\|_{L^1_TL^2}\bigr).
\end{align*}
Here and what follows, we assume $1\le r\le \infty$.
Noting that $1-\f np<s\le \f {2n}p-\f n2+1$, we apply Proposition \ref{Prop:Product} (b)
with $s=\f n p, t=s-\f n p+\f n {p'}-1, \widetilde{s}=s_p-1, \widetilde{t}=\f n2,\gamma=\f {n}p-\f n2 $ to get
\begin{align*}
&\sum_{2^j\le R_0}2^{j(s_p-1)}\|\Delta_j(v\na a)\|_{L^1_TL^2}\nonumber\\&\le C
\|v\|_{\widetilde{L}^{\tilde{p}}_T\dot\cB^{\f np,\f {2n}p-\f n2}_{2,p}}
\|\na a\|_{\widetilde{L}^{\tilde{p}'}_T\dot\cB^{s-\f np+\f n{p'}-1,s_p-1+\gamma}_{2,p}}+C
\|\na a\|_{\widetilde{L}^2_T\dot\cB^{s_p-1,s-1}_{2,p}}\|v\|_{\widetilde{L}^2_T\dot\cB^{\f n2,\f np}_{2,p}}\nonumber\\
&\le C\|v\|_{\widetilde{L}^{\tilde{p}}_T\dot\cB^{\f np,\f {2n}p-\f n2}_{2,p}}
\|a\|_{\widetilde{L}^{\tilde{p}'}_T\dot\cB^{s-\f np+\f n{p'},s}_{2,p}}+
C\|a\|_{\widetilde{L}^2_T\dot\cB^{s_p,s}_{2,p}}\|v\|_{\widetilde{L}^2_T\dot\cB^{\f n2,\f np}_{2,p}},
\end{align*}
and with $ \gamma=0$,
\begin{align*}
&\sum_{2^j\le R_0}2^{j(s_p-1)}\|\Delta_j(v\na d)\|_{L^1_TL^2}\nonumber\\&\le C
\|v\|_{\widetilde{L}^{\tilde{p}}_T\dot\cB^{\f np,\f {2n}p-\f n2}_{2,p}}
\|\na d\|_{\widetilde{L}^{\tilde{p}'}_T\dot\cB^{s-\f np+\f n{p'}-1,s_p-1+\gamma}_{2,p}}+C
\|\na d\|_{\widetilde{L}^2_T\dot\cB^{s_p-1,s-1}_{2,p}}\|v\|_{\widetilde{L}^2_T\dot\cB^{\f n2,\f np}_{2,p}}\nonumber\\
&\le C\|v\|_{\widetilde{L}^{\tilde{p}}_T\dot\cB^{\f np,\f {2n}p-\f n2}_{2,p}}
\|d\|_{\widetilde{L}^{\tilde{p}'}_T\dot\cB^{s-\f np+\f n{p'},s_p}_{2,p}}+C
\|d\|_{\widetilde{L}^2_T\dot\cB^{s_p,s}_{2,p}}\|v\|_{\widetilde{L}^2_T\dot\cB^{\f n2,\f np}_{2,p}}.
\end{align*}
Then we conclude that for any $1\le r\le \infty$,
\begin{align}\label{equ:(a,d)-Lowest}
&\sum_{2^j\le R_0}2^{j(s_p-1+\f2r)}(\|a_j\|_{L^r_TL^2}+\|d_j\|_{L^r_TL^2})\nonumber\\&\le C
\sum_{2^j\le  R_0}2^{j(s_p-1)}(\|a_j^0\|_{L^2}+\|d_j^0\|_{L^2})\nonumber\\ &\quad+C
\Big\{\|F\|_{ L^1_T\dot \cB^{s_p-1,s}_{2,p}}+\|G\|_{ L^1_T\dot \cB^{s_p-1,s-1}_{2,p}}
+\|v\|_{\widetilde{L}^2_T\dot \cB^{\f n 2,\f n p}_{2,p}}\|(a,d)\|_{\widetilde{L}^2_T\dot \cB^{{s_p,s}}_{2,p}}
\nonumber\\&\qquad\quad+
\|v\|_{\widetilde{L}^{\tilde{p}}_T\dot\cB^{\f np,\f {2n}p-\f n2}_{2,p}}
\bigl(\|a\|_{\widetilde{L}^{\tilde{p}'}_T\dot\cB^{s-\f np+\f n{p'},s}_{2,p}}+
\|d\|_{\widetilde{L}^{\tilde{p}'}_T\dot\cB^{s-\f np+\f n{p'},s_p}_{2,p}}\bigr)\Big\}.
\end{align}

\subsection{High frequency estimate: $2^j>R_0$}
 In this case, since the Green's matrix has no smoothing effect on $a$,
we will lose one derivative about $a$ in the energy estimate if we still regard $v\cdot\na a$  as the perturbation term.
 To avoid the loss of the derivative, we will work in the Lagrangian coordinate in this subsection.

Firstly, we rewrite \eqref{equ:linsys-local} as
\begin{align}\label{equ:linsys-local-1}
\left\{
\begin{aligned}
&\p_ta_j+S_{j-1}v\cdot\na a_j+\Lambda d_j=(S_{j-1}v\cdot\na a_j-\Delta_j(v\cdot\na a))+F_j,\\
&\p_td_j+S_{j-1}v\cdot\na d_j-\bar{\nu}\Delta d_j-\Lambda a_j=(S_{j-1}v\cdot\na d_j
-\Delta_j(v\cdot\na d))+G_j, \\
&(a_j,d_j)|_{t=0}=(a_j^0,d_j^0).
\end{aligned}\right.
\end{align}
Let  $\psi_j(t,x)$  be the solution to
\begin{align*}
&\frac{d}{dt}\psi_j(t,x)=S_{j-1}v(t,\psi_j(t,x)),\quad\psi_j(0,x)=x.
\end{align*}
We set
$$\ba_j=a_j(t,\psi_j(t,x)),\quad\bd_j=d_j(t,\psi_j(t,x)),\quad
\bar{F}_j=F_j(t,\psi_j(t,x)),\quad\bar{G}_j=G_j(t,\psi_j(t,x)).$$
Then the new unknown $(\ba_j,\bd_j)$ satisfies
\begin{align}\label{eq:linsys-locLancor}\left\{
\begin{aligned}
&\p_t\ba_j+\Lambda \bd_j={\bar F}_j+\bar{\mathcal {R}}_j,\\
&\p_t\bd_j-\bar{\nu}\Delta \bd_j-\Lambda \ba_j={\bar G}_j+\bar{\mathcal {Q}}_j,\\
&(\bar{a}_j,\bar{d}_j)|_{t=0}=(\Delta_ja_0,\Delta_jd_0),
\end{aligned}\right.
\end{align}
where
\begin{align*}
&\bar{\mathcal {R}}_j\eqdefa \big(S_{j-1}v\cdot\na a_j-\Delta_j(v\cdot\na a)\big)\circ\psi_j
+\Lambda(d_j\circ\psi_j)-(\Lambda d_j)\circ\psi_j,\\
&\bar{\mathcal {Q}}_j\eqdefa \big((S_{j-1}v\cdot\na d_j)-\Delta_j(v\cdot\na d)\big)\circ\psi_j-\bar{\nu}\Delta(d_j\circ\psi_j)
+\bnu(\Delta d_j)\circ\psi_j\\
&\qquad\quad+(\Lambda a_j)\circ\psi_j-\Lambda(a_j\circ\psi_j).
\end{align*}

We get by Lemma \ref{Lem:com-est} that
\begin{align}\label{equ:a-Lp}
&\|a_j\|_{L^p}\le e^{\overline{V}(t)}\|\ba_j\|_{L^p},
\end{align}
where we used the fact
\begin{align*}
\int_0^t\|\na v(\tau)\|_{L^\infty}d\tau\le C\int_0^t\|\na v(\tau)\|_{\dot{\cB}^{ \f n2,\f np}_{2,p}}d\tau
\le C\overline{V}(t).
\end{align*}
Fix $N\in\N$ to be chosen later. Noting that
$$
\ba_j=\displaystyle\sum_{k}\Delta_k\ba_j=\displaystyle\sum_{|k-j|\le N}\Delta_k\ba_j+
\displaystyle\sum_{k-j>N}\Delta_k\ba_j+S_{j-N}\ba_j,
$$
we have by (\ref{equ:a-Lp}) that
\begin{align*}
\|a_j\|_{L^p}\le e^{\overline{V}(t)}\Big(\sum_{|k-j|\le N}\!\|\Delta_k\ba_j\|_{L^p}+\!\sum_{k-j>N}\|\Delta_k\ba_j\|_{L^p}
+\|S_{j-N}\ba_j\|_{L^p}\Big).
\end{align*}
Thanks to  Lemma \ref{Lem:com-localest}, we have
\begin{align*}
&\!\sum_{k-j>N}\|\Delta_k\ba_j\|_{L^p}\le Ce^{C\overline{V}(t)}\sum_{k-j>N}2^{-(k-j)}\|a_j\|_{L^p}
\le C2^{-N}e^{C\overline{V}(t)}\|a_j\|_{L^p},\\
&\|S_{j-N}\ba_j\|_{L^p}\le Ce^{C\overline{V}(t)}
\bigl(\overline{V}(t)\|a_j\|_{L^p}+2^{-N}\|a_j\|_{L^p}\bigr).
\end{align*}

Summing up the above estimates yields
\begin{align*}
\|a_j\|_{L^p}
\le Ce^{C\overline{V}(t)}\Bigl(\sum_{|k-j|\le N}\|\Delta_k\ba_j\|_{L^p}+2^{-N}\|a_j\|_{L^p}+\overline{V}(t)
\|a_j\|_{L^p}\Bigr).
\end{align*}
Choosing $N$ large enough such that
\begin{align*}
Ce^{C\overline{V}(T)}2^{-N}\le \f12 \quad \textrm{i.e.}, N\sim e^{C\overline{V}(T)},
\end{align*}
then we obtain
\begin{align}\label{equ:linsy-ahigh}
\|a_j\|_{L^r_TL^p}\le Ce^{C\overline{V}(T)}\Bigl(
\sum_{|k-j|\le N}\|\Delta_k\ba_j\|_{L^r_TL^p}+\overline{V}(T)\|a_j\|_{L^r_TL^p}\Bigr).
\end{align}
Similarly, we can obtain
\begin{align}\label{equ:linsy-dhigh}
\|d_j\|_{L^r_TL^p}\le Ce^{C\overline{V}(T)}\Bigl(
\sum_{|k-j|\le N}\|\Delta_k\bd_j\|_{L^r_TL^p}+\overline{V}(T)
\|d_j\|_{L^r_TL^p}\Bigr).
\end{align}

Next, we need to estimate $\|\Delta_k\ba_j\|_{L^r_TL^p}$ and $\|\Delta_k\bd_j\|_{L^r_TL^p}$ for $|k-j|\le N$.
In order to do this,
we first present the estimates of $\bar{F}_j, \bar{G}_j, \bar{\mathcal {R}}_j$ and $\bar{\mathcal {Q}}_j$.
We get by using Lemma \ref{Lem:com-localest} that
\begin{equation}\label{equ:FG-Lancor}\begin{split}
&\|\Delta_k\bar{F}_j\|_{L^1_TL^p}
\le  C2^{-(k-j)}e^{C\overline{V}(T)}\|F_j\|_{L^1_TL^p},\\
&\|\Delta_k\bar{G}_j\|_{L^1_TL^p}\le  C2^{-(k-j)}e^{C\overline{V}(T)}\|G_j\|_{L^1_TL^p}.
\end{split}
\end{equation}
Set
\begin{align*}
&\bar{\mathcal {R}}^1_j=\big(S_{j-1}v\cdot\na a_j-\Delta_j(v\cdot\na a)\big)\circ\psi_j,\\&
\bar{\mathcal {R}}^2_j=\Lambda(d_j\circ\psi_j)-(\Lambda d_j)\circ\psi_j,\\&
\bar{\mathcal {Q}}^1_j=\big((S_{j-1}v\cdot\na d_j)-\Delta_j(v\cdot\na d)\big)\circ\psi_j,\\&
\bar{\mathcal {Q}}^2_j=\bnu(\Delta d_j)\circ\psi_j-\bar{\nu}\Delta(d_j\circ\psi_j)
+(\Lambda a_j)\circ\psi_j-\Lambda(a_j\circ\psi_j),
\end{align*}
Using Proposition \ref{Prop:Commutator} with $s=s_p, \sigma=s$ and Lemma \ref{Lem:Bernstein}, we get
\begin{align*}
&\|S_{j-1}v\cdot\na a_j-\Delta_j(v\cdot\na a)\|_{L^1_TL^p}\\
&\le \|[v,\Delta_j]\cdot\na a\|_{L^1_TL^p}+\|(S_{j-1}v-v)\cdot\na a_j\|_{L^1_TL^p}\nonumber\\&
\le Cc(j)2^{-s j}
\|v\|_{ L^1_T\dot{\cB}^{\frac n 2+1,\frac n p+1}_{2,p}}\|a\|_{\widetilde{L}^\infty_T\dot{\cB}^{s_p,s}_{2,p}}
+C\sum_{j'\ge j-1}2^{j'}\|\Delta_{j'}v\|_{L^1_TL^\infty}2^{j-j'}\|a_j\|_{L^\infty_TL^p}\nonumber\\&\le C
c(j)2^{-s j}\|v\|_{ L^1_T\dot{\cB}^{\frac n 2+1,\frac n p+1}_{2,p}}
\|a\|_{\widetilde{L}^\infty_T\dot{\cB}^{s_p,s}_{2,p}},
\end{align*}from which, Lemma \ref{Lem:com-est}, and Lemma \ref{Lem:Besov-properties} (a), we deduce that
\begin{align}\label{equ:R_1-Lancor}
\big\|\Delta_k\bar{\mathcal {R}}^1_j\big\|_{L^1_TL^p}
&\le Ce^{C\overline{V}(T)}\|S_{j-1}v\cdot\na a_j-\Delta_j(v\cdot\na a)\|_{L^1_TL^p}
\nonumber\\&\le Cc(j)
2^{-s j}e^{C\overline{V}(T)}\overline{V}(T)\|a\|_{\widetilde{L}^\infty_T\dot{\cB}^{s_p-1,s}_{2,p}}.
\end{align}
From Lemma \ref{Lem:com-commutator}, it follows that
\begin{align}\label{equ:R_2-Lancor}
\big\|\Delta_k\bar{\mathcal {R}}^2_j(t)\big\|_{L^p}\le C2^je^{C\overline{V}(t)}\overline{V}(t)^{\f12}\|d_j\|_{L^p}.
\end{align}
Similarly, we have
\begin{align}\label{equ:Q_1-Lancor}
&\big\|\Delta_k\bar{\cQ}^1_j\big\|_{L^1_TL^p}
\le Cc(j)2^{-(s-1)j}e^{C\overline{V}(T)}\overline{V}(T)
\|d\|_{\widetilde{L}^\infty_T\dot{\cB}^{s_p-1,s-1}_{2,p}},\\
&\big\|\Delta_k\bar{\cQ}^2_j(t)\big\|_{L^p}
\le Ce^{C\overline{V}(t)}\overline{V}(t)^{\f12}\bigl(2^{2j}\|d_j\|_{L^p}+2^{j}\|a_j\|_{L^p}\bigr).\label{equ:Q_2-Lancor}
\end{align}

Now we are in position to estimate $\|(\Delta_k\ba_j,\Delta_k\bd_j)\|_{L^r_TL^p}$.
 Since $(\ba_j,\bd_j)$ satisfies (\ref{eq:linsys-locLancor}),
it can be expressed in terms of  Green's matrix as
\begin{align*}
\left(\begin{array}{ll} \Delta_k\ba_j \vspace{.15cm} \\ \Delta_k\bd_j\end{array}\right)=
\mathcal{G}(x,t)\ast\left(\begin{array}{ll} \Delta_k a_j^0
\vspace{.15cm} \\ \Delta_k d_j^0\end{array}\right)
+\int_0^t\mathcal{G}(x,t-\tau)\ast\left(\begin{array}{ll} \Delta_k\bar{F}_j
+\Delta_k\bar{\mathcal {R}}_j\vspace{.15cm} \\ \Delta_k\bar{G}_j
+\Delta_k\bar{\mathcal {Q}}_j\end{array}\right)(\tau)d\tau.
\end{align*}
We consider two cases:
\vspace{0.3cm}

{\bf Case 1.}  $2^k>R_0$.
 By Lemma \ref{Lem:pointwise estimate of Green function} (b),  we have
\begin{align*}\Delta_k\ba_j(t,x)=&e^{-\bnu^{-1} t}\Delta_k a_j^0+\mathcal{G}^1*\Delta_kd_j^0
+\mathcal{G}^2*\Delta_ka_j^0
\\&+\int_0^t\big\{e^{-\bnu^{-1} (t-\tau)}+\mathcal{G}^2(x,t-\tau)\ast\big\}
\big(\Delta_k\bar{F}_j+\Delta_k\bar{\mathcal {R}}_j\big)(\tau)d\tau
\\&+\int_0^t\mathcal{G}^1(x,t-\tau)\ast\big(\Delta_k\bar{G}_j+\Delta_k
\bar{\mathcal {Q}}_j\big)(\tau)d\tau, \end{align*}
where the scale function ${\mathcal G}^2$ corresponds to the first nonzero in the diagonal matrix \eqref{4-add1}.
 Noting that $2^k>R_0$, we apply Proposition \ref{Prop:Green-Lpest} (c) to get
\begin{align*}
\|\Delta_k\ba_j(t)\|_{L^p}\le & Ce^{-ct}\|\Delta_k a_j^0\|_{L^p}+C2^{-k}e^{-ct}\|\Delta_kd_j^0\|_{L^p}
\nonumber\\&+C\int_0^te^{-c(t-\tau)}\big(\|\Delta_k\bar{F}_j(\tau)\|_{L^p}+\|\Delta_k\bar{\mathcal {R}}_j(\tau)\|_{L^p}\big)d\tau
\nonumber\\&+C2^{-k}\int_0^te^{-c(t-\tau)}
\big(\|\Delta_k\bar{G}_j(\tau)\|_{L^p}+\|\Delta_k\bar{\mathcal {Q}}_j(\tau)\|_{L^p}\big)d\tau.
\end{align*}
Taking $L^r$ norm with respect to $t$ and using Young's inequality give
\begin{align}\label{eq:linhhfrc}
\|\Delta_k\ba_j\|_{L^r_TL^p}\le
&C\bigl(\|\Delta_k a_j^0\|_{L^p}+2^{-k}\|\Delta_kd_j^0\|_{L^p}+
\|\Delta_k\bar{F}_j\|_{L^1_TL^p}+\|\Delta_k\bar{\mathcal {R}}_j\|_{L^1_TL^p}\bigr)\nonumber\\&\quad+C
2^{-k}\big(\|\Delta_k\bar{G}_j\|_{L^1_TL^p}+\|\Delta_k\bar{\mathcal {Q}}_j\|_{L^1_TL^p}\big)
.
\end{align}
On the other hand, we have
\begin{align*}\Delta_k\bd_j(t,x)=&\mathcal{G}^1*\Delta_k a_j^0+e^{\bnu\Delta t}\Delta_kd_j^0
+\mathcal{G}^2*\Delta_kd_j^0\\&
+\int_0^t\mathcal{G}^1(x,t-\tau)\ast\big(\Delta_k\bar{F}_j+\Delta_k\bar{\mathcal {R}}_j\big)(\tau)d\tau
\\&+\int_0^t\big(e^{\bnu\Delta (t-\tau)}+\mathcal{G}^2(x,t-\tau)\ast\big)
\big(\Delta_k\bar{G}_j+\Delta_k\bar{\mathcal {Q}}_j\big)(\tau)d\tau,\end{align*}
where the scale function ${\mathcal G}^2$ corresponds to the second nonzero in the diagonal matrix \eqref{4-add1}.
In analogy as leading to \eqref{eq:linhhfrc}, we get
\begin{align}\label{eq:linhhfrd}
\|\Delta_k\bd_j\|_{L^r_TL^p}\le &C\Bigl(
2^{-k}\|\Delta_k a_j^0\|_{L^p}+2^{-\f 2rk}\|\Delta_kd_j^0\|_{L^p}
\nonumber\\
&+2^{-k}\big(\|\Delta_k\bar{F}_j\|_{L^1_TL^p}+\|\Delta_k\bar{\mathcal {R}}_j^1\|_{L^1_TL^p}\big)
+2^{-k}\|\Delta_k\bar{\mathcal {R}}_j^2\|_{L^r_TL^p}\nonumber\\
&+2^{-\f2rk}
\big(\|\Delta_k\bar{G}_j\|_{L^1_TL^p}+\|\Delta_k\bar{\mathcal {Q}}_j^1\|_{L^1_TL^p}\big)
+2^{-2k}\|\Delta_k\bar{\mathcal {Q}}_j^2\|_{L^r_TL^p}\Bigr).
\end{align}
Substituting \eqref{equ:FG-Lancor}-\eqref{equ:Q_2-Lancor} into
\eqref{eq:linhhfrc} and  \eqref{eq:linhhfrd}, then summing over $k$, we obtain
\begin{align}\label{equ:a-hhest}
&\sum_{|k-j|\le N,2^k>R_0}\!
\|\Delta_k\ba_j\|_{L^r_TL^p}\le Ce^{C\overline{V}(T)}
\Big\{\|a_j^0\|_{L^p}+2^{-j}\|d_j^0\|_{L^p}\nonumber\\&\qquad\qquad+
\overline{V}(T)^\f12\bigl(\|a_j\|_{L^1_TL^p}+2^j\|d_j\|_{L^1_TL^p}\bigr)
+\|F_j\|_{L^1_TL^p}+2^{-j}\|G_j\|_{L^1_TL^p}\nonumber\\&\qquad\qquad\qquad+2^{-sj}c(j)\overline{V}(T)
\big(\|a\|_{\widetilde{L}^\infty_T\dot{\cB}^{s_p-1,s}_{2,p}}
+\|d\|_{\widetilde{L}^\infty_T\dot{B}^{s_p-1,s-1}_{2,p}}\big)\Big\},
\end{align}
and
\begin{align}\label{equ:d-hhest}
&\sum_{|k-j|\le N,2^k>R_0}\!\!\!
\|\Delta_k\bd_j\|_{L^r_TL^p}\le Ce^{C\overline{V}(T)}\Big\{
2^{-j}\|a_j^0\|_{L^p}+2^{-\f 2rj}\|d_j^0\|_{L^p}\nonumber\\&\qquad+\overline{V}(T)^\f12\bigl(
2^{-j}\|a_j\|_{L^r_TL^p}+\|d_j\|_{L^r_TL^p}\bigr)+2^{-j}
\|F_j\|_{L^1_TL^p}+2^{-\f 2rj}\|G_j\|_{L^1_TL^p}
\nonumber\\&\qquad\qquad\qquad+c(j)2^{-(s-1+\f 2r)j}\overline{V}(T)
\big(\|a\|_{\widetilde{L}^\infty_T\dot{\cB}^{s_p-1,s}_{2,p}}
+\|d\|_{\widetilde{L}^\infty_T\dot{\cB}^{s_p-1,s-1}_{2,p}}\big)\Big\}.
\end{align}
Here we used $2^N\sim e^{C\overline{V}(T)}$ and the summation is finite(at most $2N+1$) for fixed $j$.

\vspace{0.2cm}
{\bf Case 2.} $2^k\le R_0$.
Noting that
\begin{align*}&2^k\ge 2^{j-N}\ge R_02^{-N}\sim R_0e^{-C\overline{V}(T)}.
\end{align*}
Then we apply Proposition \ref{Prop:Green-Lpest} (b) with $b\sim e^{-C\overline{V}(T)}$ to get
\begin{align*}
\|\Delta_k\ba_j(t)\|_{L^p}&+\|\Delta_k\bd_j(t)\|_{L^p}\le C
e^{C\overline{V}(T)}e^{-c2^{2k} t}\bigl(\|\Delta_k a_j^0\|_{L^p}+\|\Delta_k d_j^0\|_{L^p}\bigr)\\
&+Ce^{C\overline{V}(T)}\int_0^te^{-c2^{2k} (t-\tau)}\bigl(
\|\Delta_k\bar{F}_j(\tau)\|_{L^p}+\|\Delta_k\bar{G}_j(\tau)\|_{L^p}\\
&\qquad\qquad\qquad+\|\Delta_k\bar{\mathcal {R}}_j(\tau)\|_{L^p}
+\|\Delta_k\bar{\mathcal {Q}}_j)(\tau)\|_{L^p}\bigr)d\tau.
\end{align*}
Then taking $L^r$ norm with respect to $t$, we get by Young's inequality that
\begin{align*}
&\|\Delta_k\ba_j\|_{L^r_TL^p}+\|\Delta_k\bd_j\|_{L^r_TL^p}\\&\le C
e^{C\overline{V}(T)}\Big\{2^{-\f2rk}\big(\|\Delta_k a_j^0\|_{L^p}+\|\Delta_k d_j^0\|_{L^p}
+\|\Delta_k\bar{F}_j\|_{L^1_TL^p}+\|\Delta_k\bar{\mathcal {R}}_j^1\|_{L^1_TL^p}\\
&\qquad+\|\Delta_k\bar{G}_j\|_{L^1_TL^p}+\|\Delta_k\bar{\mathcal {Q}}_j^1\|_{L^1_TL^p}\big)
+2^{-2k}\bigl(\|\Delta_k\bar{\mathcal {R}}_j^2\|_{L^r_TL^p}+\|\Delta_k\bar{\mathcal {Q}}_j^2\|_{L^r_TL^p}\bigr)\Big\}.
\end{align*}
Substituting \eqref{equ:FG-Lancor}-\eqref{equ:Q_2-Lancor} into
the above inequality, then summing over $k$, we obtain
\begin{align}\label{equ:a-hlest}
&\sum_{|k-j|\le N,2^k\le R_0}\|\Delta_k\ba_j\|_{L^r_TL^p}+
\|\Delta_k\bd_j\|_{L^r_TL^p}
\le Ce^{C\overline{V}(T)}\Bigl\{2^{-\f 2rj}
\big(\|a_j^0\|_{L^p}+\|d_j^0\|_{L^p}\big)\nonumber\\&\qquad\quad+2^{-\f 2rj}\big(\|F_j\|_{L^1_TL^p}+\|G_j\|_{L^1_TL^p}\big)
+\overline{V}(T)^\f12\bigl(2^{-j}\|a_j\|_{L^r_TL^p}+\|d_j\|_{L^r_TL^p}\bigr)
\nonumber\\&\qquad\qquad
+c(j)\overline{V}(T)\bigl(2^{-(s+\f 2r)j}
\|a\|_{\widetilde{L}^\infty_T\dot{\cB}^{s_p-1,s}_{2,p}}
+2^{-(s-1+\f 2r) j}\|d\|_{\widetilde{L}^\infty_T\dot{\cB}^{s_p-1,s-1}_{2,p}}\bigr)
\Big\}.
\end{align}

Summing up \eqref{equ:a-hhest}$-$\eqref{equ:a-hlest} and noting that for $2^k\le R_0$,
\begin{align*}
R_0<2^j\le2^{k+N}\le R_02^{N}\sim R_0e^{C\overline{V}(T)},
\end{align*}
we deduce that for any $1\le r\le \infty$,
\begin{align*}
&\sum_{|k-j|\le N}2^{js}\|\Delta_k\ba_j\|_{L^r_TL^p}+2^{j(s-1+\f 2r)}\|\Delta_k\bd_j\|_{L^r_TL^p}\\
&\le Ce^{C\overline{V}(T)}\Bigl\{2^{js}\big(\|a_j^0\|_{L^p}+2^{-j}\|d_j^0\|_{L^p}\big)
\\&\qquad+\overline{V}(T)^\f12\big(2^{j(s-1+\f 2r)}\|d_j\|_{L^r_TL^p}
+2^{js}\|a_j\|_{L^r_TL^p\cap L^1_TL^p}+2^{j(s+1)}\|d_j\|_{L^1_TL^p}
\big)\\&\qquad\qquad+2^{js}\|F_j\|_{L^1_TL^p}
+2^{j(s-1)}\|G_j\|_{L^1_TL^p}+c(j)\overline{V}(T)\|(a,d)\|_{\cE^s_T}\Big\},
\end{align*}
which together with (\ref{equ:linsy-ahigh}) and (\ref{equ:linsy-dhigh}) implies that
\begin{align}\label{equ:(a,d)-Highest}
&\sum_{2^j> R_0}\bigl(2^{js}\|a_j\|_{L^\infty_TL^p}+2^{js}\|a_j\|_{L^1_TL^p}
+2^{j(s-1)}\|d_j\|_{L^\infty_TL^p}+2^{j(s+1)}\|d_j\|_{L^1_TL^p}\bigr)\nonumber\\
&\le Ce^{C\overline{V}(T)}\Big\{\sum_{2^j> R_0}
2^{js}\big(\|a_j^0\|_{L^p}+2^{-j}\|d_j^0\|_{L^p}\big)+
\|F\|_{ L^1_T\dot{\cB}^{s_p-1, s}_{2,p}}+\|G\|_{ L^1_T\dot{\cB}^{s_p-1,s-1}_{2,p}}\nonumber\\
&\qquad\qquad\qquad\qquad\qquad+
\big(\overline{V}(T)+\overline{V}(T)^\f12\big)\|(a,d)\|_{\cE^{s}_T}\Big\}.
\end{align}

\subsection{Proof of Theorem \ref{thm:hyper-prabolic estimate}}

Thanks to (\ref{equ:(a,d)-Lowest}) and (\ref{equ:(a,d)-Highest}), we obtain
\begin{align*}\label{}
\|a\|&_{\widetilde{L}^\infty_T\dot{\cB}^{s_p-1, s}_{2,p}\cap  L^1_T\dot{\cB}^{s_p+1, s}_{2,p}}
+\|d\|_{\widetilde{L}^\infty_T\dot{\cB}^{s_p-1,s-1}_{2,p}\cap  L^1_T\dot{\cB}^{s_p+1,s+1}_{2,p}}
\nonumber\\ \le&
Ce^{C\overline{V}(T)}\Big\{\|a_0\|_{\dot{\cB}^{s_p-1, s}_{2,p}}
+\|d_0\|_{\dot{\cB}^{s_p-1,s-1}_{2,p}}
+\|F\|_{ L^1_T\dot{\cB}^{s_p-1, s}_{2,p}}
+\|G\|_{ L^1_T\dot{\cB}^{s_p-1,s-1}_{2,p}}\nonumber\\&
+\|v\|_{\widetilde{L}^{\tilde{p}}_T\dot\cB^{\f np,\f {2n}p-\f n2}_{2,p}}
\big(\|a\|_{\widetilde{L}^{\tilde{p}'}_T\dot\cB^{s-\f np+\f n{p'},s}_{2,p}}+
\|d\|_{\widetilde{L}^{\tilde{p}'}_T\dot\cB^{s-\f np+\f n{p'},s_p}_{2,p}}\big)\nonumber\\&\quad+
\|v\|_{\widetilde{L}^2_T\dot \cB^{\f n2,\f n p}_{2,p}}\big(\|a\|_{\widetilde{L}^2_T\dot \cB^{s_p,s}_{2,p}}+
\|d\|_{\widetilde{L}^2_T\dot \cB^{s_p,s}_{2,p}}\big)+\big(\overline{V}(T)+\overline{V}(T)^\f12\big)\|(a,d)\|_{\cE^{s}_T}\Big\}.
\end{align*}
This  together with the following interpolations inequalities
\begin{align*}
&\|v\|_{\widetilde{L}^2_T\dot \cB^{\f n2,\f np}_{2,p}}\le
\|v\|_{\widetilde{L}^{\infty}_T\dot \cB^{\f n2-1,\f {n}p-1}_{2,p}}^{\f{1}2}
\|v\|_{ L^1_T\dot \cB^{\f n2+1,\f np+1}_{2,p}}^{\f{1}{2}},\\&
\|v\|_{\widetilde{L}^{\tilde{p}}_T\dot\cB^{\f np,\f {2n}p-\f n2}_{2,p}}\le
\|v\|_{\widetilde{L}^{\infty}_T\dot\cB^{\f n2-1,\f {n}p-1}_{2,p}}^{\f1{\tilde{p}'}}
\|v\|_{{L}^{1}_T\dot\cB^{\f n2+1,\f {n}p+1}_{2,p}}^{\f1{\tilde{p}}},\\
&
\|d\|_{\widetilde{L}^{\tilde{p}'}_T\dot\cB^{s-\f np+\f n{p'},s_p}_{2,p}}\le
\|d\|_{\widetilde{L}^{\infty}_T\dot\cB^{s_p-1,s-1}_{2,p}}^{\f 1{\tilde{p}}}
\|d\|_{{L}^{1}_T\dot\cB^{s_p+1,s+1}_{2,p}}^{\f 1{\tilde{p}'}},\\
&\|d\|_{\widetilde{L}^2_T\dot \cB^{s_p,s}_{2,p}}\le
\|d\|_{\widetilde{L}^{\infty}_T\dot\cB^{s_p-1,s-1}_{2,p}}^{\f 1{2}}
\|d\|_{{L}^{1}_T\dot\cB^{s_p+1,s+1}_{2,p}}^{\f 1{2}},\\
&\|a\|_{\widetilde{L}^{\tilde{p}'}_T\dot\cB^{s-\f np+\f n{p'},s}_{2,p}}\le
\|a\|_{\widetilde{L}^{\infty}_T\dot\cB^{s_p-1,s}_{2,p}}^{\f 1{\tilde{p}}}
\|a\|_{{L}^{1}_T\dot\cB^{s_p+1,s}_{2,p}}^{\f 1{\tilde{p}'}},\\
&\|a\|_{\widetilde{L}^2_T\dot \cB^{s_p,s}_{2,p}}\le
\|a\|_{\widetilde{L}^{\infty}_T\dot\cB^{s_p-1,s}_{2,p}}^{\f 1{2}}
\|a\|_{{L}^{1}_T\dot\cB^{s_p+1,s}_{2,p}}^{\f 1{2}},
\end{align*}
and Young's inequality yield  Theorem \ref{thm:hyper-prabolic estimate}. \ef

\subsection{Proof of Theorem \ref{thm:hyper-prabolic estimate-2}}
Since the proof is essentially the same as that of Theorem \ref{thm:hyper-prabolic estimate},
here we only present main differences. We use Proposition \ref{Prop:Product} (c) to obtain
\begin{align*}\label{}
&\sum_{j\in\Z}2^{j(s-1)}\|\Delta_j(v\na a)\|_{L^1_TL^2}\le C\bigl(
\|v\|_{\widetilde{L}^2_T\dot\cB^{\f n2,\f {n}p}_{2,p}}
\|a\|_{\widetilde{L}^2_T\dot B^{s}_{2,1}}+\|a\|_{\widetilde{L}^2_T\dot\cB^{\f n2,\f {n}p}_{2,p}}
\|v\|_{\widetilde{L}^2_T\dot B^{s}_{2,1}}\bigr),\\
&\sum_{j\in\Z}2^{j(s-1)}\|\Delta_j(v\na d)\|_{L^1_TL^2}
\le C\bigl(
\|v\|_{\widetilde{L}^2_T\dot\cB^{\f n2,\f {n}p}_{2,p}}
\|d\|_{\widetilde{L}^2_T\dot B^{s}_{2,1}}+\|d\|_{\widetilde{L}^2_T\dot\cB^{\f n2,\f {n}p}_{2,p}}
\|v\|_{\widetilde{L}^2_T\dot B^{s}_{2,1}}\bigr).
\end{align*}
Then in the case of low frequency, we can conclude that
\begin{align}\label{equ:(a,d)-high}
&\sum_{2^j\le R_0}2^{j(s-1+\f2r)}(\|a_j\|_{L^r_TL^2}+\|d_j\|_{L^r_TL^2})\nonumber\\
&\le C\Big(\sum_{2^j\le R_0}
2^{j(s-1)}(\|a_j^0\|_{L^2}+\|d_j^0\|_{L^2})+\|v\|_{\widetilde{L}^2_T\dot \cB^{\f n 2,\f n p}_{2,p}}
\|(a,d)\|_{\widetilde{L}^2_T\dot B^{s}_{2,1}}
\nonumber\\&\qquad\quad+\|v\|_{\widetilde{L}^2_T\dot B^{s}_{2,1}}\|(a,d)\|_{\widetilde{L}^2_T\dot\cB^{\f n2,\f {n}p}_{2,p}}+
\|F\|_{L^1_T\dot \cB^{s-1,s}}+\|G\|_{L^1_T\dot B^{s-1}_{2,1}}\Big).
\end{align}

From  Proposition \ref{Prop:Commutator} and Lemma \ref{Lem:com-est}, we infer that
\begin{align*}
&\big\|\Delta_k\bar{\mathcal {R}}^1_j\big\|_{L^1_TL^2}
\le Cc(j)2^{-s j}e^{\overline{V}(T)}\overline{V}(T)\|a\|_{\widetilde{L}^\infty_T\dot{B}^{s}_{2,1}},\\
&\big\|\Delta_k\bar{\mathcal {Q}}^1_j\big\|_{L^1_TL^2}
\le Cc(j)
2^{-(s-1)j}e^{\overline{V}(T)}\overline{V}(T)\|d\|_{\widetilde{L}^\infty_T\dot{B}^{s-1}_{2,1}}.
\end{align*}
Then in the case of high frequency, we can conclude that
\begin{align}
&\sum_{2^j> R_0}\bigl(2^{js}\|a_j\|_{L^\infty_TL^2}+2^{js}\|a_j\|_{L^1_TL^2}
+2^{j(s-1)}\|d_j\|_{L^\infty_TL^2}+2^{j(s+1)}\|d_j\|_{L^1_TL^2}\bigr)\nonumber\\
&\le Ce^{C\overline{V}(T)}\Big\{\sum_{2^j> R_0}
2^{js}\big(\|a_j^0\|_{L^2}+2^{-j}\|d_j^0\|_{L^2}\big)+
\|F\|_{{L}^1_T\dot{\cB}^{s-1,s}}+\|G\|_{ L^1_T\dot{B}^{s-1}_{2,1}}\nonumber\\&
\qquad\qquad\qquad\qquad+
\big(\overline{V}(T)+\overline{V}(T)^\f12\big)\|(a,d)\|_{E^{s}_T}\Big\}.\nonumber
\end{align}
This together with (\ref{equ:(a,d)-high}) implies Theorem \ref{thm:hyper-prabolic estimate-2}.\ef

\section{Proof of theorem 1.1}

This section is devoted to the proof of Theorem 1.1.
We denote
\begin{align*}
a(t,x)\eqdefa\f{\rho(\varpi^{-2}t,\varpi^{-1} x)}{\bar{\rho}}-1,
\quad v(t,x)\eqdefa\varpi^{-1} u(\varpi^{-2}t,\varpi^{-1}x),
\end{align*}
where $\varpi\eqdefa (P'(\bar{\rho}))^{\f12}$. Then the system \eqref{equ:cNS} can be rewritten as
\ben\label{equ:cNS-anotherform}
\left\{
\begin{array}{ll}
\p_ta+v\cdot\na a+\dv v=-a\dv v,\\
\p_tv+v\cdot\na v-\mathcal {A}v+\na a=-L(a)\mathcal {A}v-K(a)\na a,\\
(a,v)|_{t=0}=(a_0,v_0),
\end{array}
\right.
\een
where $$\mathcal {A}=\bar{\mu}\Delta+(\bar{\lambda}+\bar{\mu})\na \dv,
\quad K(a)=\f{P'(\bar{\rho}(1+a))}{(1+a)P'(\bar{\rho})}-1\quad\hbox{and}\quad L(a)=\frac{a}{1+a},$$
with $\bar{\mu}=\f{\mu}{\bar{\rho}}$
and $\bar{\lambda}=\f{\lambda}{\bar{\rho}}$. We denote
$$
d\eqdefa \Lambda^{-1}\dv v,\quad \Omega\eqdefa \Lambda^{-1}(\textrm{curl} v)_j^i.
$$
In view of the equality $$\Delta v^i=\p_i\dv v+\p_j(\textrm{curl}\, v)_j^i, \quad\hbox{where}\,\,
(\textrm{curl} v)_j^i=\p_j v^i-\p_i v^j,$$
we find that $(a, d,\Omega)$ satisfies
\begin{equation}\label{equ:cNStrans}
\left\{
\begin{array}{ll}
\p_ta+\Lambda d+v\cdot\na a=F,\\
\p_td-\bar{\nu}\Delta d-\Lambda a+v\cdot\na d=G, \\
\p_t\Omega-\bar{\mu}\Delta \Omega=H,\\
(a,d,\Omega)|_{t=0}=(a_0,\Lambda^{-1}\dv v_0,\Lambda^{-1}\textrm{curl} v_0),
\end{array}
\right.
\end{equation}
where $\bar\nu=2\bar\mu+\bar\lambda$ and
\begin{align*}
&F=-a\dv v,\\
&G=v\cdot\na d-\Lambda^{-1}\dv\big(v\cdot\na v+L(a)
\mathcal{A}v+K(a)\na a\big),\\
&H=-\Lambda^{-1}\textrm{curl}\big(v\cdot\na v+L(a)\mathcal{A}v\big),\\
&v=-\Lambda^{-1}\na d-\Lambda^{-1}\dv \Omega.
\end{align*}
\subsection{A priori estimates}

\bthm{Proposition}\label{Prop:Prioriestimate-v1} Let $2\le p<2n$, $p\le\min\big(
4,\f {2n}{n-2}\big)$.
Assume that $(a,d,\Omega)$ is a smooth solution of the system (\ref{equ:cNStrans})
on $[0,T]$ with
\begin{align*}
\|a\|_{L^\infty([0,T]\times \mathbb{R}^n)}\le \f 12.
\end{align*}
Then we have
\begin{align}\label{equ:Prioriestimate}
\|(a,v)\|_{\cE^{\f np}_T}\le C&e^{C\|(a,v)\|_{\cE^{\f np}_T}}
\Big\{\|(a_0,v_0)\|_{\cE^{\f np}_0}
+\|(a,v)\|_{\cE^{\f np}_T}^\f 32\bigl(1+\|(a,v)\|_{\cE^{\f np}_T}\bigr)^{\f n2+2}
\Big\}.
\end{align}
Here
$
\|(a_0,v_0)\|_{\cE_0^{\f n p}}\eqdefa \|a_0\|_{\dot \cB^{ \frac{n}{2}-1, \f n p}_{2,p}}+
\|v_0\|_{\dot \cB^{\frac{n}{2}-1, \frac{n}{p}-1}_{2,p}}.
$
\ethm

\noindent{\bf Proof.}\, Thanks to Theorem \ref{thm:hyper-prabolic estimate}, we get
\begin{align*}
\|(a,d)\|_{\cE^{\f np}_T}\le C&e^{C\|(a,v)\|_{\cE^{\f np}_T}}
\Big\{\|(a_0,d_0)\|_{\cE_0^{\f n p}}
+\big(\|(a,v)\|_{\cE^{\f np}_T}^\f12+
\|(a,v)\|_{\cE^{\f np}_T}
\big)\nonumber\\&\times
\|(a,d)\|_{\cE^{\f np}_T}+
\|F\|_{ L^1_T\dot{\cB}^{\frac n 2-1, \frac n p}_{2,p}}+\|G\|_{ L^1_T\dot{\cB}^{\frac n 2-1,\frac n p-1}_{2,p}}
\Big\},
\end{align*}
and from Proposition \ref{Prop:Parabolic-est},
\begin{align*}
\|\Omega\|_{\widetilde{L}^\infty_T\dot \cB^{\frac{n}{2}-1, \frac{n}{p}-1}_{2,p}}+
\|\Omega\|_{ L^1_T\dot \cB^{\frac{n}{2}+1, \frac{n}{p}+1}_{2,p}}
\le C(\|\Omega_0\|_{\dot \cB^{\frac{n}{2}-1, \frac{n}{p}-1}_{2,p}}+\|H\|_{ L^1_T\dot{\cB}^{\frac n 2-1, \frac n p}_{2,p}}).
\end{align*}
Then we obtain
\begin{align}\label{equ:Priori-est}
\|(a,v)\|_{\cE^{\f np}_T}&\le Ce^{C\|(a,v)\|_{\cE^{\f np}_T}}
\Big\{\|(a_0,v_0)\|_{\cE_0^{\f n p}}
+\big(\|(a,v)\|_{\cE^{\f np}_T}^\f12+\|(a,v)\|_{\cE^{\f np}_T}\big)\nonumber\\&\qquad\qquad\times
\|(a,v)\|_{\cE^{\f np}_T}+
\|F\|_{ L^1_T\dot{\cB}^{\frac n 2-1, \frac n p}_{2,p}}+\|(G,H)\|_{ L^1_T\dot{\cB}^{\frac n 2-1,\frac n p-1}_{2,p}}\Big\}.
\end{align}
Here we used the fact
\beno
\|v\|_{\widetilde{L}^r_T\dot{\cB}^{s,\sigma}_{2,p}}\approx
\|d\|_{\widetilde{L}^r_T\dot{\cB}^{s,\sigma}_{2,p}}+
\|\Omega\|_{\widetilde{L}^r_T\dot{\cB}^{s,\sigma}_{2,p}}
\eeno
for any $s,\sigma\in \R, r\in [1,\infty]$.

Now we need to estimate the nonlinear terms $F, G$ and $H$ in (\ref{equ:Priori-est}).
We define $\widetilde{p}$ as $\f 1{\tilde{p}}={\f n{2p}-\f n4+\f12}$,
and denote $\widetilde{p}'$ by the conjugate index of $\widetilde{p}$.
Due to $2\le p<2n$, $p\le \f {2n} {n-2}$, using Proposition \ref{Prop:Product} (b)
with $s=\f np, t=\f n{p'}-1, \widetilde{s}=\f n2-1, \widetilde{t}=\f n2, \ga=0 $ yields
\begin{align*}
&\sum_{2^j\le R_0}2^{j(\f n2-1)}\|\Delta_jF\|_{L^1_TL^2}\\&\le C
\|a\|_{\widetilde{L}^{\tilde{p}}_T\dot\cB^{\f np,\f {2n}p-\f n2}_{2,p}}
\|\dv v\|_{\widetilde{L}^{\tilde{p}'}_T\dot\cB^{\f n{p'}-1,\f n2-1}_{2,p}}+C
\|\dv v\|_{\widetilde{L}^2_T\dot\cB^{\f n2-1,\f np-1}_{2,p}}
\|a\|_{\widetilde{L}^2_T\dot\cB^{\f n2,\f np}_{2,p}}\nonumber\\
&\le C\|a\|_{\widetilde{L}^{\tilde{p}}_T\dot\cB^{\f np,\f {n}p}_{2,p}}
\|v\|_{\widetilde{L}^{\tilde{p}'}_T\dot\cB^{\f n{p'},\f n2}_{2,p}}+C
\|v\|_{\widetilde{L}^2_T\dot\cB^{\f n2,\f np}_{2,p}}
\|a\|_{\widetilde{L}^2_T\dot\cB^{\f n2,\f np}_{2,p}},
\end{align*}
and from Proposition \ref{Prop:Product} (a) with $\sigma=\tau=\f np,$
\begin{align*}
&\sum_{2^j>R_0}2^{j\f np}\|\Delta_jF\|_{L^1_TL^p}
\le C\|a\|_{\widetilde{L}^\infty_T\dot\cB^{\f n2,\f np}_{2,p}}
\|\dv v\|_{ L^1_T\dot\cB^{\f n2,\f np}_{2,p}}\le C\|a\|_{\widetilde{L}^\infty_T\dot\cB^{\f n2-1,\f np}_{2,p}}
\|v\|_{ L^1_T\dot\cB^{\f n2+1,\f np+1}_{2,p}},
\end{align*}
from which, we infer that
\begin{align}\label{equ:Nonest-F}
\|F\|_{ L^1_T\dot{\cB}^{\frac n 2-1, \frac n p}_{2,p}}\le C\Bigl(&\|a\|_{\widetilde{L}^{\tilde{p}}_T\dot\cB^{\f np,\f {n}p}_{2,p}}
\|v\|_{\widetilde{L}^{\tilde{p}'}_T\dot\cB^{\f n{p'},\f n2}_{2,p}}+
\|v\|_{\widetilde{L}^2_T\dot\cB^{\f n2,\f np}_{2,p}}
\|a\|_{\widetilde{L}^2_T\dot\cB^{\f n2,\f np}_{2,p}}\nonumber\\&\quad+
\|a\|_{\widetilde{L}^\infty_T\dot\cB^{\f n2-1,\f np}_{2,p}}
\|v\|_{ L^1_T\dot\cB^{\f n2+1,\f np+1}_{2,p}}\Bigr).
\end{align}
Similarly, we have
\begin{align*}
&\|v\na v\|_{ L^1_T\dot{\cB}^{\frac n 2-1, \frac n p-1}_{2,p}}+
\|v\na d\|_{ L^1_T\dot{\cB}^{\frac n 2-1, \frac n p-1}_{2,p}}\nonumber\\
&\le C\|v\|_{\widetilde{L}^{\tilde{p}}_T\dot\cB^{\f np,\f {2n}p-\f n2}_{2,p}}
\|v\|_{\widetilde{L}^{\tilde{p}'}_T\dot\cB^{\f n{p'},\f n2}_{2,p}}+
C\|v\|_{\widetilde{L}^2_T\dot\cB^{\f n2,\f np}_{2,p}}\|v\|_{\widetilde{L}^2_T\dot\cB^{\f n2,\f np}_{2,p}}.
\end{align*}
From Proposition \ref{Prop:Product} (b) with $\gamma=\f np-\f n2$ and Proposition \ref{Prop:composition},
it follows that
\begin{align*}
&\sum_{2^j\le R_0}2^{j(\f n2-1)}\|\Delta_j(K(a)\na
a)\|_{L^1_TL^2}\nonumber\\ &\le  C\|K(a)\|_{\widetilde{L}^{\tilde{p}}_T
\dot\cB^{\f np,\f {2n}p-\f n2}_{2,p}} \|\na a\|_{\widetilde{L}^{\tilde{p}'}_T\dot\cB^{\f n{p'}-1,\f
np-1}_{2,p}}+C\|\na a\|_{\widetilde{L}^2_T\dot\cB^{\f n2-1,\f
np-1}_{2,p}}\|K(a)\|_{\widetilde{L}^2_T\dot\cB^{\f n2,\f np}_{2,p}}
\nonumber\\&\le C
\big(\|a\|_{\widetilde{L}^{\tilde{p}}_T\dot\cB^{\f np,\f {n}p}_{2,p}}\|a\|_{\widetilde{L}^{\tilde{p}'}_T
\dot\cB^{\f n{p'},\f np}_{2,p}}+\|a\|_{\widetilde{L}^2_T\dot\cB^{\f n2,\f
np}_{2,p}}\|a\|_{\widetilde{L}^2_T\dot\cB^{\f n2,\f np}_{2,p}}\big)
\bigl(1+\|a\|_{\widetilde{L}^\infty_T\dot\cB^{\f np,\f np}_{2,p}}\bigr)^{\f n 2+1} ,
\end{align*}
and with $\gamma=-1$
\begin{align*}
&\sum_{2^j\le R_0}2^{j(\f n2-1)}\|\Delta_j(L(a)\mathcal {A} v)\|_{L^1_TL^2}\\
&\le C\|L(a)\|_{\widetilde{L}^{\tilde{p}}_T\dot\cB^{\f np,\f {2n}p-\f n2}_{2,p}}
\|\cA v\|_{\widetilde{L}^{\tilde{p}'}_T\dot\cB^{\f n{p'}-1,\f n2-1+\gamma}_{2,p}}
+C\|\cA v\|_{ L^1_T\dot\cB^{\f n2-1,\f np-1}_{2,p}}
\|L(a)\|_{\widetilde{L}^\infty_T\dot\cB^{\f n2,\f np}_{2,p}}\nonumber\\ &\le C
\big(\|a\|_{\widetilde{L}^{\tilde{p}}_T\dot\cB^{\f np,\f np}_{2,p}}
\|v\|_{\widetilde{L}^{\tilde{p}'}_T\dot\cB^{\f n{p'},\f n2}_{2,p}}
+\|a\|_{\widetilde{L}^\infty_T\dot\cB^{\f n2,\f np}_{2,p}}
\|v\|_{ L^1_T\dot\cB^{\f n2+1,\f np+1}_{2,p}}\big)\bigl(1+
\|a\|_{\widetilde{L}^\infty_T\dot\cB^{\f np,\f np}_{2,p}}\bigr)^{\f n 2+1}.
\end{align*}
On the other hand, from Proposition \ref{Prop:Product} (a)  and Proposition \ref{Prop:composition},
\begin{align*}
&\sum_{2^j>R_0}2^{j(\f np-1)}(\|\Delta_j(K(a)\na a)\|_{L^1_TL^p}
+\|\Delta_j(L(a)\mathcal {A} v)\|_{L^1_TL^p})\\
&\le C\|K(a)\|_{\widetilde{L}^2_T\dot\cB^{\f n2,\f np}_{2,p}}
\|\na a\|_{\widetilde{L}^2_T\dot\cB^{\f n2-1,\f np-1}_{2,p}}
+C\|L(a)\|_{\widetilde{L}^\infty_T\dot\cB^{\f n2,\f np}_{2,p}}
\|A v\|_{ L^1_T\dot\cB^{\f n2-1,\f np-1}_{2,p}}\\
&\le C\big(1
+\|a\|_{\widetilde{L}^\infty_T\dot\cB^{\f np,\f np}_{2,p}}\big)^{\f n 2+1}\big(\|a\|_{\widetilde{L}^2_T\dot\cB^{\f n2,\f np}_{2,p}}
\|a\|_{\widetilde{L}^2_T\dot\cB^{\f n2,\f np}_{2,p}}+\|a\|_{\widetilde{L}^\infty_T\dot\cB^{\f n2,\f np}_{2,p}}
\|v\|_{ L^1_T\dot\cB^{\f n2+1,\f np+1}_{2,p}}\big)
.
\end{align*}
Thus, we deduce that
\begin{align}\label{equ:Nonest-GH}
&\|(G,H)\|_{ L^1_T\dot{\cB}^{\frac n 2-1,\frac n p-1}_{2,p}}\le
C\bigl(1+\|a\|_{\widetilde{L}^\infty_T\dot\cB^{\f np,\f np}_{2,p}}\bigr)^{\f n 2+1}
\Bigl(\|a\|_{\widetilde{L}^2_T\dot\cB^{\f n2,\f
np}_{2,p}}\|a\|_{\widetilde{L}^2_T\dot\cB^{\f n2,\f np}_{2,p}}\nonumber\\&\quad+\|a\|_{\widetilde{L}^{\tilde{p}}_T\dot\cB^{\f np,\f {n}p}_{2,p}}\bigl(\|a\|_{\widetilde{L}^{\tilde{p}'}_T
\dot\cB^{\f n{p'},\f np}_{2,p}}+\|v\|_{\widetilde{L}^{\tilde{p}'}_T\dot\cB^{\f n{p'},\f n2}_{2,p}}\bigr)
+\|a\|_{\widetilde{L}^\infty_T\dot\cB^{\f n2,\f np}_{2,p}}\|v\|_{ L^1_T\dot\cB^{\f n2+1,\f np+1}_{2,p}}\Bigr).
\end{align}
Noting that the interpolation inequalities
\begin{align*}
&\|a\|_{\widetilde{L}^{\tilde{p}}_T\dot\cB^{\f np,\f {n}p}_{2,p}}\le
\|a\|_{\widetilde{L}^{\infty}_T\dot\cB^{\f n2-1,\f {n}p}_{2,p}}^{\f1{\tilde{p}'}}
\|a\|_{{L}^{1}_T\dot\cB^{\f n2+1,\f {n}p}_{2,p}}^{\f1{\tilde{p}}},\\
&\|a\|_{\widetilde{L}^{\tilde{p}'}_T\dot\cB^{\f n{p'},\f np}_{2,p}}\le
\|a\|_{\widetilde{L}^{\infty}_T\dot\cB^{\f n{2}-1,\f np}_{2,p}}^{\f1{\tilde{p}}}
\|a\|_{{L}^{1}_T\dot\cB^{\f n{2}+1,\f np}_{2,p}}^{\f1{\tilde{p}'}},\\&
\|a\|_{\widetilde{L}^2_T\dot\cB^{\f n2,\f np}_{2,p}}\le
\|a\|_{\widetilde{L}^\infty_T\dot\cB^{\f n2-1,\f np}_{2,p}}^{\f12}
\|a\|_{ L^1_T\dot\cB^{\f n2+1,\f np}_{2,p}}^{\f12},\\
&\|v\|_{\widetilde{L}^{\tilde{p}'}_T\dot\cB^{\f n{p'},\f n2}_{2,p}}\le
\|v\|_{\widetilde{L}^{\infty}_T\dot\cB^{\f n{2}-1,\f np-1}_{2,p}}^{\f1{\tilde{p}}}
\|v\|_{{L}^{1}_T\dot\cB^{\f n{2}+1,\f np+1}_{2,p}}^{\f1{\tilde{p}'}},\\&
\|v\|_{\widetilde{L}^{\tilde{p}}_T\dot\cB^{\f np,\f {2n}p-\f n2}_{2,p}}\le
\|v\|_{\widetilde{L}^{\infty}_T\dot\cB^{\f n2-1,\f {n}p-1}_{2,p}}^{\f1{\tilde{p}'}}
\|v\|_{{L}^{1}_T\dot\cB^{\f n2+1,\f {n}p+1}_{2,p}}^{\f1{\tilde{p}}},\\
&
\|v\|_{\widetilde{L}^2_T\dot\cB^{\f n2,\f np}_{2,p}}\le
\|v\|_{\widetilde{L}^\infty_T\dot\cB^{\f n2-1,\f np-1}_{2,p}}^{\f12}
\|v\|_{ L^1_T\dot\cB^{\f n2+1,\f np+1}_{2,p}}^{\f12},
\end{align*}
and from $p\le \f {2n} {n-2}$,
\begin{align*}
\|a\|_{\widetilde{L}^\infty_T\dot\cB^{\f np,\f np}_{2,p}}
\le C\|a\|_{\widetilde{L}^\infty_T\dot\cB^{\f n2-1,\f np}_{2,p}},
\end{align*}
we plug (\ref{equ:Nonest-F})-(\ref{equ:Nonest-GH}) into \eqref{equ:Priori-est} to obtain (\ref{equ:Prioriestimate}).\ef

\bthm{Proposition}\label{Prop:Prioriestimate-v2} Under the assumption of Proposition \ref{Prop:Prioriestimate-v1},
 we have
\begin{align}\label{equ:Prioriestimate2}
\|(a,v)\|_{E^{\f n2}_T}\le Ce^{C\|(a,v)\|_{\cE^{\f np}_T}}
\Big\{\|(a_0,v_0)\|_{E^{\f n2}_0}
+&\|(a,v)\|_{E^{\f n2}_T}\|(a,v)\|_{\cE^{\f np}_T}^\f12\nonumber
\\&\quad\times\bigl(1+\|(a,v)\|_{\cE^{\f np}_T}\bigr)^{\f n2+2}
\Big\}.
\end{align}
Here
$
\|(a_0,v_0)\|_{E_0^{\f n 2}}\eqdefa \|a_0\|_{\dot \cB^{\frac{n}{2}-1, \f n 2}}+
\|v_0\|_{\dot B^{\frac{n}{2}-1}_{2,1}}.
$
\ethm

\noindent{\bf Proof.}\,\,Since the proof is similar to that of Proposition \ref{Prop:Prioriestimate-v1}, here we only indicate main differences.
Firstly, we get by using Theorem \ref{thm:hyper-prabolic estimate-2} and Proposition \ref{Prop:Parabolic-est} that
\begin{align}\label{equ:Priori-est2}
\|(a,v)\|_{E^{\f n2}_T}\le Ce^{C\|(a,v)\|_{\cE^{\f np}_T}}
\Big\{&\|(a_0,v_0)\|_{E_0^{\f n 2}}+\big(\|(a,v)\|_{\cE^{\f np}_T}^\f12+\|(a,v)\|_{\cE^{\f np}_T}\big)
\|(a,v)\|_{E^{\f n2}_T}\nonumber\\&+\|F\|_{L^1_T{\dot\cB}^{\f n2-1, \f n2}}+\|(G,H)\|_{L^1_T\dot{B}^{\f n2-1}_{2,1}}
\Big\}.
\end{align}
Applying Proposition \ref{Prop:Product} (c) yields
\begin{align*}
&\sum_{2^j \le R_0}2^{j(\f n2-1)}\|\Delta_jF\|_{L^1_TL^2}\le
C\|a\|_{L^\infty_T\dot\cB^{\f n2-1,\f {n}p-1}_{2,p}}
\|\dv v\|_{L^1_T\dot B^{\f n2}_{2,1}}+
C\|\dv v\|_{\widetilde{L}^2_T\dot\cB^{\f n2-1,\f np-1}_{2,p}}
\|a\|_{\widetilde{L}^2_T\dot B^{\f n2}_{2,1}},\nonumber\\
&\sum_{2^j> R_0}2^{j\f n2}\|\Delta_jF\|_{L^1_TL^2}
\le C\|a\|_{\widetilde{L}^\infty_T\dot\cB^{\f n2,\f np}_{2,p}}
\|\dv v\|_{L^1_T\dot B^{\f n2}_{2,1}}+C\|\dv v\|_{L^1_T\dot\cB^{\f n2,\f np}_{2,p}}
\|a\|_{\widetilde{L}^\infty_T\dot B^{\f n2}_{2,1}},
\end{align*}
which imply by the interpolation  and Lemma \ref{Lem:Besov-properties} (a) that
\begin{align}\label{equ:F-est3}
\|F\|_{L^1_T\dot{\cB}^{\frac n 2-1,\f n2}}&\le C\|(a,v)\|_{\cE^{\f np}_T}\|(a,v)\|_{E^{\f n2}_T}.
\end{align}
Similarly, we have
\begin{align*}
&\sum_{j\in\mathbb{Z}}2^{j(\f n2-1)}(\|\Delta_j(v\na d)\|_{L^1_TL^2}+\|\Delta_j(v\na v)\|_{L^1_TL^2})
\nonumber\\&\le C\|v\|_{\widetilde{L}^\infty_T\dot\cB^{\f n2-1,\f {n}p-1}_{2,p}}
\|v\|_{L^1_T\dot B^{\f n2+1}_{2,1}}+C
\|v\|_{L^1_T\dot\cB^{\f n2+1,\f np+1}_{2,p}}\|v\|_{\widetilde{L}^\infty_T\dot B^{\f n2-1}_{2,1}}.
\end{align*}
From Proposition \ref{Prop:Product} (c) and Proposition \ref{Prop:composition}, we infer that
\begin{align*}\label{}
&\sum_{j\in\mathbb{Z}}2^{j(\f n2-1)}\|\Delta_j(K(a)\na a)\|_{L^1_TL^2}\\
&\le C\|K(a)\|_{\widetilde{L}^2_T\dot\cB^{\f n2,\f np}_{2,p}}
\|\na a\|_{\widetilde{L}^2_T\dot B^{\f n2-1}_{2,1}}+C\|\na a\|_{\widetilde{L}^2_T\dot \cB^{\f n2-1,\f np-1}_{2,p}}
\|K(a)\|_{\widetilde{L}^2_T\dot B^{\f n2}_{2,1}}\\
&\le C\bigl(1+\|a\|_{\widetilde{L}^\infty_T\dot\cB^{\f np,\f np}_{2,p}}\bigr)^{\f n 2+1}\|a\|_{\widetilde{L}^2_T\dot\cB^{\f n2,\f np}_{2,p}}
\|a\|_{\widetilde{L}^2_T\dot B^{\f n2 }_{2,1}},
\end{align*}and
\begin{align*}
&\sum_{j\in\mathbb{Z}}2^{j(\f n2-1)}\|\Delta_j(L(a)\mathcal {A} v)\|_{L^1_TL^2}\\
&\le C\|L(a)\|_{\widetilde{L}^\infty_T\dot\cB^{\f n2,\f np}_{2,p}}
\|\cA v\|_{L^1_T\dot B^{\f n2-1}_{2,1}}+C\|\cA v\|_{L^1_T\dot \cB^{\f n2-1,\f np-1}_{2,p}}
\|L(a)\|_{\widetilde{L}^\infty_T\dot B^{\f n2 }_{2,1}}\\ &\le  C\bigl(1+\|a\|_{\widetilde{L}^\infty_T\dot\cB^{\f np,\f np}_{2,p}}\bigr)^{\f n 2+1}
\|(a,v)\|_{\cE^{\f np}_T}\|(a,v)\|_{E^{\f n2}_T}.
\end{align*}
Thus, we obtain
\begin{align*}
\|G\|_{{L}^1_T\dot{B}^{\frac n 2-1}_{2,1}}+
\|H\|_{{L}^1_T\dot{B}^{\frac n 2-1}_{2,1}}\le C\bigl(1+\|a\|_{\widetilde{L}^\infty_T\dot\cB^{\f np,\f np}_{2,p}}\bigr)^{\f n 2+1}
\|(a,v)\|_{\cE^{\f np}_T}\|(a,v)\|_{E^{\f n2}_T}.
\end{align*}
This  together with (\ref{equ:Priori-est2}) and (\ref{equ:F-est3}) implies (\ref{equ:Prioriestimate2}).\ef

\subsection{Approximate solutions and uniform estimates}

The construction of approximate solutions is based on the following local existence theorem.

\bthm{Theorem} \cite{Danchin-cpde07}\label{thm:local} Assume that $\rho_0-\bar \rho\in \dot B^{\f n 2}_{2,1}$ and $u_0\in \dot B^{\f n 2-1}_{2,1}$
with $\rho_0$ bounded away from zero. There exists a positive time $T$ such that
the system (\ref{equ:cNS}) has a unique solution $(\rho, u)$ with $\rho$ bounded away from zero,
\begin{align*}
\rho-\bar\rho\in C([0,T);\dot B^{\f n 2}_{2,1})\quad \textrm{and} \quad
u\in C([0,T);\dot B^{\f n 2-1}_{2,1})\cap L^1(0,T;\dot B^{\f n 2+1}_{2,1}).
\end{align*}
Moreover, the solution $(\rho,u)$ can be continued beyond $T$ if the following three conditions hold:

(a) The function $\rho-\bar \rho$ belongs to $L^\infty(0,T;\dot B^{\f n2}_{2,1})$;

(b) The function $\rho$ is bounded away from zero;

(c) $\int_0^T\|\na u(\tau)\|_{L^\infty}d\tau<\infty$.

\ethm

To apply Theorem \ref{thm:local} , we need
\bthm{Lemma}\label{Lem:approximate} Let $p\ge2$.
For any $(\rho_0-\bar{\rho},u_0)\in \dot \cB^{\f {n} 2-1,\f np}_{2,p}\times
\bigl(\dot \cB^{\f {n} 2-1,\f n p-1}_{2,p}\bigr)^n$ satisfying $\rho_0\ge \bar\rho$, there exists a sequence
$\{(\rho_{0}^k,u_{0}^k)\}_{k\in \N}$ with
$(\rho_{0}^k-\bar \rho, u_{0}^k)\in \dot \cB^{\f n 2-1,\f n2}\times \bigl(\dot B^{\f n 2-1}_{2,1}\bigr)^n$
such that
\begin{align}\label{equ:convergence}
\|\rho_{0}^k-\rho_0\|_{\dot \cB^{\f {n} 2-1,\f np}_{2,p}}\longrightarrow 0, \quad
\|u_{0}^k-u_0\|_{\dot \cB^{\f {n} 2-1,\f np-1}_{2,p}}\longrightarrow 0
\end{align}
as $k\rightarrow \infty$, and $\rho_{0,k}\ge \f {\bar\rho} 2$ for any $k\in \N$.
\ethm

\noindent{\bf Proof.}\, This lemma can be proved by following the proof of Lemma 4.2 in \cite{Abidi}.\ef

Let $(\rho_{0}^k, u_{0}^k)$ be as stated in Lemma \ref{Lem:approximate}. Then Theorem \ref{thm:local} ensures that
there exists a maximal existence time $T_k>0$ such that the system (\ref{equ:cNS}) with the initial data $(\rho_{0,k}, u_{0,k})$ has
a unique solution $(\rho^k, u^k)$ with $\rho^k$ bounded away from zero,
\begin{align*}
\rho^k-\bar\rho\in C([0,T_k);\dot B^{\f n 2}_{2,1})\quad \textrm{and} \quad
u^k\in C([0,T_k);\dot B^{\f n 2-1}_{2,1})\cap L^1(0,T_k;\dot B^{\f n 2+1}_{2,1}).
\end{align*}
Using the definition of the Besov space and Lemma \ref{Lem:Bernstein}, it is easy to check that
\begin{align*}
\rho^k-\bar\rho\in C([0,T_k);\dot \cB^{\f {n} 2-1,\f np}_{2,p})\quad \textrm{and} \quad
u^k\in C([0,T_k);\dot \cB^{\f {n} 2-1,\f np-1}_{2,p})\cap L^1(0,T_k;\dot \cB^{\f {n} 2+1,\f np+1}_{2,p}).
\end{align*}

We set
\begin{align*}
a^k(t,x)=\f{\rho^k(\varpi^{-2}t,\varpi^{-1} x)}{\bar{\rho}}-1,
\quad v^k(t,x)=\varpi^{-1} u^k(\varpi^{-2}t,\varpi^{-1} x).
\end{align*}
Then  from (\ref{equ:initial condition}) and (\ref{equ:convergence}), we find
\begin{align*}
\|(a_0^k,v_0^k)\|_{\cE^{\f np}_0}\le C_0\eta,
\end{align*}
for some constant $C_0$.
Given a constant $M$ to be chosen later on,
let us define
$$T^*_k\eqdefa \sup\Big\{t\in[0,T_k); \,\,\|(a^k,v^k)\|_{\cE_t^{\f np}}\le M\eta\Big\}.$$

Firstly, we claim that
\begin{align*}
T^*_k=T_k,\qquad \forall k\in\mathbb{N}.
\end{align*}
Using the continuity argument, it suffices to show that for all $k\in\mathbb{N}$,
\begin{align}\label{equ:uniform-est}
\|(a^k,v^k)\|_{\cE^{\f np}_{T^*_k}}\le \f 34M\eta.
\end{align}
In fact, noticing that $\|a^k\|_{L^\infty}\le C_1\|a^k\|_{\dot \cB^{\f n2-1,\f np}_{2,p}}$, we choose $\eta$ as
\begin{align*}
M\eta\le \f1{2C_1}
\end{align*}
such that
\begin{align*}
\|a^k\|_{L^\infty([0,T^*_k)\times \mathbb{R}^n)}\le \f 12.
\end{align*}
Thus, we can apply Proposition \ref{Prop:Prioriestimate-v1} to obtain
\begin{align}\label{equ:appro-est1}
\|(a^{k},v^{k})&\|_{\cE^{\f np}_{T^*_k}}\le Ce^{CM\eta}
\Big\{C_0\eta
+(M\eta)^\f32(1+M\eta)^{\f n2+2}
\Big\}.
\end{align}
Set $M=4CC_0$, then choose $\eta$ small enough such that
\begin{align*}
e^{CM\eta}\le \frac32,\quad \quad C(M\eta)^\f12(1+M\eta)^{\f n2+2}\le \f 14,
\end{align*}
then the inequality (\ref{equ:uniform-est}) follows from (\ref{equ:appro-est1}). In conclusion,
we construct a sequence of approximate solution $(\rho^k,u^k)$ of (\ref{equ:cNS}) on $[0,T_k)$ satisfying
\ben\label{equ:uniform-est1}
\|(a^k,v^k)\|_{\cE^{\f np}_{T_k}}\le M\eta,
\een
for any $k\in \N$. Next, we claim that
\begin{align*}
T_k=+\infty,\quad \forall\, k\in\mathbb{N}.
\end{align*}
Thanks to Theorem \ref{thm:local} and (\ref{equ:uniform-est1}), it remains to prove $a^k\in L^\infty(0,T_k;\dot B^{\f n2}_{2,1})$.
While thanks to Proposition \ref{Prop:Prioriestimate-v2} and (\ref{equ:uniform-est1}), we have
\begin{align*}
\|(a^k,v^k)\|_{E^{\f n2}_{T_k}}\le C&e^{CM\eta}
\Big\{\|(a_0^k,v_0^k)\|_{E^{\f n2}_0}
+\|(a^k,v^k)\|_{E^{\f n2}_{T_k}}(M\eta)^\f12(1+M\eta)^{\f n2+2}
\Big\},
\end{align*}
which implies that
$$
\|(a^k,v^k)\|_{E^{\f n2}_{T_k}}\le C\|(a_0^k,v_0^k)\|_{E^{\f n2}_0}.
$$

\subsection{Existence}

We will use the compact argument to prove the existence of the solution. Due to (\ref{equ:uniform-est1}),
it is easy to verify that
\begin{align*}
&\bullet\,\, a^k\textrm{ is uniformly bounded in }L^\infty(0,\infty; \dot B^{\f np}_{p,1});\\
&\bullet\,\, v^k\textrm{ is uniformly bounded in }L^\infty(0,\infty; \dot B^{\f np-1}_{p,1})\cap L^1(0,\infty;\dot B^{\f np+1}_{p,1}).
\end{align*}
By the interpolation, we also have
\begin{align*}
v^k\textrm{ is uniformly bounded in }L^{\f {2}{1-\epsilon}}(0,\infty;\dot B^{\f np-\epsilon}_{p,1}),
\end{align*}
for any $\epsilon\in [-1,1]$.

Let $v^k_L$ be a solution of
\begin{align*}
\p_t v^k_L-\mathcal {A}v^k_L=0,\quad v^k_L(0)=v_{0}^{k}.
\end{align*}
It is easy to check that $v^k_L$ tends to the solution of
\begin{align*}
\p_t v_L-\mathcal {A}v_L=0,\quad v_L(0)=v_{0},
\end{align*}
in $L^\infty(0,\infty; \dot B^{\f np-1}_{p,1})\cap L^1(0,\infty;\dot B^{\f np+1}_{p,1})$.

We set $\widetilde{a}^k\eqdefa a^k-a_{0}^k$ and $\widetilde{v}^k\eqdefa v^k-v^k_{L}$.
Firstly, we claim that $(\widetilde{a}^k, \widetilde{v}^k)$ is uniformly bounded in $C^{\f12}_{loc}(\R^+;\dot B^{\f n p-1}_{p,1})\times
C^{\f{2-\varsigma}{2}}_{loc}(\R^+;\dot B^{\f n p-1-\varsigma}_{p,1})$ with
$$\varsigma=\min\Big(\f{2n}p-1, 1\Big).$$
Recall that
$$\p_t\widetilde{a}^k=-v^k\cdot\na a^k-\dv v^k-a^k\dv v^k,$$
which together with Lemma \ref{Lem:Product} implies $\p_t\widetilde{a}^{k}\in L^2\dot B^{\f n p-1}_{p,1}$, thus
$\widetilde{a}^k$ is uniformly bounded in $C^{\f12}(\R^+; \dot B^{\f n p-1}_{p,1})$. On the other hand,
$$
\p_t\widetilde{v}^k=-v^k\cdot\na v^k+\mathcal {A}\widetilde{v}^k-\na a^k-L(a^k)\mathcal {A}v^k-K(a^k)\na a^k.
$$
Thanks to Lemma \ref{Lem:Product} and Proposition \ref{Prop:composition}, we have
\begin{align*}\|&v^k\cdot\na v^{k}+L(a^k)\cA v^k)\|_{L^{\f {2}{2-\varsigma}}\dot B^{\f np-1-\varsigma}_{p,1}}
\\&\le C\|v^k\|_{L^2\dot B^{\f np}_{p,1}}\|\na v^k\|_{L^{\f {2}{1-\varsigma}}\dot B^{\f np-\varsigma-1}_{p,1}}
+C(\|a^k\|_{L^\infty\dot B^{\f np}_{p,1}})
\|\cA v^k\|_{L^{\f {2}{2-\varsigma}}\dot B^{\f np-1-\varsigma}_{p,1}}.
\end{align*}
The inclusion map
$\dot\cB^{\f n2-1,\f np}_{2,p}\hookrightarrow \dot B^{\f np-\varsigma}_{p,1}$ implies
$\Lambda a^{k}\in L^\infty\dot B^{\f np-1-\varsigma}_{p,1}$, thus $K(a^k)\na a^k$ is bounded in
$L^\infty\dot B^{\f np-1-\varsigma}_{p,1}$. Therefore, we have
$\p_t \widetilde{v}^k\in L^{\f {2}{2-\varsigma}}_{loc}\dot B^{\f np-1-\varsigma}_{p}$, thus
$\widetilde{v}^k$ is uniformly bounded in $C^{\f{2-\varsigma}{2}}_{loc}(\R^+;\dot B^{\f n p-1-\varsigma}_{p,1})$.
\vspace{.25cm}

Let $\{\chi_j\}_{j\in \N}$ be a sequence of smooth functions supported in the ball $B(0,j+1)$
and equal to 1 on $B(0,j)$. The
claim ensures that for any $j\in \N$, $\{\chi_j \widetilde{a}^k\}_{k\in
\N}$ is uniformly bounded in $C^\f12_{loc}(\R^+;\dot B^{
\frac{n}{p}-1}_{p,1})$, and $\{\chi_j\widetilde{v}^k\}_{k\in \N}$ is uniformly
bounded in $C^{\f{2-\varsigma}{2}}_{loc}(\R^+;\dot B^{\f n p-1-\varsigma}_{p,1})$. Observe that for
any $\chi\in C_0^\infty(\mathbb{R}^n)$, the map: $(\widetilde{a}^k, \widetilde{v}^k)\mapsto  (\chi \widetilde{a}^k, \chi \widetilde{v}^k)$ is
compact from $$\big(\dot B^{\f n p-1}_{p,1}\cap \dot B^{\f n p}_{p,1}\big)
\times\big(\dot B^{\f n p-1-\varsigma}_{p,1}\cap \dot B^{\f n p-1}_{p,1}\big)\quad\hbox{into}
\quad \dot B^{\f n p-1}_{p,1}\times \dot B^{\f n p-1-\varsigma}_{p,1}.$$
By applying Ascoli's theorem and Cantor's diagonal process ,
there exists some function $(\widetilde{a},\widetilde{v})\in \cE^{\f n p}$ such that for
any $j\in\N$, \ben\label{eq:conver1}
\begin{split}
&\chi_j \widetilde{a}^k\longrightarrow \chi_j \widetilde{a} \quad \textrm{in}\quad C_{loc}(\R^+;\dot B^{\f n p-1}_{p,1}),\\
&\chi_j \widetilde{v}^k\longrightarrow \chi_j \widetilde{v} \quad \textrm{in}\quad C_{loc}(\R^+;\dot B^{\f n p-1-\varsigma}_{p,1}),
\end{split}
\een
as $k$ tends to $\infty$(up to a subsequence). By the interpolation, we also have
\begin{align}\label{eq:conver2}
\begin{split}
&\chi_j a^k\longrightarrow \chi_j \widetilde{a} \quad \textrm{in}\quad C_{loc}(\R^+;\dot B^{ \frac{N}{p}-s}_{p,1}),\quad \forall\,0<s\le 1,\\
&\chi_j v^k\longrightarrow \chi_j \widetilde{v} \quad \textrm{in}\quad
L^1_{loc}(\R^{+};\dot B^{ \frac{N}{p}+s}_{p,1}),\quad \forall\,-1 \le
s<1.
\end{split}
\end{align}

With (\ref{eq:conver1})-(\ref{eq:conver2}), it is a routine process to
verify that $(\widetilde{a}+a_0,\widetilde{v}+v_L)$ satisfies the system (\ref{equ:cNS-anotherform}) in the
sense of distribution(see also \cite{Danchin-inven}). Finally, following
the argument in \cite{Danchin-inven}, we can show that $(a,v)\in
C([0,\infty);\dot \cB^{\frac{n}{2}-1,\f n p}_{2,p})\times C([0,\infty);\dot \cB^{\frac{n}{2}-1,
\frac{n}{p}-1}_{2,p})$.

\subsection{Uniqueness}
In this subsection, we prove the uniqueness of the solution. Assume that
$(\rho^1,u^1)\in \cE^{\f n p}$ and $(\rho^2,u^2)\in \cE^{\f np}$ are  two solutions of the system (\ref{equ:cNS})
with the same initial data. Without loss of generality, we may assume that
 $(\rho^1,u^1)$ satisfies
\begin{align}
\|(\rho^1-\bar{\rho},u^1)\|_{\cE^{\f np}}\le M\eta.
\end{align}
Since $\rho^2-\bar\rho\in C([0,T];\dot \cB^{\f n2-1,
\frac{N}{p}}_{2,p})$ and $\rho^2(0,x)\ge \f {\bar \rho} 2$, there exists a
positive time $T$ such that
\begin{align*}
\rho^2(t,x)\ge \f {\bar \rho} 3,\quad \textrm{for } (t,x)\in [0,T]\times \R^n.
\end{align*}

We set
\begin{align*}
&a^k(t,x)=\f{\rho^k(\varpi^{-2}t,\varpi^{-1} x)}{\bar{\rho}}-1,
\quad v^k(t,x)=\varpi^{-1} u^k(\varpi^{-2}t,\varpi^{-1} x),\quad k=1,2,\\
&\delta a=a^1-a^2,\quad \delta v=v^1-v^2.
\end{align*}
Thanks to (\ref{equ:cNS-anotherform}), we find that $(\delta a,\delta u)$ satisfies
\begin{align}\label{eq:linearcNS-diffe} \left\{
\begin{array}{ll}
\p_t\delta a+v^2\cdot\na \delta a=\delta F,\\
\p_t\delta v-\mathcal {A}\delta v=\delta G, \\
(\delta a,\delta v)|_{t=0}=(0,0),
\end{array}
\right. \end{align}
where
\begin{align*}
\delta F=&-\delta v\cdot\na a^1-\dv \delta v-a^1\dv \delta v-\delta a\dv v^2,\\
\delta G=&-\na \delta a-(v^1\cdot \na v^1-v^2\cdot\na v^2)-L(a^1)\cA v^1+L(a^2)\cA v^2\\
&-K(a^1)\na a^1+K(a^2)\na a^2.
\end{align*}

In what follows, we denote $V^i(t)=\int_0^t\|v^i(\tau)\|_{\dot{B}^{
\frac{n}{p}+1}_{p,1}}d\tau$ for $i=1,2$, and denote by $A_T$ a
constant depending  on $\|a^1\|_{\widetilde{L}^\infty_T(\dot B^{
\frac{n}{p}}_{p,1})}$ and $\|a^2\|_{\widetilde{L}^\infty_T(\dot B^{
\frac{n}{p}}_{p,1})}$. Due to the inclusion relation $\cE^{\f n p}\subseteq
\cE^{1}(p\le n)$, it suffices to prove the uniqueness of the solution in
$\cE^1$. So, in the following we take $p=n$.

Applying Proposition \ref{Prop:transport} gives
\begin{align}\label{eq:dens-diffe2}
\|\delta a(t)\|_{\dot{B}^{0}_{p,\infty}}\le e^{CV^2(t)}\int_0^t \|\delta F(\tau)\|_{\dot{B}^{0}_{p,\infty}}d\tau.
\end{align}
By Lemma \ref{Lem:Product}, we have
\begin{align*}
\|\delta F(\tau)\|_{\dot{B}^{0}_{p,\infty}}
\le C\|v^2\|_{\dot{B}^{2}_{p,1}}\|\delta a\|_{\dot{B}^{0}_{p,\infty}}+C(1+\|a^1\|_{\dot{B}^{1}_{p,1}})\|\delta v\|_{\dot{B}^{1}_{p,1}}.
\end{align*}
Plugging it into (\ref{eq:dens-diffe2}), we get by Gronwall's inequality that
\begin{align}\label{eq:dens-diffe3}
\|\delta a(t)\|_{\dot{B}^{0}_{p,\infty}}\le e^{CV^2(t)}\int_0^t(1+\|a^1\|_{\dot{B}^{1}_{p,1}})\|\delta v\|_{\dot{B}^{1}_{p,1}}d\tau.
\end{align}
Applying Proposition \ref{Prop:momentum} to the second equation of \eqref{eq:linearcNS-diffe} gives
\begin{align}\label{eq:velo-diff2}
\|\delta v(t)\|_{\widetilde{L}^1_t(\dot B^{1}_{p,\infty})}+\|\delta v(t)\|_{\widetilde{L}^2_t(\dot B^{0}_{p,\infty})}
\le C\int_0^t \|\delta G(\tau)\|_{\dot B^{-1}_{p,\infty}}d\tau.
\end{align}
From Lemma \ref{Lem:Product} and  Proposition \ref{Prop:composition}, we infer that
\begin{align}
\|\delta G(t)\|_{\dot B^{-1}_{p,\infty}}
\le& C\|(v^1,v^2)\|_{\dot B^{1}_{p,1}}\|\delta v\|_{\dot B^{0}_{p,\infty}}
+A_t\|a^1\|_{\dot B^{1}_{p,1}}
\|\delta v\|_{\dot B^{1}_{p,\infty}}\nonumber\\
&+A_t(1+\|v^2\|_{\dot B^{2}_{p,1}})
\|\delta a\|_{\dot B^{0}_{p,\infty}}.\label{eq:G-diff1}
\end{align}
We take $T$ small enough such that
\begin{align*}
\|(v^1,v^2)\|_{\widetilde{L}^1_T(\dot B^{2}_{p,1})\cap \widetilde{L}^2_T(\dot B^{1}_{p,1})} \ll 1.
\end{align*}
Thus, plugging (\ref{eq:G-diff1}) into (\ref{eq:velo-diff2}),
we infer that for any $t\in [0,T]$,
\begin{align}\label{eq:velo-diff3}
\|\delta v\|_{\widetilde{L}^1_t(\dot B^{1}_{p,\infty})}
\le A_T\int_0^t\bigl(1+\|(v^1,v^2)\|_{\dot B^{2}_{p,1}}\bigr)\|\delta a\|_{\dot B^{0}_{p,\infty}}d\tau.
\end{align}
\bthm{Lemma}\cite{Danchin-PRSE}\label{Lem:loginequ} Let $s\in \R$. Then for any $1\le p,r\le+\infty$ and $0<\epsilon\le 1$,
we have
\begin{align*}
\|f\|_{\widetilde{L}^r_T(\dot{B}^s_{p,1})}\le
C\frac{\|f\|_{\widetilde{L}^r_T(\dot{B}^s_{p,\infty})}}{\epsilon}
\log\Bigl(e+\frac{\|f\|_{\widetilde{L}^r_T(\dot{B}^{s-\epsilon}_{p,\infty})}
+\|f\|_{\widetilde{L}^r_T(\dot{B}^{s+\epsilon}_{p,\infty})}}
{\|f\|_{\widetilde{L}^r_T(\dot{B}^s_{p,\infty})}}\Bigr).
\end{align*}
\ethm
From Lemma \ref{Lem:loginequ}, it follows that
\begin{align*}
\|\delta v\|_{L^1_t(\dot B^{1}_{p,1})}\le C\|\delta v\|_{\widetilde{L}^1_t(\dot{B}^1_{p,\infty})}
\log\Bigl(e+\frac{\|\delta v\|_{\widetilde{L}^1_t(\dot{B}^{0}_{p,\infty})}
+\|\delta v\|_{\widetilde{L}^1_t(\dot{B}^{2}_{p,\infty})}}
{\|\delta v\|_{\widetilde{L}^1_t(\dot{B}^1_{p,\infty})}}\Bigr),
\end{align*}
which together with (\ref{eq:dens-diffe3}) and (\ref{eq:velo-diff3}) yields that for any $t\in [0,{T}]$,
\begin{align*}
\|\delta v\|_{\widetilde{L}^1_t(\dot B^{1}_{p,\infty})}
\le e^{CV^2(t)}A_T\int_0^t\bigl(1+\|(v^1,v^2)\|_{\dot B^{2}_{p,1}}\bigr)
\|\delta v\|_{\widetilde{L}^1_\tau(\dot{B}^1_{p,\infty})}
\log\bigl(e+C_T\|\delta v\|_{\widetilde{L}^1_\tau(\dot B^{1}_{p,\infty})}^{-1}\bigr)d\tau,
\end{align*}
where $C_T=\|\delta v\|_{\widetilde{L}^1_T(\dot{B}^{0}_{p,\infty})}
+\|\delta v\|_{\widetilde{L}^1_T(\dot{B}^{2}_{p,\infty})}$.
Noticing that $\|(v^1,v^2)(t)\|_{\dot B^{2}_{p,1}}$ is integrable on $[0,T]$ and
\begin{align*}
\int_0^1\f {dr} {r\log(r+C_Tr^{-1})}dr=+\infty,
\end{align*}
Osgood lemma applied concludes that $(\delta a,\delta v)=0$ on $[0,T]$, and
a continuity argument ensures that $(a^1,v^1)=(a^2,v^2)$ on $[0,\infty)$.

\section{Appendix}

In what follows, we denote $\chi_{\{\cdot\}}$ by the characteristic function defined in $\Z$, and $\{c(j)\}_{j\in \Z}$ by
a sequence in $\ell^1$ with the norm $\|\{c(j)\}\|_{\ell^1}\le 1$.

\bthm{Lemma}\label{Lem:Paraproduct} Let $s, \sigma, t,\tau\in\mathbb{R}$, $2\le p\le 4$,
$p'$ is the conjugate index of $p$,
and $1\le r,r_1,r_2\le \infty$ with $\f 1 r=\f 1 {r_1}+\f 1 {r_2}$.
Then there hold\vspace{0.1cm}

(a) if $s\le\f n2,\,\sigma\le\f np,$
then for $2^j>R_0$,
\begin{align}\label{equ:para-high}
&\|\Delta_j(T_fg)\|_{L_T^{r}L^p}\nonumber\\
&\le Cc(j)
\big(2^{j(\f n{p'}-s-t)}+2^{j(\f n{2}-s-\tau)}+
2^{j(\f n{p}-\sigma-\tau)}\big)
\|f\|_{\widetilde{L}_T^{r_1}\dot\cB^{s,\sigma}_{2,p}}\|g\|_{\widetilde{L}_T^{r_2}\dot\cB^{t,\tau}_{2,p}};
\end{align}

(b) if $ s\le\f np$, $\sigma\le\f{2n}{p}-\f n2$, then for $2^j\le R_0$,
\begin{align}\label{equ:para-low}
\|\Delta_j(T_fg)\|_{L_T^{r}L^2}
\le Cc(j)&
\big(2^{j(\f n2-s-t)}+\chi_{\{2^j\sim R_0\}}
(2^{j(\f np-s-\tau)}+2^{j(\f{2n}p-\f n2-\sigma-\tau)})\big)
\nonumber\\ &\times\|f\|_{\widetilde{L}_T^{r_1}\dot\cB^{s,\sigma}_{2,p}}\|g\|_{\widetilde{L}_T^{r_2}\dot\cB^{t,\tau}_{2,p}};
\end{align}

(c) if $s\le\f n2,\,\sigma\le\f np,$ then
\begin{align}
\|\Delta_j(T_fg)\|_{L_T^{r}L^2}\le Cc(j)\big(2^{j(\f {n}{2}-s-t)
}+2^{j(\f np-\sigma-t)}\big)\|f\|_{\widetilde{L}^{r_1}_T\dot\cB^{s,\sigma}_{2,p}}
\|g\|_{\widetilde{L}^{r_2}_T\dot B^{t}_{2,1}}.\label{equ:para-high-L2}
\end{align}

\ethm

\noindent{\bf Proof.} (a) Thanks to \eqref{orth}, we have
\begin{align*}
\Delta_j(T_fg)=\sum_{|k-j|\le 4}\Delta_j(S_{k-1}f\Delta_{k}g)
=\sum_{|k-j|\le 4}\sum_{k'\le k-2}\Delta_j(\Delta_{k'}f\Delta_{k}g).
\end{align*}
We denote $J\eqdefa \big\{(k',k);\,|k-j|\le 4, k'\le k-2\big\}$,
then for $2^j>R_0$,
\begin{align*}
\|\Delta_j(T_fg)\|_{L_T^r L^p}&\le\sum_{J}\|\Delta_j(\Delta_{k'}f\Delta_{k}g)\|_{L^r_TL^p}\nonumber\\
&\le\Big(\sum_{J_{\ell m}}+\sum_{J_{\ell h}}+\sum_{J_{hh}}\Big)
\|\Delta_j(\Delta_{k'}f\Delta_{k}g)\|_{L_T^r L^p}\\
&\eqdefa I_1+I_2+I_3,
\end{align*}
where
\begin{align*}
&J_{\ell m}=\{(k',k)\in J,\,\, 2^{k'}\le R_0, 2^{-4}R_0\le2^k\le R_0\},\\&
J_{\ell h}=\{(k',k)\in J,\,\, 2^{k'}\le R_0, 2^k>R_0\},\\
&J_{hh}=\{(k',k)\in J,\,\, 2^{k'}>R_0, 2^k>R_0\}.\end{align*}

We get by using Lemma \ref{Lem:Bernstein} and $s\le \f n2$ that
\begin{align*}
I_1&\le C\sum_{(k',k)\in J_{\ell m}}\!\!
\|\Delta_{k'}f\|_{L_T^{r_1}L^\infty}2^{k(\f n2-\f np)}\|\Delta_{k}g\|_{L_T^{r_2}L^2}\nonumber\\&\le C
\sum_{(k',k)\in J_{\ell m}}2^{k's}\|\Delta_{k'}f\|_{L_T^{r_1}L^2}
2^{k'(\f n2-s)}2^{kt}\|\Delta_{k}g\|_{L_T^{r_2}L^2}
2^{k(\f n2-\f np-t)}\\
&\le Cc(j)2^{j(\f n{p'}-s-t)}\|f\|_{\widetilde{L}_T^{r_1}\dot\cB^{s,\sigma}_{2,p}}\|g\|_{\widetilde{L}_T^{r_2}\dot\cB^{t,\tau}_{2,p}}.
\end{align*}
Similarly, we have
\begin{align*}
I_2&\le \sum_{(k',k)\in J_{\ell h}}\|\Delta_{k'}f\|_{L_T^{r_1}L^\infty}\|\Delta_{k}g\|_{L_T^{r_2}L^p}\\
&\le C\sum_{(k',k)\in J_{\ell h}}2^{k's}\|\Delta_{k'}f\|_{L_T^{r_1}L^2}
2^{k'(\f n2-s)}2^{k\tau}
\|\Delta_{k}g\|_{L_T^{r_2}L^p}2^{-k\tau}\nonumber\\&
\le Cc(j)2^{j(\f n{2}-s-\tau)}\|f\|_{\widetilde{L}_T^{r_1}\dot\cB^{s,\sigma}_{2,p}}
\|g\|_{\widetilde{L}_T^{r_2}\dot\cB^{t,\tau}_{2,p}},
\end{align*}
and noting $\sigma\le \f n p$,
\begin{align*}
I_3&\le \sum_{(k',k)\in J_{hh}}\|\Delta_{k'}f\|_{L_T^{r_1}L^\infty}\|\Delta_{k}g\|_{L_T^{r_2}L^p}\\
&\le C\sum_{(k',k)\in J_{hh}}2^{k'\sigma}\|\Delta_{k'}f\|_{L_T^{r_1}L^p}
2^{k'(\f np-\sigma)}2^{k\tau}
\|\Delta_{k}g\|_{L_T^{r_2}L^p}2^{-k\tau}\nonumber\\&
\le Cc(j)2^{j(\f n{p}-\sigma-\tau)}
\|f\|_{\widetilde{L}_T^{r_1}\dot\cB^{s,\sigma}_{2,p}}\|g\|_{\widetilde{L}_T^{r_2}\dot\cB^{t,\tau}_{2,p}}.
\end{align*}
Then the inequality \eqref{equ:para-high} can be deduced from the estimates of $I_1, I_2$ and $I_3$.

(b) We denote
$$J\eqdefa\big\{(k',k);\,|k-j|\le 4, k'\le k-2\big\},$$
then for $2^j\le R_0$, there holds
$$
J=J_{\ell \ell}\cup J_{\ell m}\cup J_{hm},
$$
where
\begin{align*}&J_{\ell \ell}=\{(k',k)\in J,\,\, 2^{k'}\le R_0, 2^k\le R_0\},\\
&{J}_{\ell m}=\{(k',k)\in J,\,\, 2^{k'}\le R_0, R_0<2^k\le2^4R_0\},
\\& J_{hm}=\{(k',k)\in J,\,\, 2^{k'}>R_0, R_0<2^k\le2^4R_0\}.\end{align*}
Then thanks to (\ref{orth}), we have for $2^j\le R_0$,
\begin{align*}
\|\Delta_j(T_fg)\|_{L_T^{r}L^2}&\le
\sum_{J}\|\Delta_j(\Delta_{k'}f\Delta_{k}g)\|_{L_T^{r}L^2}
\nonumber\\&\le \Big(\sum_{J_{\ell \ell}}+\sum_{{J}_{\ell m}}+\sum_{J_{hm}}\Big)
\|\Delta_j(\Delta_{k'}f\Delta_{k}g)\|_{L_T^{r}L^2}\\
&\eqdefa II_1+II_2+II_3.
\end{align*}

We get by using Lemma \ref{Lem:Bernstein} and $s\le \f n p$ that
\begin{align*}
II_1&\le C\sum_{(k',k)\in J_{\ell \ell}}2^{k's}\|\Delta_{k'}f\|_{L^{r_1}_TL^2}2^{k'(\f n2-s)}
2^{kt}\|\Delta_{k}g\|_{L^{r_2}_TL^2}2^{-kt}\\&\le Cc(j)
2^{j(\f n{2}-s-t)}\|f\|_{\widetilde{L}_T^{r_1}\dot\cB^{s,\sigma}_{2,p}}
\|g\|_{\widetilde{L}_T^{r_2}\dot\cB^{t,\tau}_{2,p}},
\nonumber\\
II_2&\le
C\sum_{(k',k)\in {J}_{\ell m}}2^{k's}\|\Delta_{k'}f\|_{L^{r_1}_TL^2}2^{k'(\f np-s)}
2^{k\tau}\|\Delta_{k}g\|_{L^{r_2}_TL^p}2^{-k\tau}\\& \le C c(j)\chi_{\{2^j\sim R_0\}}2^{j(\f n{p}-s-\tau)}
\|f\|_{\widetilde{L}_T^{r_1}\dot\cB^{s,\sigma}_{2,p}}
\|g\|_{\widetilde{L}_T^{r_2}\dot\cB^{t,\tau}_{2,p}},
\end{align*}
and noting  $\sigma\le\f{2n}{p}-\f n2,$
\begin{align*}
II_3&\le C
\sum_{(k',k)\in J_{h m}}2^{k'\sigma}\|\Delta_{k'}f\|_{L^{r_1}_TL^p}2^{k'(\f {2n}p-\f n2-\sigma)}
2^{k\tau}\|\Delta_{k}g\|_{L^{r_2}_TL^p}2^{-k\tau}\nonumber\\&\le Cc(j)
\chi_{\{2^j\sim R_0\}}2^{j(\f {2n}{p}-\f n{2}-\sigma-\tau)}
\|f\|_{\widetilde{L}_T^{r_1}\dot\cB^{s,\sigma}_{2,p}}\|g\|_{\widetilde{L}_T^{r_2}\dot\cB^{t,\tau}_{2,p}}.
\end{align*}
The inequality \eqref{equ:para-low} follows from the above estimates.

(c) Since the proof of (\ref{equ:para-high-L2}) is similar, here we omit it.\ef

\bthm{Lemma}\label{Lem:Remainder}
Let $s, \sigma, t,\tau\in \R$, $2\le p\le4$, $p'$ is the conjugate index of $p$,
and $1\le r,r_1,r_2\le \infty$ with $\f 1 r=\f 1 {r_1}+\f 1 {r_2}$. Assume that
$s+t>0,$ $s+\tau>0,$ $\sigma+t>0,$ and $\sigma+\tau>0$. Then there hold
\begin{align}\label{equ:Remainder-high}
\|\Delta_jR(f,g)\|_{L^r_TL^p}\le& Cc(j)
\bigl(2^{j(\f n{p'}-s-t)}
+2^{j(\f n2-s-\tau)}\nonumber\\&\quad+2^{j(\f n2-\sigma-t)}+2^{j(\f {n}{p}-\sigma-\tau)}\bigr)
\|f\|_{\widetilde{L}_T^{r_1}\dot\cB^{s,\sigma}_{2,p}}\|g\|_{\widetilde{L}_T^{r_2}\dot\cB^{t,\tau}_{2,p}};\\
\|\Delta_jR(f,g)\|_{L^r_TL^2} \le& Cc(j)\bigl(2^{j(\f n2-s-t)}
+2^{j(\f np-s-\tau)}\nonumber\\&
\quad+2^{j(\f np-\sigma-t)}+2^{j(\f {2n}p-\f n2-\sigma-\tau)}\bigr)\|f\|_{\widetilde{L}_T^{r_1}\dot\cB^{s,\sigma}_{2,p}}
\|g\|_{\widetilde{L}_T^{r_2}\dot\cB^{t,\tau}_{2,p}};\label{equ:Remainder-low}\\
\|\Delta_jR(f,g)\|_{L^r_TL^2}\le& Cc(j)\big(2^{j(\f n2-s-t)}+2^{j(\f np-s-\tau)}\big)\|f\|_{\widetilde{L}_T^{r_1}\dot B^{s}_{2,1}}
\|g\|_{\widetilde{L}_T^{r_2}\dot\cB^{t,\tau}_{2,p}}\label{equ:Remainder-high-L2}.
\end{align}
where $p'$ is the conjugate index of $p$.
\ethm

\noindent{\bf Proof.}\,\,Thanks to \eqref{orth}, we have
\begin{align*}
\Delta_j(R(f,g))=\sum_{k\ge
j-3}\sum_{|k'-k|\le1}\Delta_j(\Delta_{k}f{\Delta}_{k'}g), \end{align*}
Set $J\eqdefa \{(k,k');\,k\ge j-3, |k'-k|\le1\}$, then
$$J=J_{\ell\ell}\cup J_{\ell m}\cup J_{hm}\cup J_{hh},$$
where
\begin{align*}
&J_{\ell\ell}=\{(k,k')\in J,\, 2^k\le R_0, 2^{k'}\le R_0\},\\
&J_{\ell m}=\{(k,k')\in J,\, 2^k\le R_0, R_0<2^{k'}\le 2R_0\},\\&
J_{hm}=\{(k,k')\in J,\, 2^k>R_0, 2^{-1}R_0\le2^{k'}\le R_0\},\\
&J_{hh}=\{(k,k')\in J,\, 2^k>R_0, 2^{k'}>R_0\}.\end{align*}
Thus, we have
\begin{align*}
\Delta_j(R(f,g))&=\Bigl(\sum_{J_{\ell\ell}}+\sum_{J_{\ell m}}+\sum_{J_{hm}}+\sum_{J_{hh}}\Bigr)
\Delta_j(\Delta_{k}f{\Delta}_{k'}g)\\
&\eqdefa I_1+I_2+I_3+I_4.
\end{align*}

We get by Lemma  \ref{Lem:Bernstein} and $s+t>0$ that
\begin{align*}
\|I_1\|_{L^r_TL^p}&\le C2^{j\f n{p'}}\!\sum_{(k,k')\in J_{\ell\ell}}\|\Delta_kf\Delta_{k'}g\|_{L^{r}_TL^1}\\
&\le C2^{j\f n{p'}}\!\sum_{(k,k')\in J_{\ell\ell}}\!2^{ks}\|\Delta_kf\|_{L^{r_1}_TL^2}
2^{-ks}2^{k't}\|\Delta_{k'}g\|_{L^{r_2}_TL^2}
2^{-k't}\nonumber\\&\le C
c(j)2^{j(\f n{p'}-s-t)}\|f\|_{\widetilde{L}_T^{r_1}\dot\cB^{s,\sigma}_{2,p}}\|g\|_{\widetilde{L}_T^{r_2}\dot\cB^{t,\tau}_{2,p}},
\end{align*}
and
\begin{align*}
\|I_1\|_{L^{r}_TL^2}&\le C2^{j\f n{2}}\!\sum_{(k,k')\in J_{\ell\ell}}\|\Delta_kf\Delta_{k'}g\|_{L^{r}_TL^1}\\
&\le Cc(j)2^{j(\f n{2}-s-t)}\|f\|_{\widetilde{L}_T^{r_1}\dot\cB^{s,\sigma}_{2,p}}\|g\|_{\widetilde{L}_T^{r_2}\dot\cB^{t,\tau}_{2,p}}.
\end{align*}
Similarly, due to $s+\tau>0$, we obtain
\begin{align*}
\|I_2\|_{L^r_TL^p}&\le C2^{j\f n{2}}\!\sum_{(k,k')\in J_{\ell m}}\|\Delta_kf\Delta_{k'}g\|_{L^{r}_TL^\f{2p}{2+p}}\\
&\le C2^{j\f n{2}}\!\sum_{(k,k')\in J_{\ell m}}\!2^{ks}\|\Delta_kf\|_{L^{r_1}_TL^2}
2^{-ks}2^{k'\tau}\|\Delta_{k'}g\|_{L^{r_2}_TL^p}
2^{-k'\tau}\nonumber\\&\le
Cc(j)2^{j(\f n{2}-s-\tau)}\|f\|_{\widetilde{L}_T^{r_1}\dot\cB^{s,\sigma}_{2,p}}\|g\|_{\widetilde{L}_T^{r_2}\dot\cB^{t,\tau}_{2,p}},
\end{align*}
and
\begin{align*}
\|I_2\|_{L^{r}_TL^2}&\le C2^{j\f n{p}}\!\sum_{(k,k')\in J_{\ell m}}\|\Delta_kf\Delta_{k'}g\|_{L^{r}_T{L^\f{2p}{2+p}}}\\
&\le Cc(j)2^{j(\f n{p}-s-\tau)}\|f\|_{\widetilde{L}_T^{r_1}\dot\cB^{s,\sigma}_{2,p}}\|g\|_{\widetilde{L}_T^{r_2}\dot\cB^{t,\tau}_{2,p}}.
\end{align*}
Thanks to $\sigma+t>0$, we have
\begin{align*}
\|I_3\|_{L^r_TL^p}&\le C2^{j\f n{2}}\!\sum_{(k,k')\in J_{h m}}\|\Delta_kf\Delta_{k'}g\|_{L^{r}_TL^\f{2p}{2+p}}\\
&\le C2^{j\f n{2}}\!\sum_{(k,k')\in J_{hm}}\!2^{k\sigma}\|\Delta_kf\|_{L^{r_1}_TL^p}
2^{-k\sigma}2^{k't}\|\Delta_{k'}g\|_{L^{r_2}_TL^2}
2^{-k't}\nonumber\\&\le
Cc(j)2^{j(\f n{2}-\sigma-t)}\|f\|_{\widetilde{L}_T^{r_1}\dot\cB^{s,\sigma}_{2,p}}\|g\|_{\widetilde{L}_T^{r_2}\dot\cB^{t,\tau}_{2,p}},
\end{align*}
and
\begin{align*}
\|I_3\|_{L^{r}_TL^2}&\le C2^{j\f n{p}}\!\sum_{(k,k')\in J_{hm}}\|\Delta_kf\Delta_{k'}g\|_{L^{r}_T{L^\f{2p}{2+p}}}\\
&\le Cc(j)2^{j(\f n{p}-\sigma-t)}\|f\|_{\widetilde{L}_T^{r_1}\dot\cB^{s,\sigma}_{2,p}}\|g\|_{\widetilde{L}_T^{r_2}\dot\cB^{t,\tau}_{2,p}}.
\end{align*}
Finally, due to $\sigma+\tau>0$ and $2\le p\le 4$, we have
\begin{align*}
\|I_4\|_{L^r_TL^p}&\le C2^{j\f n{p}}\!\sum_{(k,k')\in J_{hh}}\|\Delta_kf\Delta_{k'}g\|_{L^{r}_TL^\f{p}{2}}\\
&\le C2^{j\f n{p}}\!\sum_{(k,k')\in J_{hh}}\!2^{k\sigma}\|\Delta_kf\|_{L^{r_1}_TL^p}
2^{-k\sigma}2^{k'\tau}\|\Delta_{k'}g\|_{L^{r_2}_TL^p}
2^{-k'\tau}\nonumber\\&\le
Cc(j)2^{j(\f n{p}-\sigma-\tau)}\|f\|_{\widetilde{L}_T^{r_1}\dot\cB^{s,\sigma}_{2,p}}\|g\|_{\widetilde{L}_T^{r_2}\dot\cB^{t,\tau}_{2,p}},
\end{align*}
and
\begin{align*}
\|I_4\|_{L^{r}_TL^2}&\le C2^{j(\f{2n}p-\f n2)}\!\sum_{(k,k')\in J_{hh}}\|\Delta_kf\Delta_{k'}g\|_{L^{r}_T{L^\f{p}{2}}}\\
&\le Cc(j)2^{j(\f{2n}p-\f n2-\sigma-\tau)}\|f\|_{\widetilde{L}_T^{r_1}\dot\cB^{s,\sigma}_{2,p}}
\|g\|_{\widetilde{L}_T^{r_2}\dot\cB^{t,\tau}_{2,p}}.
\end{align*}
Then the inequalities (\ref{equ:Remainder-high}) and (\ref{equ:Remainder-low})
can be deduced from the above estimates. Finally, the inequality (\ref{equ:Remainder-high-L2}) can be deduced from
\begin{align*}
\|I_1+I_3\|_{L^{r}_TL^2}&\le C2^{j\f n{2}}\!\sum_{(k,k')\in J_{\ell\ell}\cup J_{hm}}
\|\Delta_kf\Delta_{k'}g\|_{L^{r}_TL^1}\\
&\le Cc(j)2^{j(\f n{2}-s-t)}\|f\|_{\widetilde{L}_T^{r_1}\dot B^{s}_{2,1}}
\|g\|_{\widetilde{L}_T^{r_2}\dot\cB^{t,\tau}_{2,p}},
\end{align*}
and
\begin{align*}
\|I_2+I_4\|_{L^{r}_TL^2}&\le C2^{j\f n{p}}\!\sum_{(k,k')\in J_{\ell m}\cup J_{hh}}
\|\Delta_kf\Delta_{k'}g\|_{L^{r}_T{L^\f{2p}{2+p}}}\\
&\le Cc(j)2^{j(\f n{p}-s-\tau)}\|f\|_{\widetilde{L}_T^{r_1}\dot B^{s}_{2,1}}
\|g\|_{\widetilde{L}_T^{r_2}\dot\cB^{t,\tau}_{2,p}}.
\end{align*}

This completes the proof of Lemma \ref{Lem:Remainder}.\ef

\bthm{Proposition}\label{Prop:Product}
Let $s, t, \tilde{s},  \tilde{t},\sigma, \tau\in\mathbb{R}$, $2\le p\le 4$,
and $1\le r,r_1,r_2\le \infty$ with $\f 1 r=\f 1 {r_1}+\f 1 {r_2}$.
Then we have\vspace{0.1cm}

(a)\,\,if $\sigma, \tau\le\f np$ and $\sigma+\tau>0$, then
\begin{align}\label{equ:Product-high}
\sum_{2^j> R_0}2^{j(\sigma+\tau-\f np)}\|\Delta_j(fg)\|_{L^r_TL^p}\le C
\|f\|_{\widetilde{L}_T^{r_1}\dot\cB^{\f n2-\f np+\sigma,\sigma}_{2,p}}\|g\|_{\widetilde{L}_T^{r_2}\dot\cB^{\f n2-\f np+\tau,\tau}_{2,p}};
\end{align}

(b)\,\, if $s, \tilde{s}\le\f np,$  $s+t>n-\f {2n} p$ with $s+t=\tilde{s}+\tilde{t}$, and $\gamma\in\mathbb{R}$, then
\begin{align}\label{equ:Product-low}
&\sum_{2^j\le R_0}2^{j(s+t-\f n2)}\|\Delta_j(fg)\|_{L^r_TL^2}\nonumber\\&\le C\big(
\|f\|_{\widetilde{L}_T^{r_1}\dot\cB^{s,s-\f n2+\f n
p}_{2,p}}\|g\|_{\widetilde{L}_T^{r_2}\dot\cB^{t,t-\f n2+\f n p+\gamma}_{2,p}}+
\|g\|_{\widetilde{L}_T^{r_2}\dot\cB^{\tilde{s},\tilde{s}-\f n2+\f n
p}_{2,p}}\|f\|_{\widetilde{L}_T^{r_1}\dot\cB^{\tilde{t},\tilde{t}-\f n2+\f n p}_{2,p}}\big).
\end{align}

(c)\,\,if $s, \tilde{s}\le\f n2,$  $s+t>\f n2-\f {n} p$ with $s+t=\tilde{s}+\tilde{t}$, then
\begin{align}\label{equ:Product-L2}
&\sum_{j\in\Z}2^{j(s+t-\f n2)}\|\Delta_j(fg)\|_{L^r_TL^2}\nonumber\\&\le C\big(
\|f\|_{\widetilde{L}_T^{r_1}\dot\cB^{s,s-\f n2+\f n
p}_{2,p}}\|g\|_{\widetilde{L}_T^{r_2}\dot B^{t}_{2,1}}+
\|g\|_{\widetilde{L}_T^{r_2}\dot\cB^{\tilde{s},\tilde{s}-\f n2+\f n
p}_{2,p}}\|f\|_{\widetilde{L}_T^{r_1}\dot B^{\tilde{t}}_{2,1}}\big).
\end{align}

\ethm

\noindent{\bf Proof.}\,\,Thanks to \eqref{equ:para-high}, we get for $2^j>R_0$,
\begin{align*}\label{}
\|\Delta_j(T_fg)\|_{L^r_TL^p}+\|\Delta_j(T_gf)\|_{L^r_TL^p}
\le Cc(j)2^{j(\f np-\sigma-\tau)}\|f\|_{\widetilde{L}_T^{r_1}\dot\cB^{\f n2-\f np+\sigma,\sigma}_{2,p}}
\|g\|_{\widetilde{L}_T^{r_2}\dot\cB^{\f n2-\f np+\tau,\tau}_{2,p}},
\end{align*}
and from  \eqref{equ:Remainder-high}, we infer that
\begin{align*}\label{}
\|\Delta_j(R(f,g))\|_{L^r_TL^p}\le Cc(j)2^{j(\f np-\sigma-\tau)}\|f\|_{\widetilde{L}_T^{r_1}\dot\cB^{\f n2-\f np+\sigma,\sigma}_{2,p}}
\|g\|_{\widetilde{L}_T^{r_2}\dot\cB^{\f n2-\f np+\tau,\tau}_{2,p}},
\end{align*}
from which and Bony' decomposition \eqref{Bonydecom}, we obtain \eqref{equ:Product-high}.

From \eqref{equ:para-low}, we get for $2^j\le R_0$,
\begin{equation*}\label{}\begin{split}
&\|\Delta_j(T_fg)\|_{L^r_TL^2}\le Cc(j)2^{j(\f n2-s-t)}
\|f\|_{\widetilde{L}_T^{r_2}\dot\cB^{s,s-\f n2+\f n
p}_{2,p}}\|g\|_{\widetilde{L}_T^{r_2}\dot\cB^{t,t-\f n2+\f n p+\gamma}_{2,p}},\\
&\|\Delta_j(T_gf)\|_{L^r_TL^2}\le Cc(j)2^{j(\f n2-s-t)}\|g\|_{\widetilde{L}_T^{r_2}\dot\cB^{\tilde{s},\tilde{s}-\f n2+\f n
p}_{2,p}}\|f\|_{\widetilde{L}_T^{r_1}\dot\cB^{\tilde{t},\tilde{t}-\f n2+\f n p}_{2,p}},
\end{split}\end{equation*}
and from \eqref{equ:Remainder-low},
\begin{align*}
\|\Delta_j(R(f,g))\|_{L^r_TL^2}\le Cc(j)2^{j(\f n2-s-t)}\|f\|_{\widetilde{L}_T^{r_1}\dot\cB^{\tilde{t},\tilde{t}-\f n2+\f n
p}_{2,p}}\|g\|_{\widetilde{L}_T^{r_2}\dot\cB^{\tilde{s},\tilde{s}-\f n2+\f n
p}_{2,p}},
\end{align*}
which imply \eqref{equ:Product-low}. In the same manner, the inequality (\ref{equ:Product-L2})
can be deduced from (\ref{equ:para-high-L2}) and (\ref{equ:Remainder-high-L2}).\ef

\bthm{Proposition}\label{Prop:Commutator} Let $2\le p\le 4$,
$-\frac{n}{p}< s\le\frac{n}{2}+1$, and
$-\frac{n}{p}<\sigma\le\frac{n}{p}+1$, and $1\le
r,r_1,r_2\le \infty$ with $\f 1 r=\f 1 {r_1}+\f 1
{r_2}$. Then for $2^j\ge R_0$, there holds
\begin{align}\label{equ:Commutator-high}
&\|[v,\Delta_j]\cdot \na f\|_{L^r_TL^p}\le Cc(j)\big(2^{-j\sigma}
+2^{j(\f n2-\f np-s)}\big)\|v\|_{\widetilde{L}_T^{r_1}\dot\cB^{\f
n2+1,\f np+1}_{2,p}}
\|f\|_{\widetilde{L}_T^{r_2}\dot\cB^{s,\sigma}_{2,p}}.
\end{align}
Moreover, if $-\frac{n}{p}< s\le\frac{n}{p}+1$, then
\begin{align}
&\|[v,\Delta_j]\cdot \na f\|_{L^r_TL^2}\le
Cc(j)2^{-js}\|v\|_{\widetilde{L}_T^{r_1}\dot\cB^{\f n2+1,\f
np+1}_{2,p}} \|f\|_{\widetilde{L}_T^{r_2}\dot
B^{s}_{2,1}}.\label{equ:Commutator-L2}
\end{align}
\ethm

\noindent{\bf Proof.}\,\,Using the Bony's decomposition (\ref{Bonydecom}), we write
\begin{align*}
[v,\Delta_j]\cdot \na f=&[T_{v^i},\Delta_j]\p_if+T_{\p_i\Delta_jf}v^i+R(v^i,\p_i\Delta_jf)\\
&-\Delta_j(T_{\p_if}v^i)-\Delta_jR(v^i,\p_if).
\end{align*}
We get by using (\ref{equ:para-high}), $s\le\f n2+1$ and $\sigma\le\f np+1$ that
\begin{align*}
\|\Delta_j(T_{\p_if}v^i)\|_{L^r_TL^p}\le Cc(j)
(2^{-j\sigma}+2^{j(\f n2-\f np-s)})\|\na f\|_{\widetilde{L}_T^{r_2}\dot\cB^{s-1,\sigma-1}_{2,p}}
\|v\|_{\widetilde{L}_T^{r_1}\dot\cB^{\f n2+1,\f np+1}_{2,p}}.
\end{align*}
And from (\ref{equ:Remainder-high}), $s>-\f np$ and $\sigma>-\f n{p}$, we infer that
\begin{align*}
\|\Delta_jR(v^i,\p_if)\|_{L^r_TL^p}\le Cc(j)
(2^{-j\sigma}+2^{j(\f n2-\f np-s)})\|\na f\|_{\widetilde{L}_T^{r_2}\dot\cB^{s-1,\sigma-1}_{2,p}}
\|v\|_{\widetilde{L}_T^{r_1}\dot\cB^{\f n2+1,\f np+1}_{2,p}}.
\end{align*}
Noticing that
\begin{align*}
T_{\p_i\Delta_jf}'v^i\eqdefa
T_{\p_i\Delta_jf}v^i+R(v^i,\p_i\Delta_jf)=\sum_{k\ge
j-2}S_{k+2}\Delta_j\p_if\Delta_{k}v^i,
\end{align*}
then we get by Lemma \ref{Lem:Bernstein} that
\begin{align}\label{equ:para-lowhigh}
&\|T_{\p_i\Delta_jf}'v^i\|_{L^r_TL^p}\le
C\|\Delta_j\na f\|_{L^{r_2}_TL^\infty}\sum_{k\ge j-2}\|\Delta_{k}v\|_{L^{r_1}_TL^p}\nonumber\\
&\le C2^{j(1+ \frac{n}{p})}
\|\Delta_jf\|_{L^{r_2}_TL^p}\Big(\!\!\sum_{k\ge j-2,2^k>R_0}\!
\|\Delta_{k}v\|_{L^{r_1}_TL^p}+\!\!\!\sum_{k\ge j-2,\f {R_0} {4}\le 2^k\le R_0}\!\!\!2^{k(\f n2-\f np)}
\|\Delta_{k}v\|_{L^{r_1}_TL^2}\Big)
\nonumber\\
&\le Cc(j)2^{-j\sigma}\|v\|_{\widetilde{L}_T^{r_1}\dot \cB^{ \frac{n}{2}+1,\frac{n}{p}+1}_{2,p}}
\|f\|_{\widetilde{L}_T^{r_2}\dot
\cB^{s,\sigma}_{2,p}}.
\end{align}

Now we turn to estimate $\displaystyle[T_{v^i},
\Delta_j]\p_i f=\sum_{|k-j|\le4}[S_{k-1}v^i, \Delta_j]\pa_i \Delta_{k}f$.
Set $h(x)=(\cF^{-1}\phi)(x)$, we get by integration by parts that
\begin{align*}
[T_{v^i}, \Delta_j]\pa_i f
&=\sum_{|k-j|\le4}2^{nj}\int_{\mathbb{R}^n}h(2^j(x-y))(S_{k-1}v^i(x)-S_{k-1}v^i(y))\pa_i \Delta_{k}f(y)dy\\
&=\sum_{|k-j|\le4}2^{(n+1)j}\int_{\mathbb{R}^n}\int_0^1y\cdot\na
S_{k-1}v^i(x-\tau
y)d\tau\pa_ih(2^jy)\Delta_{k}f(x-y)dy\nonumber\\&\qquad\quad+
2^{nj}\int_{\mathbb{R}^n}h(2^j(x-y))\pa_iS_{k-1}v^i(y)\Delta_{k}f(y)dy,
\end{align*}
which together with the Minkowski inequality and
$\dot B^{\f n2,\f np}_{2,p}\hookrightarrow L^\infty$ implies
\begin{align*}
&\|[T_{v^i}, \Delta_j]\pa_i f\|_{L^r_TL^p}\\ &\le C
\|\na v\|_{L^{r_1}_TL^\infty}\Big(\!\!\sum_{|k-j|\le4, 2^k>R_0}\!\!
\|\Delta_{k}f\|_{L^{r_2}_TL^p}+\!\!\sum_{|k-j|\le4, \f {R_0}{2^4}\le2^k\le R_0}\!\!\!2^{k(\f n2-\f np)}
\|\Delta_{k}f\|_{L^{r_2}_TL^2}\Big)\nonumber\\
&\le Cc(j)\big(2^{-j\sigma}+2^{j(\f n2-\f np-s)}\big)
\|v\|_{\widetilde{L}_T^{r_1}\dot \cB^{\frac{n}{2}+1, \frac{n}{p}+1}_{2,p}}\|f\|_{\widetilde{L}_T^{r_2}\dot \cB^{s,\sigma}_{2,p}}.
\end{align*}

Summing up the above estimates, we conclude the inequality (\ref{equ:Commutator-high}). On the other hand,
we have from (\ref{equ:Remainder-high-L2}),
\begin{align*}
\|\Delta_jR(v^i,\p_if)\|_{L^r_TL^2}\le Cc(j)
2^{-js}\|f\|_{\widetilde{L}_T^{r_2}\dot B^{s}_{2,1}}\|v\|_{\widetilde{L}_T^{r_1}\dot\cB^{\f n2+1,\f np+1}_{2,p}}.
\end{align*}
With a slightly modification of \eqref{equ:para-lowhigh}, we get
\begin{align*}
\|\Delta_j(T_{\p_if}v^i)\|_{L^r_TL^2}+\|T_{\p_i\Delta_jf}'v^i\|_{L^r_TL^2}\le
Cc(j)
2^{-js}\|f\|_{\widetilde{L}_T^{r_2}\dot B^{s}_{2,1}}\|v\|_{\widetilde{L}_T^{r_1}\dot\cB^{\f n2+1,\f np+1}_{2,p}},
\end{align*}
and thanks to the representation of $[T_{v^i}, \Delta_j]\pa_i f$, we have
\begin{align*}
\|[T_{v^i}, \Delta_j]\pa_i f\|_{L^r_TL^2} &\le C
\|\na v\|_{L^{r_1}_TL^\infty}\sum_{|k-j|\le4}\|\Delta_{k}f\|_{L^{r_2}_TL^2}
\nonumber\\ &\le Cc(j)
2^{-js}\|v\|_{\widetilde{L}_T^{r_1}\dot\cB^{\f n2+1,\f np+1}_{2,p}}\|f\|_{\widetilde{L}_T^{r_2}\dot B^{s}_{2,1}}.
\end{align*}
Then the inequality (\ref{equ:Commutator-L2}) follows from the above three estimates.\ef

\bthm{Proposition}\label{Prop:composition}
Let $2\le p\le 4$, $s,\sigma>0$, and $s\ge\sigma-\f n2+\f np, r\ge 1$. Assume that
$F\in W^{[s]+2,\infty}_{\textrm{loc}}\cap W^{[\sigma]+2,\infty}_{\textrm{loc}}$ with  $F(0)=0$. Then
there holds
\begin{align}
\|F(f)\|_{\widetilde{L}_T^{r}\dot\cB^{s,\sigma}_{2,p}}\le C(1+\|f\|_{\widetilde{L}^\infty_T\dot\cB^{\f np,\f np}_{2,p}})^{\max([s],[\sigma])+1}
\|f\|_{\widetilde{L}_T^{r}\dot\cB^{s,\sigma}_{2,p}}.\label{equ:F-est}
\end{align}
For any $s>0$ and $p\ge 1$, there holds
\begin{align}
\|F(f)\|_{\widetilde{L}_T^{r}\dot B^{s}_{p,1}}\le C(1+\|f\|_{{L}^\infty_TL^\infty})^{[s]+1}
\|f\|_{\widetilde{L}_T^{r}\dot B^{s}_{p,1}}.\label{equ:F-est2}
\end{align}
\ethm

\no{\bf Proof.}\,\,The inequality (\ref{equ:F-est2}) is classical, see \cite{Danchin-cpde01}.
We only present the proof of (\ref{equ:F-est}).
Decompose $F(f)$ as
\begin{align}
F(f)=\sum_{k'\in\mathbb{Z}}F(S_{k'+1}f)-F(S_{k'}f)&=\sum_{k'\in\mathbb{Z}}\Delta_{k'}f\int_0^1F'(S_{k'}f+\tau\Delta_{k'}f)d\tau\nonumber\\
&\eqdefa \sum_{k'\in\mathbb{Z}}\Delta_{k'}f\, m_{k'},\nonumber
\end{align}
where $m_{k'}\eqdefa \int_0^1F'(S_{k'}f+\tau\Delta_{k'}f)d\tau$.
We denote
$$J_{\ell}=\{k; 2^k\le R_0\},\quad J_h=\{k; 2^k>R_0\}.$$
Then we have
$$
J_{\ell}=J_{\ell\ell}\cup J_{\ell m}\cup J_{\ell h},\quad J_{h}=J_{h\ell}\cup {J}_{h m}\cup {J}_{h h}.
$$
where
\begin{align*}&J_{\ell\ell}=\big\{(k,k'); k\in J_{\ell}, k'\le k\big\},\\
&J_{\ell m}=\big\{(k,k'); k\in J_\ell, k'>k, 2^{k'}\le R_0\big\},\\
&J_{\ell h}=\big\{(k,k'); k\in J_\ell, k'>k, 2^{k'}>R_0\big\},\\
&{J}_{h\ell}=\big\{(k,k'); k\in {J}_h, k'\le k, 2^{k'}\le R_0\big\},\\
&{J}_{hm}=\big\{(k,k'); k\in {J}_h, k'\le k, 2^{k'}>R_0\big\},\\
&{J}_{hh}=\big\{(k,k'); k\in {J}_h, k'>k\big\}.
\end{align*}
By Lemma \ref{Lem:Bernstein}, we have
\begin{align*}
\|\Delta_k(\Delta_{k'}f\,m_{k'})\|_{L^2}
\le C2^{-k|\al|}\sup_{|\gamma|=|\al|}\|D^\gamma\Delta_k(\Delta_{k'}f\,m_{k'})\|_{L^2},
\end{align*}
for any $\al\in \N^n$ and for $|\gamma|\ge0$,
$$\|D^\gamma m_{k'}\|_{L^\infty}\le C2^{k'|\gamma|}(1+\|f\|_{L^\infty})^{|\gamma|},$$
which imply
\begin{align}\label{equ:F-L2est}
\|\Delta_k(\Delta_{k'}f\,m_{k'})\|_{L^2}\le C2^{(k'-k)|\al|}
\|\Delta_{k'}f\|_{L^2}(1+\|f\|_{L^\infty})^{|\al|}.
\end{align}

We apply (\ref{equ:F-L2est}) with $|\al|=[s]+1$ to get
\begin{align}\label{equ:F-lowlow}
&\sum_{(k,k')\in J_{\ell\ell}}2^{ks}\|\Delta_k(\Delta_{k'}f\, m_{k'})\|_{L^r_TL^2}\nonumber\\
&\le C\sum_{2^{k'}\le R_0}2^{k's}\|\Delta_{k'}f\|_{L^2}
\sum_{k\ge k'}2^{(k-k')(s-[s]-1)}(1+\|f\|_{L^\infty_TL^\infty})^{|\al|}\nonumber\\
&\le C(1+\|f\|_{L^\infty})^{[s]+1}\|f\|_{\widetilde{L}_T^{r}\dot\cB^{s,\sigma}_{2,p}}.
\end{align}
And we apply (\ref{equ:F-L2est}) with $|\al|=0$ to get
\begin{align}\label{equ:F-lowmed}
\sum_{(k,k')\in J_{\ell m}}2^{ks}\|\Delta_k(\Delta_{k'}f\, m_{k'})\|_{L^r_TL^2}
&\le C\sum_{2^{k'}\le R_0}2^{k's}\|\Delta_{k'}f\|_{L^r_TL^2}
\sum_{k<k'}2^{(k-k')s}
\nonumber\\&\le C\|f\|_{\widetilde{L}_T^{r}\dot\cB^{s,\sigma}_{2,p}}.
\end{align}

Setting $m_0\eqdefa F'(0) $, we write
$$
\Delta_k(\Delta_{k'}f\, m_{k'})=\Delta_k(\Delta_{k'}f\, m_{0})+\Delta_k(\Delta_{k'}f (m_{k'}-m_0)).
$$
It is easy to find that
\begin{align*}
\sum_{(k,k')\in J_{\ell h}}2^{ks}\|\Delta_k(\Delta_{k'}f\, m_0)\|_{L^r_TL^2}
\le C\|f\|_{\widetilde{L}_T^{r}\dot\cB^{s,\sigma}_{2,p}},
\end{align*}
and using the formula
\begin{align*}
m_{k'}-m_0&=\int_0^1m_{k'}'(\tau(S_{k'-1}f+\Delta_{k'}f))d\tau(S_{k'-1}f+\Delta_{k'}f)\\
&\eqdefa (S_{k'-1}f+\Delta_{k'}f)\widetilde{m}_{k'},
\end{align*}
we obtain
\begin{align*}
&\sum_{(k,k')\in J_{\ell h}}2^{ks}\|\Delta_k(\Delta_{k'}f (m_{k'}-m_0))\|_{L^r_TL^2}\\
&\le C\sum_{(k,k')\in J_{\ell h}}2^{ks}\bigl(\|\Delta_k(\Delta_{k'}f\widetilde{m}_{k'} S_{k'-1}f)\|_{L^r_TL^2}
+2^{k(\f {2n}{p}-\f n2)}
\|(\Delta_{k'}f)^2\|_{L^r_TL^\f p2}\bigr)\\
&\eqdefa I_1+I_2.
\end{align*}
Now we have
\begin{align*}
I_1 \le& C\sum_{(k,k')\in J_{\ell h}}2^{ks}
\|\Delta_{k'}f\|_{L^r_TL^p}\sum_{k''\le k'-2, 2^{k''}\le R_0}\|\Delta_{k''}f\|_{L^\infty_TL^\f{2p}{p-2}}\nonumber\\&
+C\sum_{(k,k')\in J_{\ell h}}2^{ks}2^{k(\f {2n}{p}-\f n2)}\sum_{k''\le k'-2, 2^{k''}> R_0}
\|\Delta_{k'}f\Delta_{k''}f\|_{L^r_TL^\f p2}.
\end{align*}
The first term of the right hand side is bounded by
\begin{align*}
&\sum_{2^{k'}\ge R_0}2^{k'\sigma}\|\Delta_{k'}f\|_{L^r_TL^p}2^{-k'\sigma}
\sum_{2^{k}\le R_0}2^{ks}
\sum_{2^{k''}\le R_0}2^{k''\f np}\|\Delta_{k''}f\|_{L^\infty_TL^2}\nonumber\\&\le C
\|f\|_{\widetilde{L}_T^{\infty}\dot\cB^{\f np,\f np}_{2,p}}\|f\|_{\widetilde{L}_T^{r}\dot\cB^{s,\sigma}_{2,p}},
\end{align*}
and the second term is bounded by
\begin{align*}
&\sum_{J_{\ell h}}2^{k(s+\f {2n}{p}-\f n2)}
2^{k'\sigma}\|\Delta_{k'}f\|_{L^r_TL^p}2^{-k'\sigma}
\sum_{2^{k''}> R_0}2^{k''\f np}
\|\Delta_{k''}f\|_{L^\infty_TL^p}2^{-k''\f np}\nonumber\\&\le C
\|f\|_{\widetilde{L}^\infty_T\dot\cB^{\f np,\f np}_{2,p}}\|f\|_{\widetilde{L}_T^{r}\dot\cB^{s,\sigma}_{2,p}}.
\end{align*}
Similarly, we have
\begin{align*}
&I_2\le C\|f\|_{\widetilde{L}^\infty_T\dot\cB^{\f np,\f np}_{2,p}}\|f\|_{\widetilde{L}_T^{r}\dot\cB^{s,\sigma}_{2,p}}.
\end{align*}
Thus we obtain
\begin{align}\label{equ:F-lowhigh}
\sum_{(k,k')\in J_{\ell h}}2^{ks}\|\Delta_k(\Delta_{k'}f\, m_{k'})\|_{L^r_TL^2}\le C
\big(1+\|f\|_{\widetilde{L}^\infty_T\dot\cB^{\f np,\f np}_{2,p}}\big)\|f\|_{\widetilde{L}_T^{r}\dot\cB^{s,\sigma}_{2,p}}.
\end{align}

Using the argument as leading to \eqref{equ:F-lowlow} and $s\ge\sigma-\f n2+\f np$, we have
\begin{align}\label{equ:F-highlow}
&\sum_{(k,k')\in J_{h\ell}}2^{k\sigma}\|\Delta_k(\Delta_{k'}f\, m_{k'})\|_{L^r_TL^p}\nonumber\\
&\le C(1+\|f\|_{L^\infty_TL^\infty})^{[\sigma]+1}\sum_{2^{k'}\le
R_0}2^{k's}\|\Delta_{k'}f\|_{L^r_TL^2} \sum_{k\ge
{k'}}2^{(k-k')(\sigma-[\sigma]-1)}2^{k'(\f n2-\f np-s+\sigma)}
\nonumber\\& \le
C(1+\|f\|_{L^\infty_TL^\infty})^{[\sigma]+1}\|f\|_{\widetilde{L}_T^{r}\dot\cB^{s,\sigma}_{2,p}},
\end{align}
and
\begin{align}\label{equ:F-highmed}
&\sum_{(k,k')\in J_{hm}}2^{k\sigma}\|\Delta_k(\Delta_{k'}f\, m_{k'})\|_{L^r_TL^p}\nonumber\\
&\le C(1+\|f\|_{L^\infty_TL^\infty})^{[\sigma]+1}\sum_{2^{k'}>R_0}2^{k'\sigma}\|\Delta_{k'}f\|_{L^r_TL^p}
\sum_{k>{k'}}2^{(k-k')(\sigma-[\sigma]-1)}\nonumber\\
&\le C(1+\|f\|_{L^\infty_TL^\infty})^{[\sigma]+1}\|f\|_{\widetilde{L}_T^{r}\dot\cB^{s,\sigma}_{2,p}}.
\end{align}
Finally due to $\sigma>0$, we have
\begin{align}\label{equ:F-highhigh}
\sum_{(k,k')\in J_{hh}}2^{k\sigma}\|\Delta_k(\Delta_{k'}f\, m_{k'})\|_{L^r_TL^p}
&\le C\sum_{2^{k'}\ge R_0}2^{k'\sigma}\|\Delta_{k'}f\|_{L^r_TL^p}
\sum_{k'\ge k}2^{(k-k')\sigma}\nonumber\\
&\le C\|f\|_{\widetilde{L}_T^{r}\dot\cB^{s,\sigma}_{2,p}}.
\end{align}

Summing up (\ref{equ:F-lowlow})--(\ref{equ:F-highhigh}) and noting $\dot\cB^{\f np,\f np}_{2,p} \hookrightarrow L^\infty$,
we obtain (\ref{equ:F-est}).\ef

\section{Acknowledgments}
Q. Chen and C.
Miao  were partially  supported by the NSF of China (grants 10701012 and 10725102).
Z. Zhang was  supported by the NSF of China (grant 10601002).


\begin{thebibliography}{50}

\bibitem{Abidi} H. Abidi, {\it \'{E}quation de Navier-Stokes avec densit\'{e} et viscosit\'{e}
variables dans l'espace critique.} Rev. Mat. Iberoam.  23(2007), 537--586.

%\bibitem{Abidi-Hmidi} H. Abidi and T. Hmidi,
%{\it On the global well-posedness of the critical quasi-geostrophic equation}.
%SIAM J. Math. Anal., 40(2008), 167--185.

\bibitem{Bony} J.-M. Bony, {\it Calcul symbolique et propagation
des singulariti\'{e}s pour les \'{e}quations aux d\'{e}riv\'{e}es
partielles non lin\'{e}aires}. Ann. de l'Ecole Norm. Sup.,
14(1981), 209-246.

\bibitem{Cannone} M. Cannone,
{\it A generalization of a theorem by Kato on Naiver-Stokes equations}.
Revista Matem\"{a}tica Iberoamericana, 13(1997)515-541.


\bibitem{Cannon-Meyer-Planchon}
M. Cannone, Y. Meyer and F. Planchon,
{\it Solutions autosimilaires des {\'e}quations de Navier-Stokes}.
{S{\'e}minaire ``{\'E}quations aux D{\'e}riv{\'e}es Partielles" de l'{\'E}cole polytechnique},
Expos{\'e} VIII, 1993--1994.

\bibitem{Cannon-Miao-Wu} M. Cannone, C. Miao and G. Wu, {\it On the inviscid limit of the two-dimensional Navier-Stokes equations with
fractional diffusion}. Advances in Mathematical Science and Application, 18(2008),607-624.

\bibitem{Chemin-book} J.-Y. Chemin,
{\it Perfect incompressible fluids}. Oxford University Press, New York, 1998.

\bibitem{Chemin-Lecture} J.-Y. Chemin, {\it Localization in Fourier space and Navier-Stokes system.}
Phase space analysis of partial differential equations. Vol. I, 2004, 53-135.

\bibitem{Chemin-Gal}
J.-Y. Chemin and I. Gallagher, {\it On  the global wellposedness
of the 3-D Navier-Stokes equations with large initial data}.
Ann. de l'Ecole Norm. Sup., (39)2006, 679--698.

\bibitem{Chemin-Gal-Paicu}
J.-Y. Chemin, I. Gallagher and M. Paicu,
{\it Global regularity for  some classes of large solutions to the Navier-Stokes equations}. preprint, 2008.

\bibitem{Chemin-Zhang} J.-Y. Chemin and P. Zhang,
{\it On the global  wellposedness of the 3-D incompressible anisotropic
Navier-Stokes equations}. Comm. Math. Phys., (272)2007, 529--566.

\bibitem{Danchin-inven} R. Danchin,
{\it Global existence  in  critical spaces for compressible Navier-Stokes equations}.
Invent. Math., 141(2000), 579-614.

\bibitem{Danchin-cpde01} R. Danchin,
{\it Local theory in critical spaces for compressible viscous and heat-conductive gases},
Comm. Partial Differential Equations, 26(2001), 1183-1233.

\bibitem{Danchin-PRSE} R. Danchin,
{\it Density-dependent incompressible viscous fluids in critical spaces}.
Proc. Roy. Soc. Edinburgh Sect. A, 133(2003), 1311-1334.

\bibitem{Danchin-NDEA} R. Danchin,
{\it On the uniqueness  in critical spaces for compressible Navier-Stokes equations}.
Nonlinear Differential Equations Appl., 12(2005), 111-128.

\bibitem{Danchin-cpde07} R. Danchin,
{\it Well-posedness in critical spaces for barotropic viscous fluids with truly not constant density}.
Comm. Partial Differential Equations, 32(2007), 1-17.

\bibitem{Danchin-JHDE}R. Danchin, {\it Uniform estimates for Transport-Diffusion equations}.
J. Hyper. Differential Equations, 32(2007), 1373-1397.

\bibitem{Fuj-Kat} H. Fujita and T. Kato, {\it On the Navier-Stokes initial value problem I}.
Arch. Rational Mech. Anal., 16(1964), 269-315.

\bibitem{Hmidi-Keraani-Advance} T. Hmidi and S. Keraani, {\it Global solutions of the super-critical 2D Q-G equation in Besov spaces},
Advances in Math., 214(2007), 618-638.

\bibitem{Hmidi-Keraani-ARMA} T. Hmidi and S. Keraani, {\it Incompressible viscous flows in borderline Besov space.}
Arch. Rational Mech. Anal., 189(2008), 283--300.

\bibitem{Hoff-Zum} D. Hoff and K. Zumbrun, {\it Multi-dimensional diffusion waves for the Navier-Stokes equations of compressible flow.}
Indiana Univ. Math. Jour., 44(1995), 603-676.

\bibitem{Kato} T. Kato,
{\it Strong $L^p$-solutions of the Navier-Stokes equation in
$\R^m$ with applications to weak solutions}. Mathematische Zeitschrift, (187)1984, 471-480.

\bibitem{Lions} P.-L. Lions, {\it Mathematical Topics in Fluid Mechanics}. Vol.2, Compressible models,
Oxford University Press, 1998.

\bibitem{Mat-Nish}  A. Matsumura and T. Nishida, {\it The initial value problem for the equations of motion of compressible viscous and heat-conductive fluids}.
Proc. Japan Acad. Ser. A Math. Sci. 55(1979), 337-342.

\bibitem{Nash} J. Nash, {\it Le probl\`{e}me de Cauchy pour les \'{e}quations diff\'{e}rentielles d'un fluide g\'{e}n\'{e}ral},
Bulletin de la Soc. Math. de France, 90(1962),487-497.

\bibitem{Paicu-Zhang} M. Paicu and Z. Zhang,
{\it Global Regularity for the Navier-Stokes equations with large, slowly varying initial data in the vertical direction},
preprint, 2009.

\bibitem{Weissler} F. Weissler, {\it The Navier-Stokes Initial Value Problem in~$L^p$}.
Arch. Rational Mech. Anal., (74)1980, 219-230.


\end{thebibliography}
\end{document}